\pdfoutput=1
\documentclass[leqno,12pt]{article}

\usepackage{amsmath,amssymb,rotating,a4wide,color}
\usepackage[integrals]{wasysym}
\usepackage{enumerate}
\usepackage{mathtools}
\usepackage[font=footnotesize,labelfont=bf,margin=2cm]{caption}

\usepackage[dvipsnames]{xcolor}
\usepackage{bbm} 
\usepackage[mathscr]{eucal} 

\usepackage{xfrac} 

\usepackage{tocloft}
\cftsetindents{section}{0em}{2em}

\usepackage{array}
\newcolumntype{M}{>{\centering\arraybackslash}m{\dimexpr.25\linewidth-2\tabcolsep}}

\usepackage{tikz}
\usetikzlibrary{circuits.logic.US,circuits.logic.IEC,fit}
\usetikzlibrary{cd}
\usetikzlibrary{backgrounds}
\usetikzlibrary{matrix,snakes}
\usepackage{subcaption}
\usepackage{lscape}

\usepackage[all,cmtip]{xy}

\usepackage{verbatim}

\usepackage{bbding}             
\usepackage{marginnote}         

\usepackage{hyperref}           

\newdir^{ (}{{}*!/-3pt/\dir^{(}}
\newdir_{ (}{{}*!/-3pt/\dir_{(}}

\entrymodifiers={+!!<0pt,\fontdimen22\textfont2>}

\def\bigtimes{\mathop{\raise-2pt\hbox{\huge$\times$}}}

\newbox\circbulletbox
\setbox\circbulletbox=\hbox{\,{\large$\circledcirc$}\kern-8.7pt\raise .6pt\hbox{$\bullet$}\,}

\let\le\leqslant
\let\ge\geqslant

\def\LIC#1{{\CI\CC}_{\ICF,#1}}
\def\GLIC#1{{\CI\CC}_{F,#1}}

\def\circVbig{\hbox{\text{\it\r{V}}}}
\def\circVscript{\hbox{\scriptsize\text{\it\r{V}}}}
\def\circVscriptscript{\mbox{\tiny\text{\it\r{V}}}}

\def\circVlimits_#1^#2{{\mathchoice%
{\circVbig{}^{\kern2pt #2}_{\kern-2pt #1}}%
{\circVbig{}^{\kern2pt #2}_{\kern-2pt #1}}%
{\scriptstyle\circVscript{}^{\kern1.7pt #2}_{\kern-1pt #1}}%
{\scriptscriptstyle\circVscriptscript{}^{\kern1.5pt #2}_{\kern-1pt #1}}%
}}
\def\circVr_#1{\circVlimits_#1^r}
\def\circVs_#1{\circVlimits_#1^s}

\def\CMod{\mathscr{M}}
\def\minusoneoverh{{\kern-1pt\sfrac{-\kern-1.5pt 1}{\kern1pt h}}}
\def\oneoverh{{\sfrac{1}{\kern1.5pt h}}}
\def\minusoneoverr{{\kern-1pt\sfrac{-\kern-1.5pt 1}{\kern.8pt r}}}
\def\oneoverr{{\sfrac{1}{\kern.8pt r}}}

\def\circWbig{\hbox{\text{\it\r{W}}}}
\def\circWscript{\hbox{\scriptsize\text{\it\r{W}}}}
\def\circWscriptscript{\mbox{\tiny\text{\it\r{W}}}}

\def\circWlimits_#1^#2{{\mathchoice%
{\circWbig{}^{\kern2pt #2}_{\kern-2pt #1}}%
{\circWbig{}^{\kern2pt #2}_{\kern-2pt #1}}%
{\scriptstyle\circWscript{}^{\kern1.7pt #2}_{\kern-1pt #1}}%
{\scriptscriptstyle\circWscriptscript{}^{\kern1.5pt #2}_{\kern-1pt #1}}%
}}

\newcommand{\Rsep}{R^{\mkern+2mu\sep}}%
\newcommand{\Rxhens}{R_x^{\kern+.7pt\hens}}%
\newcommand{\Rxshens}{R_x^{\kern+1.5pt\shens}}%
\newcommand{\kxsep}{k_x^{\kern+1pt\sep}}%

\def\OM{\mathchoice
{\rlap{\kern3.2pt$\overline{\phantom{L}}$}M}
{\rlap{\kern3.2pt$\overline{\phantom{L}}$}M}
{\rlap{\kern2.4pt$\scriptstyle\overline{\phantom{L}}$}M}
{\rlap{\kern1.8pt$\scriptscriptstyle\overline{\phantom{L}}$}M}}

\def\mycirc{{\kern1pt\circ\kern2pt}}
\def\charact{\mathop{\rm char}\nolimits}

\def\Mdirtext{\kern6pt\widetilde{\phantom{t}}\kern-10pt M}%
\def\Mdirscript{\kern0pt\widetilde{\phantom{N}}\kern-9pt M}%
\def\Mdirscriptscript{\kern.5pt\widetilde{\phantom{N}}\kern-7.5pt M}%

\def\subsetdown{{\rotatebox[origin=c]{-90}{$\textstyle\subset$}}}
\def\subsetup{{\rotatebox[origin=c]{90}{$\textstyle\subset$}}}
\def\equaldown{{\rotatebox[origin=c]{90}{$\textstyle=$}}}

\def\opp{{\rm opp}}

\def\kalg{{\itoverline{k}}}
\def\Kalg{{\itoverline[4]{K}}}
\def\Ksep{{K^{\rm sep}}}

\def\eval{\mathop{\rm eval}\nolimits}
\def\pptauone{(\kern-1pt(\tau^{-1})\kern-1pt)}
\def\ev{\mathop{\rm ev}\nolimits}
\def\coev{\mathop{\rm coev}\nolimits}

\def\gr{{\rm gr}}

\def\Aut{\mathop{\rm Aut}\nolimits}
\def\Frob{\mathop{\rm Frob}\nolimits}
\def\Hom{\mathop{\rm Hom}\nolimits}

\def\Gal{\mathop{\rm Gal}\nolimits}
\def\End{\mathop{\rm End}\nolimits}

\def\Spec{\mathop{\rm Spec}\nolimits}

\def\deg{\mathop{\rm deg}\nolimits}

\def\coker{\mathop{\rm coker}\nolimits}

\def\Mat{\mathop{\rm Mat}\nolimits}

\def\Norm{\mathop{\rm Norm}\nolimits}

\def\Quot{\mathop{\rm Quot}\nolimits}

\def\trace{\mathop{\rm tr}\nolimits}

\def\ord{\mathop{\rm ord}\nolimits}
\def\rank{\mathop{\rm rank}\nolimits}
\def\proj{\mathop{\rm pr}\nolimits}
\def\GL{\mathop{\rm GL}\nolimits}

\def\ad{{\rm ad}}

\def\geom{{\rm geom}}

\def\sep{{\rm sep}}
\def\shens{{\rm sh}}
\def\hens{{\rm h}}

\def\id{{\rm id}}

\let\phi\varphi
\let\epsilon\varepsilon
\let\setminus\smallsetminus
\let\oldnmid\nmid
\def\nmid{\kern-1pt\oldnmid\kern-1pt}
\def\ndiv\nmid
\let\emptyset\varnothing
\def\barbar#1{\bar{\bar#1}}

\newcommand{\BA}{{\mathbb{A}}}

\newcommand{\BF}{{\mathbb{F}}}

\newcommand{\BQ}{{\mathbb{Q}}}
\newcommand{\BR}{{\mathbb{R}}}

\newcommand{\BZ}{{\mathbb{Z}}}

\newcommand{\Fa}{{\mathfrak{a}}}

\newcommand{\Fm}{{\mathfrak{m}}}

\newcommand{\Fp}{{\mathfrak{p}}}
\newcommand{\Fq}{{\mathfrak{q}}}

\newcommand{\FP}{{\mathfrak{P}}}

\newcommand{\CC}{{\cal C}}

\newcommand{\CE}{{\cal E}}
\newcommand{\CF}{{\cal F}}

\newcommand{\CI}{{\cal I}}

\newcommand{\CL}{{\cal L}}
\newcommand{\CM}{{\cal M}}
\newcommand{\CN}{{\cal N}}
\newcommand{\CO}{{\cal O}}
\newcommand{\CP}{{\cal P}}

\newbox\mybox
\def\arrover#1{\mathrel{
\setbox\mybox=\hbox spread 1.4em
{\hfil\raisebox{1pt}{$\scriptstyle#1$}\hfil}
\vbox{\offinterlineskip\copy\mybox
\hbox to\wd\mybox{\rightarrowfill}}}}

\def\larrover#1{\mathrel{
\setbox\mybox=\hbox spread 1.4em
{\hfil$\scriptstyle#1\vphantom{g}$\hfil}
\vbox{\offinterlineskip\copy\mybox
\hbox to\wd\mybox{\leftarrowfill}}}}

\def\ontoover#1{\mathrel{
\setbox\mybox=\hbox spread 1.4em
{\hfil$\scriptstyle#1\vphantom{g}$\hfil}
\vbox{\offinterlineskip\copy\mybox
\hbox to\wd\mybox{\rightarrowfill\hskip-2.8mm
$\rightarrow$}}}}
\def\leftontoover#1{\mathrel{
\setbox\mybox=\hbox spread 1.4em
{\hfil$\scriptstyle#1\vphantom{g}$\hfil}
\vbox{\offinterlineskip\copy\mybox
\hbox to\wd\mybox{$\leftarrow$\hskip-2.8mm
\leftarrowfill}}}}
\let\longto\longrightarrow
\let\into\hookrightarrow
\let\onto\twoheadrightarrow
\def\longonto{\ontoover{\ }}
\def\longinto{\lhook\joinrel\longrightarrow}

\def\isoto{\mathrel{
\setbox\mybox=\hbox spread 0.9em
{\hfil$\scriptstyle\sim$\hfil}
\vbox{\offinterlineskip\copy\mybox
\hbox to\wd\mybox{\rightarrowfill}}}}

\def\dirlim{\mathop{\vtop{\hbox{\rm lim}\vskip-8pt
\hbox{\hskip1pt$\scriptstyle\longrightarrow$}\vskip-1pt}}}

\def\Bigskip{\bigskip\bigskip}

\newtheorem{Thm}{Theorem}[section]
\newtheorem{Prop}[Thm]{Proposition}

\newtheorem{Lem}[Thm]{Lemma}

\newtheorem{Cor}[Thm]{Corollary}

\newtheorem{Def}[Thm]{Definition}

\newtheorem{PropDef}[Thm]{Proposition-Definition}

\newtheorem{Rem}[Thm]{Remark}
\newtheorem{Ex}[Thm]{Example}

\newtheorem{Cons}[Thm]{Construction}
\newtheorem{PropCons}[Thm]{Proposition-Construction}

\def\UseTheoremCounterForNextEquation{\setcounter{equation}{\value{Thm}}\addtocounter{Thm}{1}}
\counterwithin{equation}{section}

\def\qed{{\hskip0pt\unskip\unskip\nobreak\hfil\penalty50
\hskip1em\hbox{}\nobreak\hfil
{$\square$}
\parfillskip=0pt\finalhyphendemerits=0
\par}\medskip}

\newenvironment{Proof}
{\noindent{\bf Proof.}}
{\qed}

{\noindent{\bf Proof of #1.}\ }
{\qed}

\newcommand{\StatementTagList}{{\upshape(\it\alph{enumi}\kern1pt\upshape)}}
\newcommand{\StatementTagText}{{(\alph{enumi})}}

\newcommand{\StatementLabels}{%
  \renewcommand{\labelenumi}{\StatementTagList}%
  \renewcommand{\theenumi}{\StatementTagText}%
}

\newcommand{\isochar}{\ensuremath{\sim}}
\newcommand{\genisosign}[1]{\smash{\raisebox{-0.65ex}{#1}}}
\newcommand{\isosign}{\genisosign{\isochar}} 

\newcommand{\longisospace}{\hspace{.375ex}}
\newcommand{\longisoarrow}{\xrightarrow{\longisospace\isosign\longisospace}}
\newcommand{\longisoto}{\longisoarrow}

\reversemarginpar                       

\newcommand{\coloredit}{RoyalBlue}      

\newcommand{\editsym}[1]{\smash{\scalebox{1.2}{#1}}}
\newcommand{\editpoint}[1]{\editsym{\raisebox{-.25em}{#1}}}
\newcommand{\editpointright}{\editpoint{{\rm\HandRight}}}

  {\par\noindent\leavevmode\begingroup\color{\coloredit}\marginnote{\editpointright}\ignorespaces}%
  {\par\noindent\endgroup}

\newcommand{\itoverline}[2][3]{{}\mkern#1mu\overline{\mkern-#1mu#2}}

\newcommand{\tholine}[1]{
  \vbox{%
    \hrule height .06em
    \kern.25ex%
    \hbox{%
      \kern-.15ex%
      \ensuremath{#1}%
      \kern-.15ex%
    }
  }
}

\newcommand{\vartholine}[1]{
  \vbox{%
    \hrule height .06em
    \kern.25ex%
    \hbox{%
      \kern-.1ex%
      \ensuremath{#1}%
    }
  }
}

\newcommand{\wedgesym}{\mbox{\large$\wedge$}}
\newcommand{\wedgesymscript}{\scalebox{0.78}{$\wedge$}}
\renewcommand{\bigwedge}{\mathchoice{\wedgesym\mkern-1mu}{\wedgesym\mkern-1mu}{\wedgesymscript}{\wedgesymscript}}

\newcommand{\GK}{\varGamma_{\kern-1ptK}}        
\newcommand{\GKn}{\varGamma_{\kern-1ptn}}       
\newcommand{\Gk}{\varGamma_{\kern-1ptk}}        
\newcommand{\Gkx}{\varGamma_{\kern-1ptk_x}}     

\newcommand{\GKgeom}{\varGamma_{\kern-1pt K}^{\kern+1pt\geom}} 

\newcommand{\Wof}[1]{W_{\kern-1pt#1}}              
\newcommand{\WK}{W_{\kern-1pt K}}               
\newcommand{\Wk}{W_{\kern-1pt k}}               

\newcommand{\hatBZ}{\widehat{\BZ}}

\def\Mtildetext{\kern5.8pt\widetilde{\phantom{t}}\kern-10pt M}
\def\Mtildescript{\kern1.4pt\widetilde{\phantom{N}}\kern-9pt M}
\def\Mtildescriptscript{\kern1.1pt\widetilde{\phantom{N}}\kern-7.5pt M}

\newcommand{\Mpbar}{\itoverline[4]{M}_{\kern-1pt\Fp}} 

\newcommand{\complot}{\mathbin{\widehat{\vbox to 6.4pt{}\smash{\otimes}}}} 

\newcommand{\centrot}{\mathop{{\otimes}}\limits} 
\newcommand{\otend}[2][-3]{\kern-4pt\centrot_{\raisebox{-4pt}{$\scriptstyle\End(#2)$}}\kern#1pt#2} 

\newcommand{\Tuad}{T^{\kern.5pt\circ}{\hspace{-1.3ex}}_{\ad}}

\newcommand{\Tux}[1]{T^{\kern.5pt\circ}{\hspace{-1.2ex}}_{#1\,}}

\newcommand{\BFq}{\BF_{\mkern-2muq}}      
\def\BFqn{{\BF_{\mkern-2mu q^n\!}}}

\newcommand{\Bado}{{B^{\kern.5pt\circ}{\hspace{-1ex}}_{\ad}}}

\newcommand{\Rado}{{R^{\kern1pt\circ}{\hspace{-.9ex}}_{\ad}}}

\newcommand{\Bop}{{B^{\kern.5pt\circ}{\hspace{-1ex}}_{\Fp}}}

\newcommand{\KOKsym}[1]{\CP_{\kern-1pt#1}}

\newcommand{\Wn}[2]{W_{\kern-1pt#1}\hspace{1pt}#2}%
\newcommand{\taulin}{\tau^{\textup{lin}}}

\newcommand{\unit}{\mathbbm{1}} 
\newcommand{\IHom}{\mathop{\rm hom}\nolimits}

\newcommand{\ICF}{F_{\kern-.33pt\Fp}}

\newcommand{\ICO}{{A_\Fp}}

\newcommand{\ICk}{{k_\Fp}}
\newcommand{\ICd}{d_\Fp}
\newcommand{\ICu}{z_\Fp}

\newcommand{\ICFalt}{F_{\kern-.33pt\Fp}{\kern-1.64pt'}}

\newcommand{\ppX}[1]{(\kern-1.33pt(#1)\kern-1.33pt)}
\newcommand{\bbX}[1]{[\kern-.24pt[#1]\kern-.24pt]}
\def\ppu{\ppX{\ICu}}
\def\bbu{\bbX{\ICu}}

\newcommand{\unrsup}{\raisebox{1.24pt}{$\scriptstyle\rm nr$}}
\newcommand{\ICFnr}{\ICF{}^{\kern-1.5pt\unrsup}}
\newcommand{\ICFnrcpl}{\widehat{F}_{\kern-.33pt\Fp}{}^{\kern-1pt\unrsup}}
\newcommand{\ICFalg}{\itoverline{F}_{\kern-1.72pt\Fp}}

\newcommand{\EFp}[1]{\CE_{\ICF,#1}}
\newcommand{\EF}[1]{\CE_{F,#1}}
\newcommand{\EAp}[1]{\CE_{\ICO,#1}}
\newcommand{\EFpK}{\EFp{K}}
\newcommand{\EFpk}{\EFp{k}}
\newcommand{\EFK}{\EF{K}}
\newcommand{\EFk}{\EF{k}}

\newcommand{\SAa}[2]{S_{#1,#2}}
\newcommand{\SA}[1]{\SAa{#1}{\alpha}}

\begin{document}

\title{\strut
\vskip-80pt
Weil representations associated to\\
isocrystals over function fields
}
\author{
\begin{minipage}{.3\hsize}
Maxim Mornev\\[12pt]
\small 
ICMAT\\
C. Nicol\'as Cabrera, 13--15\\
28049 Madrid\\
Spain\\
maxim.mornev@icmat.es\\[9pt]
\end{minipage}
\qquad
\begin{minipage}{.3\hsize}
Richard Pink\\[12pt]
\small Department of Mathematics \\
ETH Z\"urich\\
8092 Z\"urich\\
Switzerland \\
pink@math.ethz.ch\\[9pt]
\end{minipage}
}
\date{\today}

\maketitle

\Bigskip

\begin{abstract}\noindent
Every Anderson $A$-motive $M$ over a field determines a compatible system of Galois representations on its Tate modules at almost all primes of~$A$. This adapts easily to $F$-isocrystals, which are rational analogues of $A$-motives for the global function field $F:=\Quot(A)$. We extend this compatible system by constructing a Weil group representation associated to $M$ for \emph{every} place of~$F$. To this end we generalize the Tate module construction to a tensor functor on $F_\Fp$-isocrystals that are not necessarily pure. To prove that this yields a compatible system, we work out how that construction behaves under reduction of~$M$. 
As an offshoot we obtain a new kind of $\wp$-adic Weil representations associated to Drinfeld modules of special characteristic~$\wp$.
\end{abstract}

{\renewcommand{\thefootnote}{}
\footnotetext{MSC classification: 11G09, 11F80}
}

\newpage
\tableofcontents

\newpage
\section{Introduction}
\label{Intro}

Let $F$ be the function field of a smooth connected projective algebraic curve $C$ over a finite field $\BFq$ with $q$ elements. Let $A\subset F$ be the subring of all functions that are regular outside a fixed closed point $\infty\in C$. Let $K$ be another field over~$\BFq$, and let $M$ be an Anderson $A$-motive of rank $r$ over~$K$, for instance that associated to a Drinfeld module of rank $r$ over~$K$. Let $c^*\colon A\to K$ be the characteristic homomorphism of~$M$ and let $\wp\subset A$ denote its kernel. Then for every place $\Fp\not=\wp,\infty$ of~$F$, the rational $\Fp$-adic Tate module $V_\Fp(M)$ is an $F_\Fp$-vector space of dimension $r$ endowed with a natural left action of the Galois group $\GK := \Gal(\Ksep/K)$. 

If $K$ is finitely generated over~$\BFq$, it is known that these Tate modules form a compatible system of Galois representations. Namely, take a normal subring $R\subset K$ that is finitely generated over~$\BFq$, such that $M$ arises from an $A$-motive over~$R$. Then $c^*$ corresponds to a morphism of schemes $c\colon \Spec R\to C$. For every $\Fp\not=\wp,\infty$ and every closed point $x\in\Spec R$ with $c(x)\not=\Fp$, 
the inertia group at $x$ acts trivially on $V_\Fp(M)$ and the characteristic polynomial of any geometric Frobenius element at $x$ has coefficients in~$A$ and is independent~of~$\Fp$.
(See Gardeyn \cite[Prop.\,3.3]{GardeynT} for the case when $\wp = (0)$.)

When $M$ arises from a Drinfeld module, J.-K.\;Yu \cite{Yu2003} and the first author \cite[\S9.2]{MornevT} have extended this compatible system to $\Fp=\infty$ by a representation of the Weil group $\WK\subset\GK$ with values in the group of units of a central division algebra of dimension $r^2$ over~$F_\infty$. In the present article we construct a corresponding extension for every $A$-motive and likewise for $\Fp=\wp$ if $\wp$ is a maximal ideal.

\medskip
In fact, we construct a compatible system of representations of $\WK$ for any system of $\sigma$-linear equations with non-zero determinant. For this consider the total ring of quotients 
\UseTheoremCounterForNextEquation
\begin{equation}\label{IntroEFKDef}
\EFK := \Quot(F \otimes_{\BFq} K),
\end{equation}
endowed with the partial Frobenius endomorphism induced by $\sigma\colon a\otimes\xi \mapsto a\otimes\xi^q$. Then by an \emph{$F$-isocrystal over $K$ (with the ground field~$\BFq$)} we mean a finitely generated projective $\EFK$-module $M$ together with an isomorphism
$$\taulin_M\colon\ \sigma^* M := \EFK\otimes_{\EFK,\sigma}M\ \longisoto\ M.$$
With a natural tensor product these objects form an $F$-linear tannakian category.

\medskip
Similarly, for every place $\Fp$ of $F$ we consider the ring
\UseTheoremCounterForNextEquation
\begin{equation}\label{IntroEFpKDef}
\EFpK := \ICF\complot_{\BFq}K,
\end{equation}
endowed with the same partial Frobenius endomorphism~$\sigma$. By an \emph{$\ICF$-isocrystal over $K$ (with the ground field~$\BFq$)} we mean a finitely generated projective $\EFpK$-module $M_\Fp$ together with an isomorphism
$$\taulin_{M_\Fp}\colon\ \sigma^* M_\Fp := \EFpK\otimes_{\EFpK,\sigma}M_\Fp\ \longisoto\ M_\Fp.$$
Again, with a natural tensor product these objects form an $\ICF$-linear tannakian category.

\medskip
Moreover, the canonical embedding $F \otimes_{\BFq} K \into \ICF\complot_{\BFq}K$ extends uniquely to a $\sigma$-equivariant homomorphism $\EFK \to \EFpK$. By base change, any $F$-isocrystal $M$ over~$K$ thus yields an $\ICF$-isocrystal $M_\Fp$ over~$K$, and this construction is functorial and commutes with the tensor product. To obtain the desired compatible system associated to~$M$, we thus construct a natural representation of $\WK$ for every $\ICF$-isocrystal $M_\Fp$ over~$K$.

\medskip
For this we use a notion of \emph{pure $\ICF$-isocrystals of slope~$\alpha\in\BQ$}, and show that $M_\Fp$ possesses a natural filtration with pure subquotients $M_{\Fp,\alpha}$ of slope~$\alpha$ for all $\alpha\in\BQ$. 
Let $\kalg$ be the algebraic closure of $\BFq$ within a fixed algebraic closure $\Kalg$ of~$K$. Then for every $\alpha\in\BQ$ there exists a simple $\ICF$-isocrystal $\smash{\SA{\kalg}}$ over $\kalg$ that is pure of slope~$\alpha$, and this isocrystal is unique up to isomorphism. 
Moreover its endomorphism ring $D_\alpha$ is a central division algebra over~$\ICF$. 
Following the method of \cite[\S6]{MornevT}, we construct a representation of $\WK$
on the right $D_\alpha$-vector space 
\UseTheoremCounterForNextEquation
\begin{equation}\label{IntroTMDef}
T_\alpha(M_\Fp)\ \cong\ \Hom\bigl(\SA{\kalg}\complot_{\kalg}\Kalg, 
M_{\Fp,\alpha} \complot_K\Kalg\bigr).
\end{equation}

Since $D_\alpha$ varies with~$\alpha$, to combine these representations into a single representation we must pass to an extension of the coefficient field~$F_\Fp$. A natural choice of extension is the ring $\smash{\ICF\complot_{\BFq}\kalg}$, which is 
isomorphic to the product of finitely many copies of the completion of a maximal unramified extension of~$\ICF$.
We are interested in
\UseTheoremCounterForNextEquation
\begin{equation}\label{IntroUMDef}
U(M_\Fp)\ \cong\ \bigoplus_{\alpha\in\BQ}\, T_\alpha(M_\Fp)\centrot_{\raisebox{-2pt}{$\scriptstyle D_\alpha$}} \SA{\kalg},
\end{equation}
which is a free module over $\smash{\ICF\complot_{\BFq}\kalg}$ of the same rank as~$M_\Fp$ and carries a natural action of~$\WK$. 
If $M_\Fp$ is pure of slope~$0$, this formula simplifies to an isomorphism 
\UseTheoremCounterForNextEquation
\begin{equation}\label{IntroUM0T0}
U(M_\Fp)\ \cong\ T_0(M_\Fp) \complot_{\BFq} \kalg,
\end{equation}
where $T_0(M_\Fp) \cong (M_\Fp\complot_K \Kalg)^\tau$ plays the role of the usual Tate module of~$M_\Fp$. 

The actual construction of $U(M_\Fp)$ is based on an explicit tensor quasi-inverse of the base change functor from $\ICF$-isocrystals over $\kalg$ to $\ICF$-isocrystals over~$\Kalg$.
Thus \eqref{IntroUMDef} defines a tensor functor. Since the category of $F_\Fp$-isocrystals is $\ICF$-linear tannakian, any element of $\WK$ has a well-defined characteristic polynomial on $U(M_\Fp)$ and this polynomial has coefficients in $F_\Fp$.
A simplified version of our main result on compatible systems of representations can now be stated as follows:

\begin{Thm}\label{IntroMainThm}
Let $R$ be a normal integral domain that is finitely generated over~$\BFq$. Consider an $F$-isocrystal $M$ over $K:=\Quot(R)$ and a place $\Fp$ of~$F$. Then for every closed point $x \in\Spec R$ outside a proper closed subset, the action of $\WK$ on $U(M_\Fp)$ is unramified at $x$, and the characteristic polynomial of any geometric Frobenius element at $x$ has coefficients in $F$ and is independent of $\Fp$.
\end{Thm}

See Theorem \ref{CompatibleSystemU} for the precise statement. We also work out specialized versions of this theorem for $A$-motives and for Drinfeld modules, see Theorems \ref{CompatibleSystemUAMot} and \ref{CompatibleSystemDM}.

\medskip
We hope that both this main result, and our construction of $U(M_\Fp)$ as a tensor functor, will find many uses. The former may, for instance, have applications to global problems when information on all places is needed,
and the latter may help to understand the algebraic monodromy group of an $A$-motive or isocrystal, especially when that group is not $\GL_n$.


\bigskip
Now we describe the content of the paper in detail. In Section~\ref{DiffRingIC} we lay out basic properties of dualizable modules over a difference ring. The rest of the paper is divided into one part on local isocrystals and one on global isocrystals.

\medskip
The first six sections of Part I concentrate on $\ICF$-isocrystals over a field~$K$. 
In Section~\ref{LocIsoCrys} we discuss their basic properties, and in Section~\ref{Slopes1} we examine the notion of purity. In Section~\ref{Slopes2} we show that every $\ICF$-isocrystal possesses a unique slope filtration, and we study how this slope filtration behaves under various constructions. 

Over an algebraically closed field, we show that an $\ICF$-isocrystal is determined up to isomorphism by its Newton polygon, much as in the Dieudonn\'e--Manin classification of isogeny classes of $p$-divisible groups. In particular, let $\Kalg$ denote an algebraic closure of~$K$ and by $\kalg$ the algebraic closure of $\BFq$ in~$\Kalg$. Then the base change functor defines an equivalence of tannakian categories from $\ICF$-isocrystals over $\kalg$ to $\ICF$-isocrystals over~$\Kalg$. In Section~\ref{EquivCat} we construct an explicit quasi-inverse $M_\Fp \mapsto U(M_\Fp)$ to this equivalence. 

To any $\ICF$-isocrystal $M_\Fp$ over~$K$ we then also associate a natural $\ICF$-isocrystal $U(M_\Fp)$ over $\kalg$ by first taking its base change to~$\Kalg$ and then applying the functor~$U$. 
In Section~\ref{WeilRep} we use the explicit construction of that functor to endow $U(M_\Fp)$ with a natural action of the Weil group $\WK$ of~$K$. 
The reason that one can only define a representation of $\WK$ is that for Galois elements acting non-trivially on~$\kalg$, one must rescale the action by an \emph{integral} power of~$\tau$.
All this is functorial in $M_\Fp$ and compatible with the tensor product. 

In Section~\ref{LocFinField} we show that for any $\ICF$-isocrystal $M_\Fp$ over a finite field~$k$, the characteristic polynomial of the geometric Frobenius element $\Frob_k^{-1} \in \Wk$ on $U(M_\Fp)$ is equal to the characteristic polynomial of $\smash{\tau_{M_\Fp}^{[k/\BFq]}}$ on $M_\Fp$ itself.

In Section~\ref{GoodRedLoc} we turn to $\ICF$-isocrystals $\CM_\Fp$ defined over an $\BFq$-algebra~$R$, which can be viewed as algebraic families of $\ICF$-isocrystals $\CM_{\Fp,x}$ over all points $x\in X := \Spec R$. After establishing their basic properties, we generalize a fundamental theorem of Watson by showing that for any embedding of a noetherian normal integral domain $R$ into its quotient field~$K$, the base change functor from $\ICF$-isocrystals over~$R$ to $\ICF$-isocrystals over~$K$ is fully faithful. 

In Section~\ref{RedSlopFilt} we discuss when an $\ICF$-isocrystal $\CM_\Fp$ over $R$ is pure and when it possesses a slope filtration. In general a slope filtration does not exist, but when $R$ is a noetherian integral domain, we show that there always exists a slope filtration over $R[\xi^{-1}]$ for some non-zero element $\xi\in R$. We also show that the Newton polygon of $\CM_{\Fp,x}$ depends in an upper semicontinuous way on the point $x\in X$.

In Section~\ref{RedUM} we assume that $R$ is a noetherian normal integral domain with quotient field~$K$, and that $\CM_\Fp$ possesses a slope filtration.
Denoting by $M_\Fp$ the base change of $\CM_\Fp$ to~$K$, for each point $x \in X$ we produce a natural isomorphism $U(M_\Fp) \cong U(\CM_{\Fp,x})$ that is equivariant under the chosen decomposition subgroup of $x$ in the Weil group $\WK$. 
This means in particular that the representation $U(M_\Fp)$ is unramified at~$x$.

\medskip
In Part II we turn to global isocrystals. The basic properties of $F$-isocrystals over a field are established in Section~\ref{GlobCrys}.
In Section~\ref{FinField} we look at an $F$-isocrystal $M$ over a finite field $k$ and show that the characteristic polynomial of $\tau_M^{[k/\BFq]}$ on $M$ equals the characteristic polynomials of the associated $\ICF$-isocrystals $M_\Fp$.

In Section~\ref{GoodRedGlob} we start with an $F$-isocrystal $M$ over a field $K$ and would in principle like to extend it to an $F$-isocrystal over a finitely generated $\BFq$-subalgebra $R$ with quotient field~$K$.
But in general there does not exist a finitely generated $F \otimes_{\BFq} R$-submodule of $M$ on which $\tau_M$ yields an isomorphism. We can, however, achieve an isomorphism by inverting all the conjugates $\sigma^\nu(e)$ of some non-zerodivisor $e\in F\otimes_{\BFq} R$. As this construction depends on the somewhat arbitrary choice of~$e$, we refrain from building a general theory of $F$-iso\-crystals over~$R$, and work instead with an ad hoc model $\CMod$ of $M$ over the indicated localization of $F \otimes_{\BFq} R$.

On the one hand, for every point $x\in X$ this model induces an $F$-isocrystal $\CMod_x$ over the residue field of~$x$. On the other hand, for every place $\Fp$ of $F$ this model induces an $\ICF$-isocrystal $\CMod_\Fp$ over $R[\epsilon_\Fp^{-1}]$ for some non-zero element $\epsilon_\Fp\in R$. Moreover, for almost all $\Fp$ these elements $\epsilon_\Fp$ can be chosen in a systematic way, such that $\epsilon_\Fp(x)\not=0$ if and only if $\Fp\times x$ is disjoint from a certain divisor in $C\times_{\BFq} X$.

In Section~\ref{CompSys} we can then formulate and prove the main result of this paper, saying that the $U(M_\Fp)$ for \emph{all} places $\Fp$ of $F$ form a compatible system of representations of~$\WK$. Here we assume that $R$ is finitely generated over $\BFq$ and normal, so that we can speak of Frobenius elements $\Frob_x\in \WK$ at all closed points ${x\in X}$. The main result says that for every place $\Fp$ of $F$ and every closed point $x$ of~$X$ with $\epsilon_\Fp(x)\not=0$, the inertia group $I_x$ acts trivially on $U(M_\Fp)$, and the characteristic polynomial of the geometric Frobenius element $\Frob_x^{-1}$ is that associated to the $F$-isocrystal $\CMod_x$ over~$k_x$. In particular this characteristic polynomial has coefficients in~$F$ and is independent~of~$\Fp$. Moreover, for almost all $\Fp$ the compatible system comes from a compatible system of representations of the Galois group $\GK$ on the Tate modules $T_0(M_\Fp)$.

\medskip
In Section~\ref{AMot} we apply the preceding results to Anderson $A$-motives. So consider an $A$-motive $\CM$ over a normal integral domain $R$ that is finitely generated over~$\BFq$, and let $M$ denote the associated $A$-motive over $K := \Quot(R)$. Then $M$ gives rise to an $F$-isocrystal $M_F$ over~$K$, and from $\CM$ we can construct a model of $M_F$ over~$R$ as in Section~\ref{GoodRedGlob}. 
Working out what the general machinery from the earlier sections says in this case, we arrive at the following main result for $A$-motives: 

Let $c\colon X\to\Spec A$ be the characteristic morphism associated to~$\CM$. Then for \emph{every} place $\Fp$ of $F$ there exists a divisor $D_\Fp\subset X$, such that $D_\Fp = c^{-1}(\Fp)$ whenever $\Fp\not\in\{\wp,\infty\}$ and that for every closed point $x$ of~$X\smallsetminus D_\Fp$, the inertia group at $x$ acts trivially on $U(M_{F,\Fp})$, and the characteristic polynomial of the geometric Frobenius element $\Frob_x^{-1}$ is that associated to the $A$-motive $\CM_x$ over~$k_x$. In particular this characteristic polynomial has coefficients in~$A$ and is independent~of~$\Fp$.

\medskip
In the final Section~\ref{Drin} we apply these results to the $A$-motive $M$ arising from a Drinfeld $A$-module $\phi$ of rank $r$ over~$K$. If $\wp$ is a maximal ideal of~$A$, let $1\le h\le r$ denote the height of~$\phi$. Then for any place $\Fp\not=\infty$ of~$F$ the rational $\Fp$-adic Tate module $V_\Fp(\phi)$ of~$\phi$ is an $F_\Fp$-vector space of dimension $r$ if $\Fp\not=\wp$, respectively of dimension $r-h$ if $\Fp=\wp$. In either case it carries a natural continuous representation of~$\GK$, and there is a natural $\GK$-equivariant isomorphism of $F_\Fp$-vector spaces
$$V_\Fp(\phi)\ \cong\ T_0(M_{F,\Fp})^\vee \otimes_A\Omega_A.$$

To \emph{every} place $\Fp$ of $F$ we therefore associate the following $\ICF$-isocrystal over~$\kalg$ with an action of~$\WK$:
$$U_\Fp(\phi)\ :=\ U(M_{F,\Fp})^\vee \otimes_A\Omega_A.$$
Then, firstly, for any $\Fp\not=\wp,\infty$ the $\ICF$-isocrystal $M_{F,\Fp}$ is pure of slope $0$ and there is a natural $\WK$-equivariant isomorphism 
$$U_\Fp(\phi)\ \cong\ \unit_{\kalg} \otimes_{F_\Fp} V_\Fp(\phi).$$
Secondly, for $\Fp=\infty$ the $F_\infty$-isocrystal $M_{F,\infty}$ is pure of slope $\minusoneoverr$ and of rank~$r$. As this slope is non-zero, the role of the $\infty$-adic Tate module of $\phi$ is played by a certain one-dimensional left vector space $V_{\infty,\oneoverr}(\phi)$ over a central division algebra $D_\infty$ with Hasse invariant $[\oneoverr]$ over~$F_\infty$. We thus obtain a $\WK$-equivariant isomorphism
$$U_\infty(\phi)\ \cong \SAa{\kalg}{\minusoneoverr}^\vee \otimes_{D_\infty} V_{\infty,\oneoverr}(\phi).$$
Thirdly, for $\Fp=\wp$ the $F_\wp$-isocrystal $M_{F,\wp}$ has slope~$0$ with multiplicity $r-h$ and slope~$\oneoverh$ with multiplicity~$h$. Here the slope $0$ part corresponds to the Tate module $V_\wp(\phi)$, and the slope $\oneoverh$ part gives rise to a certain one-dimensional left vector space $V_{\wp,\minusoneoverh}(\phi)$ over a central division algebra $D_\wp$ with Hasse invariant $[\minusoneoverh]$ over~$F_\wp$. We thus obtain a $\WK$-equivariant isomorphism
$$U_\wp(\phi)\ \cong\ 
\bigl[\unit_{\kalg} \otimes_{F_\wp} \kern-4pt V_\wp(\phi) \bigr]
\oplus
\bigl[\SAa{\kalg}{\oneoverh}^\vee \kern-2pt \otimes_{D_\wp} \kern-4pt V_{\wp,\minusoneoverh}(\phi) \bigr].$$

Finally suppose that $\phi$ is defined over a normal integral domain $R$ that is finitely generated over~$\BFq$ and has quotient field~$K$. For every~$\Fp$ we can then explicitly describe over which localization $R[\epsilon_\Fp^{-1}]$ the $F_\Fp$-isocrystal $M_{F,\Fp}$ extends and possesses a slope filtration. In the same way as for $A$-motives, the $U_\Fp(\phi)$ then form a compatible system of $\WK$-representations for all places $\Fp$ of~$F$. 

For $\Fp=\infty$, using a basis of $V_{\infty,\oneoverr}(\phi)$ over $D_\infty$, the representation of $\WK$ corresponds to a continuous homomorphism ${\WK\to(D_\infty^\opp)^\times}$. Our result for $\Fp=\infty$ thus reproduces that of J.-K. Yu \cite[Thm.\,3.4]{Yu2003}.

Similarly, for $\Fp=\wp$, using a basis of $V_{\wp,\minusoneoverh}(\phi)$ over $D_\wp$, the representation of $\WK$ on it corresponds to a continuous homomorphism ${\WK\to(D_\wp^\opp)^\times}$. 
Likewise, using a basis of the $(r-h)$-dimensional $F_\wp$-vector space $V_\wp(\phi)$, the representation of $\GK$ corresponds to a continuous homomorphism $\GK\to\GL_{r-h}(F_\wp)$. Together the representation of $\WK$ on $U_\wp(\phi)$ is therefore  given by a continuous homomorphism
$$\WK \longto \GL_{r-h}(F_\wp) \times (D_\wp^\opp)^\times.$$
In work in preparation, the second author plans to study this representation in more detail.

\bigskip
\noindent
\textbf{Acknowledgements.}
The first author is supported by an Ambizione fellowship of the Swiss National Science Foundation (project PZ00P2\_202119)
and by Grants CNS2023-145167 and CEX2023-001347-S, funded by MICIU/AEI/10.13039/501100011033.

\section{Modules over a difference ring}
\label{DiffRingIC}

Following Kedlaya \cite[Def.\,14.1.1]{KedlayaD} we consider a \emph{difference ring~$\CE$}, that is, a commutative unitary ring together with an endomorphism $\sigma\colon \CE \to \CE$. 
Let $\CE[\tau]$ denote the associated  \emph{twisted polynomial ring} 
consisting of all finite formal expressions $\sum_{i\ge0}e_i\tau^i$ with $e_i\in \CE$, with the usual addition and the unique distributive multiplication extending that on~$\CE$ which satisfies
\UseTheoremCounterForNextEquation
\begin{equation}
\tau \, e = \sigma(e) \, \tau
\end{equation}
for all $e \in \CE$. The ring $\CE[\tau]$ is associative and unitary, but not commutative unless $\sigma$ is the identity map.


\medskip
Giving a left $\CE[\tau]$-module is equivalent to giving an $\CE$-module $M$ together with an additive map $\tau_M\colon M \to M$ that satisfies $\tau_M(e \, m) = \sigma(e) \, \tau_M(m)$ for all $e\in\CE$ and $m\in M$. 
Giving $\tau_M$ is further equivalent to giving the $\CE$-linear map
\UseTheoremCounterForNextEquation
\begin{equation}\label{TauMLin}
\taulin_M\colon\ \sigma^* M := \CE\otimes_{\CE,\sigma}M\ \longto\ M
\end{equation}
with $\taulin_M(e\otimes m) = e\,\tau_M(m)$ for all $e\in\CE$ and $m\in M$. We call $\taulin_M$ the \emph{structure morphism} of~$M$. 
For the subset of $\tau$-invariants of $M$ we use the following notation:
\UseTheoremCounterForNextEquation
\begin{equation}\label{MtauDef}
M^\tau := \{m\in M\mid \tau m=m\}.
\end{equation}
This is a module over the subring of $\sigma$-invariants
\UseTheoremCounterForNextEquation
\begin{equation}\label{ESigmaDef}
\CE^\sigma := \{ e \in \CE \mid \sigma(e) = e \}.
\end{equation}
The set of homomorphisms of left $\CE[\tau]$-modules is also an $\CE^\sigma$-module and abbreviated by
\UseTheoremCounterForNextEquation
\begin{equation}\label{DMHomDef}
\Hom(M,N)\ :=\ \Hom_{\CE[\tau]}(M,N).
\end{equation}
The \emph{tensor product} of two left $\CE[\tau]$-modules is the left $\CE[\tau]$-module 
\UseTheoremCounterForNextEquation
\begin{equation}\label{DiffRingICTensorDef}
M \otimes N := M\otimes_\CE N
\qquad\hbox{with}\qquad 
\taulin_{M\otimes N} := \taulin_M \otimes_\CE \taulin_N. 
\end{equation}
The \emph{unit object} is the left $\CE[\tau]$-module
\UseTheoremCounterForNextEquation
\begin{equation}\label{DiffRingICUnitDef}
\unit_\CE := \CE
\qquad\hbox{with}\qquad 
\taulin_{\unit_\CE} = 1.
\end{equation}
Clearly the tensor product is functorial in $M$ and~$N$, and it is symmetric and associative with the unit object~$\unit_\CE$, turning the category of left $\CE[\tau]$-modules into a symmetric monoidal category. Also, for any integer $n\ge0$, the \emph{$n$-th exterior power} of a left $\CE[\tau]$-module $M$ is the left $\CE[\tau]$-module
\UseTheoremCounterForNextEquation
\begin{equation}\label{DiffRingICExtPowDef}
\bigwedge^nM := \bigwedge^n_{\CE}(M)
\qquad\hbox{with}\qquad 
 \taulin_{\bigwedge^nM} := \bigwedge^n_{\CE}(\taulin_M).
\end{equation}

\medskip
The following definition generalizes \cite[Def.\,14.1.4]{KedlayaD}:

\begin{PropDef}\label{DiffRingICDef}
For any left $\CE[\tau]$-module $M$ the following are equivalent:
\begin{enumerate}\StatementLabels
\item\label{DiffRingICDefA}
The module $M$ is finitely generated projective over~$\CE$ and $\taulin_M$ is an isomorphism.
\item\label{DiffRingICDefB}
The object $M$ possesses a dual in the sense of Deligne \cite[(2.1.2)]{DeligneT}, that is, there exist a left $\CE[\tau]$-module $M^\vee$ and morphisms $M\otimes M^\vee\arrover{\ev}\unit_\CE$ and $\unit_\CE \arrover{\coev} M^\vee \otimes M$ such that the composites 
$$\vcenter{\vskip-7pt\xymatrix@C+30pt@R-15pt{
M \ar[r]^-{\id_M\otimes\mkern+2mu\coev} & M \otimes M^\vee \otimes M \ar[r]^-{\ev\mkern-1mu\otimes\mkern+2mu\id_M} & M, \\
\kern2pt M^\vee \ar[r]^-{\coev\mkern-2mu\otimes\mkern+2mu\id_{M^\vee}} & \kern2pt M^\vee\otimes M \otimes M^\vee \ar[r]^-{\id_{M^\vee}\mkern-1mu\otimes\mkern+2mu\ev} & M^\vee}}$$
are the identity morphisms.
\end{enumerate}
Moreover the triple $(M^\vee,\ev,\coev)$ in \ref{DiffRingICDefB} is unique up to unique isomorphism. A left $\CE[\tau]$-module $M$ with the above properties is called \emph{dualizable} and $M^\vee$ is called its \emph{dual}.
\end{PropDef}


\begin{Proof}
If $M$ satisfies \ref{DiffRingICDefA}, its dual $M^\vee := \Hom_\CE(M,\CE)$ as an $\CE$-module comes with the usual homomorphisms $M\otimes_\CE M^\vee\arrover{\ev}\CE$ and $\CE \arrover{\coev} M^\vee \otimes_\CE M$ satisfying the composition properties in~\ref{DiffRingICDefB}. For any $f\in M^\vee$ define $\tau f$ as the composite map
\UseTheoremCounterForNextEquation
\begin{equation}\label{DiffRingICDualTauDef}
\vcenter{\vskip-2pt\xymatrix@C+15pt{M \ar[r]^-{(\taulin_M)^{-1}}_-{\sim} & \sigma^*M \ar[r]^-{\sigma^*f} & \sigma^*\CE \ar[r]^{\sim} & \CE.}}
\end{equation}
This turns $M^\vee$ into a left $\CE[\tau]$-module, and a direct computation shows that $\ev$ and $\coev$ commute with~$\tau$. Thus $M$ satisfies the condition~\ref{DiffRingICDefB}. 

Conversely let $(M^\vee,\ev,\coev)$ be as in~\ref{DiffRingICDefB}. Forgetting the action of~$\tau$, this data constitutes a dual of $M$ in the category of $\CE$-modules; hence $M$ is finitely generated projective over~$\CE$ by \cite[Prop.\,2.6]{DeligneT}.
On the other hand the conditions in~\ref{DiffRingICDefB} imply that $\sigma^*\ev$ and $\sigma^*\coev$ turn $\sigma^*M^\vee$ into an $\CE$-module dual of~$\sigma^*M$. Moreover the fact that $\ev$ and $\coev$ commute with $\tau$ amounts to the equalities
$$\ev \circ (\taulin_M \otimes \taulin_{M^\vee}) = \sigma^*(\ev)\quad\text{and}\quad
(\taulin_{M^\vee} \otimes \taulin_M) \circ \sigma^*(\coev) = \coev.$$
The $\CE$-module homomorphism $\taulin_{M^\vee}$ is thus a \emph{contragredient} of $\taulin_M$ in the sense of \cite[\S2.4]{DeligneT}. By \cite[\S2.4]{DeligneT} it follows that $\taulin_M$ is an isomorphism, so $M$ satisfies the conditions in~\ref{DiffRingICDefA}.

Finally, the uniqueness of $(M^\vee,\ev,\coev)$ follows from \cite[\S2.2]{DeligneT}.
\end{Proof}

\medskip
The proof of Proposition~\ref{DiffRingICDef} gives a direct construction of $M^\vee$. The condition \ref{DiffRingICDef}~\ref{DiffRingICDefA} thus implies that duals, tensor products and exterior powers of dualizable left $\CE[\tau]$-modules are again dualizable.
Moreover, to any dualizable left $\CE[\tau]$-modules $M$ and $N$ one can associate their \emph{inner hom}
\UseTheoremCounterForNextEquation
\begin{equation}\label{DiffRingICHomDef}
\IHom(M,N) := \Hom_\CE(M,N)
\qquad\hbox{with}\qquad 
\tau f := \taulin_N \circ (\sigma^*\!f) \circ (\taulin_M)^{-1},
\end{equation}
and this is again a dualizable left $\CE[\tau]$-module. By the construction of $M^\vee$ 
this comes with a natural isomorphism $\IHom(M,N) \cong M^\vee \otimes N$. 
Also, the inner hom is functorial in $M$ and~$N$. One easily verifies that $\hom(M,N)$ represents the functor $T \mapsto {\Hom(T \otimes M, \,N)}$, so there is a natural adjunction isomorphism
\UseTheoremCounterForNextEquation
\begin{equation}\label{HomhomAdj}
\Hom(T\otimes M,N)\ \cong\ \Hom(T,\hom(M,N))
\end{equation}
that is functorial in $T,M,N$. In the special case $T = \unit_\CE$ this amounts to an equality
\UseTheoremCounterForNextEquation
\begin{equation}\label{Homhomtau}
\Hom(M,N)\ =\ \hom(M,N)^\tau.
\end{equation}




\begin{Cons}\label{BC1}
(Base change) \rm  
Consider a morphism of difference rings $f\colon \CE \to \CE'$, that is, a homomorphism of  the underlying rings that commutes with the structure endo\-mor\-phisms~$\sigma$. To any left $\CE[\tau]$-module $M$ we can then associate the left $\CE'[\tau]$-module 
$$f^*M := \CE'\otimes_\CE M
\qquad\hbox{with}\qquad 
\taulin_{f^*M} := \id\otimes\taulin_M.$$
This construction is functorial in~$M$ and compatible with tensor product, unit objects and exterior powers, and therefore defines a monoidal functor. It is also preserves dualizable objects and is compatible with the inner hom and the duals of those.
\end{Cons}

\begin{Prop}\label{Tannaka}
Suppose that $\CE$ is a non-empty finite product of fields that are transitively permuted by~$\sigma$. Then:
\begin{enumerate}\StatementLabels
\item\label{TannakaFree}
Every dualizable left $\CE[\tau]$-module is free of a unique finite rank as an $\CE$-module.
\item\label{TannakaFField}
The ring of invariants $\CE^\sigma$ is a field.
\item\label{TannakaTannaka}
The category of dualizable left $\CE[\tau]$-modules is $\CE^\sigma$-linear tannakian.
\item\label{TannakaFiber}
The forgetful functor to the category of $\CE$-modules is a fiber functor and in particular is exact.
\end{enumerate}
\end{Prop}

\begin{Proof}
Write $\CE$ as a product of fields $\CE_i$ with $\sigma(\CE_i)\subset \CE_{i+1}$ for all $i\in\BZ/n\BZ$. Then any dualizable left $\CE[\tau]$-module has a decomposition $M = \prod_i M_i$ for $\CE_i$-modules $M_i$ with ${\tau M_i\subset M_{i+1}}$. The isomorphy of $\taulin_M$ then implies that $\dim_{\CE_i}(M_i)$ is independent of~$i$, proving \ref{TannakaFree}. Next $\CE' := \CE_0$ is a difference ring with respect to the endomorphism $\sigma' := \sigma^d|\CE_0$. Also the projection $\CE\onto \CE'$ induces an isomorphism $\CE^\sigma \isoto (\CE')^{\sigma'}$ where the latter is a field, proving~\ref{TannakaFField}.

Next the assumption on $\CE$ implies that $\CE$ is noetherian and every $\CE$-module is projective. Thus for any morphism of dualizable left $\CE[\tau]$-modules $f\colon M \to N$, the kernel $P$ and the cokernel $Q$ within the category of left $\CE[\tau]$-modules are again finitely generated projective over~$\CE$. It also follows that every $\CE$-module is flat and the functor $\sigma^*$ is therefore exact. Thus we obtain a commutative diagram with exact rows 
\UseTheoremCounterForNextEquation
\begin{equation}\label{TannakaDiag}
\vcenter{\xymatrix{
0 \ar[r] & \sigma^* P \ar[r] \ar[d]_{\taulin_P} & \sigma^* M \ar[r]^-{\sigma^*f}
 \ar[d]_{\taulin_M}^{\kern-1pt\wr} & \sigma^* N  \ar[d]_{\taulin_N}^{\kern-1pt\wr} \ar[r] & \sigma^* Q  \ar[d]_{\taulin_Q} \ar[r] & 0 \\
0 \ar[r] & P \ar[r] & M \ar[r]^-f & N \ar[r] & Q \ar[r] & 0.\\}}
\end{equation}
Since the structure morphisms of $M$ and $N$ are isomorphisms, by the 5-lemma the same follows for $P$ and~$Q$. This proves that $P$ and $Q$ are dualizable.
We conclude that the category of dualizable left $\CE[\tau]$-modules is abelian and that the forgetful functor to the category of $\CE$-modules is exact.

This functor is also faithful. By construction $\CE^\sigma$ is the ring of endomorphisms of the unit object~$\unit_\CE$. As this ring is a field, the category of dualizable left $\CE[\tau]$-modules is therefore $\CE^\sigma$-linear tannakian by \cite[\S2.1, \S2.8]{DeligneT}. Together this proves \ref{TannakaTannaka} and \ref{TannakaFiber}.
\end{Proof}

\medskip
Due to Proposition \ref{Tannaka} \ref{TannakaFree} the next definition makes sense:

\begin{Def}\label{LICRankDef}
In the situation of Proposition \ref{Tannaka} the \emph{rank} of a dualizable left $\CE[\tau]$-module is the rank of the underlying free $\CE$-module.
\end{Def}


%
%

\newpage
\bigskip
\noindent{\huge\bf Part I: Local isocrystals}
\addtocontents{toc}{\medskip}
\addcontentsline{toc}{section}{\large Part I: Local isocrystals}

\section{Local isocrystals over a field}
\label{LocIsoCrys}

Throughout the following we fix a finite field $\BFq$ of order~$q$.
In this part we also fix a local field $\ICF$ over $\BFq$ with ring of integers $\ICO$ and residue field $\ICk$ of degree $\ICd$ over~$\BFq$. For convenience we choose a uniformizer $\ICu \in \ICO$, but the results and definitions below are independent of this choice.

Let $K$ be an arbitrary field over~$\BFq$. The objects discussed below will have the \emph{ground field} $\BFq$, the \emph{field of coefficients}~$\ICF$, and the \emph{field of definition}~$K$. We equip $\ICF$ with its locally compact topology and $K$ with the discrete topology and consider the completed tensor product ring
\UseTheoremCounterForNextEquation
\begin{equation}\label{FPKDef}
\EFpK := \ICF\complot_{\BFq}\kern-1ptK,
\end{equation}
viewed as a difference ring via the partial Frobenius endomorphism $\sigma(a\complot\xi):= a\complot\xi^q$.

\begin{Def}\label{ICDef}
A dualizable left $\EFpK[\tau]$-module is called an \emph{$\ICF$-isocrystal over $K$ (with the ground field~$\BFq$)}.
\end{Def}

We denote the category of $\ICF$-isocrystals over $K$ by $\LIC{K}$. By \S\ref{DiffRingIC} this category carries a rigid symmetric monoidal structure.

\begin{Prop}\label{LICTannakian}
\begin{enumerate}\StatementLabels
\item\label{LICTannakianTannakian}
The category $\LIC{K}$ is $\ICF$-linear tannakian.
\item\label{LICTannakianFiber}
The forgetful functor to the category of $\EFpK$-modules is a fiber functor.
\item\label{LICTannakianConstRank}
Every $\ICF$-isocrystal is free of a unique finite rank as an $\EFpK$-module.
\end{enumerate}
\end{Prop}

\begin{Proof}
The choice of the uniformizer $\ICu$ determines an isomorphism $\ICF\cong\ICk\ppu$ and hence an isomorphism
\UseTheoremCounterForNextEquation
\begin{equation}\label{FPKDefVar}
\EFpK\ \cong\ \bigl(\ICk\otimes_{\BFq}\kern-2ptK\bigr)\ppu.
\end{equation}
Here the fact that $\ICk$ is a finite separable extension of $\BFq$ implies that $\ICk\otimes_{\BFq} K$ is a non-empty finite product of fields that are transitively permuted by the endomorphism $a\otimes\xi\mapsto a\otimes\xi^q$. Therefore $\EFpK$ is a non-empty finite product of fields that are transitively permuted by~$\sigma$. Furthermore, since $\sigma$ acts trivially on $\ICk$ and~$\ICu$, and the fixed field of $a\mapsto a^q$ on $K$ is $\BFq$, the subring of $\sigma$-invariants is the image of the field $\ICF\cong\ICk\ppu$ under the map $a\mapsto a\complot1$. Thus everything follows from Proposition \ref{Tannaka}.
\end{Proof}

\begin{Cor}\label{LICSimpleFinite}
\begin{enumerate}\StatementLabels
\item Every $\ICF$-isocrystal is a finite iterated extension of simple objects.
\item For any $\ICF$-isocrystals $M$ and $N$ the $\ICF$-vector space $\Hom(M,N)$ is finite-dimensional.
\end{enumerate}\end{Cor}

\begin{Proof}
This follows from Proposition \ref{LICTannakian} \ref{LICTannakianTannakian} by \cite[Prop.\,2.13]{DeligneT}.
\end{Proof}


\begin{Cons}\label{LICBC}
(Base change) \rm  
Any field embedding $i\colon K\into K'$ induces a morphism of difference rings $\id\complot i\colon \EFp{K} \into \EFp{K'}$ that we again abbreviate by~$i$. From Construction \ref{BC1} we thus obtain a natural $F_\Fp$-linear monoidal base change functor
$$i^*\colon\ \LIC{K} \longto \LIC{K'},\quad M\mapsto i^*M \cong M\complot_{K,i}K'.$$
This functor preserves the rank, and is exact and faithful because the ring homomorphism $i\colon\EFpK \to \EFp{K'}$ is faithfully flat.
\end{Cons}

\begin{Cons}\label{LICCoeffRes}
(Restriction of the coefficient field) \rm 
Since $\ICF$ is a field extension of degree $\ICd$ over the subfield $\ICFalt:=\BFq\ppu$, the ring $\EFpK$ is a free module of rank $\ICd$ over $\CE_{\ICFalt,K}$. Giving an $\ICF$-isocrystal $M$ over $K$ is thus equivalent to giving an $\ICFalt$-isocrystal $M'$ over $K$ together with an embedding $\ICF\into\End(M')$ over~$\ICFalt$. Since $\CE_{\ICFalt,K} \cong K\ppu$, this will simplify some computations.
\end{Cons}

\section{Purity}
\label{Slopes1}

Keeping the situation of Section \ref{LocIsoCrys}, we consider the following difference subring
of $\EFp{K}$:
\UseTheoremCounterForNextEquation
\begin{equation}\label{APKDef}
\EAp{K} := \ICO\complot_{\BFq}\kern-1ptK
\end{equation}
By a \emph{lattice} in an $\ICF$-isocrystal $M$ we mean a finitely generated free $\EAp{K}$-submodule $L$ such that $L[\ICu^{-1}] = M$. Recall that the residue field $\ICk$ has degree $\ICd$ over~$\BFq$.

\begin{Def}\label{ICSlopeDef}
An $\ICF$-isocrystal $M$ over $K$ is called \emph{pure of slope~$\alpha\in\BQ$} if there exist an integer $n>0$ with $n\alpha\in\BZ$ and a lattice $L\subset M$ such that $\tau^{n\ICd}L$ generates the lattice $\ICu^{n\alpha} L$.
\end{Def}

If this condition holds with a given~$n$, it clearly also holds with any multiple thereof.
By \cite[Props.\,3.4.13--14]{MornevT} we have:

\begin{Prop}\label{ICSlopeProp}
\begin{enumerate}\StatementLabels%
\item\label{ICSlopePropB}
Any nonzero pure $\ICF$-isocrystal has a unique slope.
\item\label{ICSlopePropA}
Any morphism between pure $\ICF$-isocrystals of different slopes is zero.
\end{enumerate}
\end{Prop}

\medskip
Next we study how purity interacts with basic constructions.
 
\begin{Lem}\label{CapLattice}
Let $i\colon K \into K'$ be a field embedding and let $M$ be an $\ICF$-isocrystal over~$K$. Then for any lattice $L'\subset i^*M$, the intersection $L' \cap M$ is a lattice in~$M$.
\end{Lem}

\begin{Proof}
Take any lattice $L\subset M$. Then $i^*L$ is a lattice in~$i^*M$; hence there exists an integer $e\ge 0$ such that $\ICu^e i^*L \subset L' \subset \ICu^{-e} i^* L$. Intersecting with $M$ we infer that $\ICu^e L \subset L' \cap M \subset \ICu^{-e} L$, so $L'\cap M$ is a lattice in~$M$.
\end{Proof}

\begin{Prop}\label{LICBCPure}
For any field embedding $i\colon K\into K'$, an $\ICF$-isocrystal $M$ over $K$ is pure of slope $\alpha$ if and only if $i^*M$ is pure of slope~$\alpha$.
\end{Prop}

\begin{Proof}
The ``only if'' part is clear. For the converse take an integer $n>0$ with $n\alpha\in\BZ$ and a lattice $L'\subset i^* M$ such that $\tau^{n\ICd}L'$ generates $\ICu^{n\alpha} L'$. Then $L_0 := L' \cap M$ is a lattice in~$M$ by Lemma \ref{CapLattice}. Since $i^* L_0$ is a lattice in $i^*M$, there exists an integer $e\ge 0$ such that $\ICu^e L' \subset i^* L_0$.

Next, the inclusion $\ICu^{-n\alpha}\tau^{n\ICd}L' \subset L'$ implies that $\ICu^{-n\alpha}\tau^{n\ICd}L_0 \subset L_0$. Writing $L_j$ for the lattice generated by $(\ICu^{-n\alpha}\tau^{n\ICd})^jL_0$, we thus obtain a descending sequence of lattices ${L_0 \supset L_1 \supset \dotsc}$ in~$M$.
As $L'$ is generated by $(\ICu^{-n\alpha}\tau^{n\ICd})^jL'$, the inclusion $\ICu^e L' \subset i^* L_0$ implies that $\ICu^{e}L'\subset i^*L_j$ for every~$j$. Since the quotient $i^*L_0/\ICu^{e} L'$ has finite dimension over~$K'$, it follows that $i^* L_j = i^* L_{j+1}$ for some $j\ge 0$. As the homomorphism $\ICO\complot_{\BFq} K \to \ICO\complot_{\BFq} K'$ is faithfully flat, we conclude that $L_j = L_{j+1}$. But by construction $L_{j+1}$ is the lattice generated by $\ICu^{-n\alpha}\tau^{n\ICd} L_j$; hence $\ICu^{n\alpha} L_j = \ICu^{n\alpha} L_{j+1}$ is generated by $\tau^{n\ICd}L_j$. Thus $M$ is pure of slope~$\alpha$, as desired.
\end{Proof}

\begin{Prop}\label{InvertSlopeInt}
Every $\ICF$-isocrystal of rank $1$ is pure of an integral slope.
\end{Prop}

\begin{Proof}
By Proposition \ref{LICBCPure} we may replace $K$ by any field extension. We may therefore assume that $\ICk$ possesses an embedding into~$K$. Then $\ICk\otimes_{\BFq} K$ is isomorphic to a direct product of $\ICd$ copies of $K$ that are transitively permuted by~$\sigma$. Correspondingly, the isomorphism \eqref{FPKDefVar} implies that $\EFpK$ is isomorphic to a direct product of $\ICd$ copies of $K\ppu$ that are transitively permuted by~$\sigma$.

Now consider any $\ICF$-isocrystal $M$ of rank $1$ over~$K$. Then the underlying $\EFpK$-module is a direct product of a $1$-dimensional $K\ppu$-vector space for each simple factor of $\EFpK$. Fix any free $K\bbu$-submodule $L_0$ of rank $1$ in one of them. Since $\tau^{\ICd}$ sends the corresponding factor of $M$ to itself, and $K\bbu$ is a discrete valuation ring with uniformizer~$\ICu$, it follows that $\tau^{\ICd}L_0 = \ICu^\alpha L_0$ for some integer~$\alpha$. For each $0<i<\ICd$ consider the free $K\bbu$-submodule $L_i$ of rank $1$ that is generated by $\tau^iL_0$. Then $L := \smash{\prod_{i=0}^{\ICd-1}L_i}$ is a lattice in~$M$ such that $\tau^{\ICd}L = \ICu^\alpha L$. Thus $M$ is pure of slope~$\alpha$, as desired.
\end{Proof}

\begin{Rem}\label{InvertSlopeIntRem}\rm
In general not every integer occurs as a slope of a pure $\ICF$-isocrystal of rank~$1$. For example, by \cite[Prop.\,4.2.1]{MornevT} the slope $1$ occurs if and only if $\ICk$ embeds into~$K$.
\end{Rem}

\begin{Prop}\label{LICCoeffResPure}
In the situation of Construction \ref{LICCoeffRes}, an $\ICF$-isocrystal $M$ is pure of slope $\alpha$ if and only if the corresponding $\ICFalt$-isocrystal $M'$ is pure of slope~$\frac{\alpha}{\ICd}$.
\end{Prop}

\begin{Proof}
If $M$ is pure of slope~$\alpha$, choose $n,L$ as in Definition \ref{ICSlopeDef} such that $\tau^{n\ICd}L$ generates the lattice $\ICu^{n\alpha} L$. Since the field $\ICFalt$ has residue degree $1$ in place of~$\ICd$, the same lattice shows that $M'$ is pure of slope $\frac{\alpha}{\ICd}$. Conversely, for any lattice $L'\subset M'$ that makes $M'$ pure of slope~$\frac{\alpha}{\ICd}$, the submodule $\ICO{\cdot}L'$ is a lattice that makes $M$ pure of slope~$\alpha$.
\end{Proof}


\begin{Prop}\label{TensorPure}
Consider any $\ICF$-isocrystals  $M$ and $N$ over~$K$, such that $M$ is pure of slope $\alpha$ and $N$ is pure of slope~$\beta$. Then:
\begin{enumerate}\StatementLabels
\item\label{TensorPureA}
$M\otimes N$ is pure of slope $\alpha+\beta$.
\item\label{TensorPureB}
$M^\vee$ is pure of slope $-\alpha$.
\item\label{TensorPureC}
$\bigwedge^nM$ is pure of slope $n\alpha$ for any $n\ge0$. 
\end{enumerate}
\end{Prop}

\begin{Proof}
To prove \ref{TensorPureA}, by Definition \ref{ICSlopeDef} we can choose $n>0$ with $n\alpha,n\beta\in\BZ$ and lattices $L\subset M$ and $\varLambda\subset N$ such that
$\tau^{n\ICd}L$ generates $\ICu^{n\alpha} L$ and $\tau^{n\ICd}\varLambda$ generates $\ICu^{n\beta}\varLambda$. Thus $L\otimes \varLambda$ is a lattice in $M\otimes N$ such that $\tau^{n\ICd} (L\otimes \varLambda)$ generates $\smash{\ICu^{n(\alpha+\beta)}(L\otimes \varLambda)}$; hence $M\otimes N$ is pure of slope $\alpha+\beta$.

Assertions \ref{TensorPureB} and \ref{TensorPureC} are proved similarly.
\end{Proof}

\medskip
Next, by \cite[Props.\,4.3.4 and 4.3.7]{MornevT} we have:

\begin{Prop}\label{StdPureICCons}
Suppose that $K$ is separably closed, and let $\alpha$ be a rational number with precise denominator~$r>0$.
\begin{enumerate}\StatementLabels
\item\label{StdPureICConsA}
There exists a simple $\ICF$-isocrystal $\SA{K}$ of rank $r$ over $K$ that is pure of slope~$\alpha$.
\item\label{StdPureICConsB}
Every pure $\ICF$-isocrystal of slope $\alpha$ over $K$ is isomorphic to a direct sum of copies of~$\SA{K}$. In particular $\SA{K}$ is unique up to isomorphism.
\item\label{StdPureICConsC}
The ring $\End(\SA{K})$ is a central divison algebra of dimension $r^2$ over $\ICF$ and with Hasse invariant $[-\alpha] \in \BQ/\BZ$.
\end{enumerate}
\rm Here we use the normalization of the Hasse invariant from Laumon \cite[Appendix (A.2)]{LaumonI}, but the opposite choice of sign is also common in the literature.
\end{Prop}


\begin{Prop}\label{StdPureICPullback}
Consider any embedding of separably closed fields $i\colon K\into K'$. Then for every $\alpha\in\BQ$ there exists an isomorphism $i^*\SA{K} \cong \SA{K'}$ and the resulting homomorphism ${\End(\SA{K})\longto\End(\SA{K'})}$ is an isomorphism.
\end{Prop}

\begin{Proof}
Let $r>0$ be the precise denominator of~$\alpha$. Then $i^*\SA{K}$ is pure of slope $\alpha$ by Proposition \ref{LICBCPure} and again of rank~$r$; hence by Proposition \ref{StdPureICCons} \ref{StdPureICConsB} over $K'$ it is isomorphic to $\SA{K'}$. Also the resulting map ${\End(\SA{K}) \to \End(\SA{K'})}$ is a homomorphism of division algebras of the same dimension $r^2$ over $\ICF$ and therefore an isomorphism.
\end{Proof}

\begin{Prop}\label{PureCritHartl}
Let $i\colon K \into K'$ be an embedding with $K'$ separably closed. Then an $\ICF$-isocrystal $M$ over $K$ is pure of slope $\alpha$ if and only if $i^* M$ decomposes into a direct sum of copies of $\SA{K'}$.
\end{Prop}
\begin{Proof}
Combine Propositions \ref{LICBCPure} and \ref{StdPureICCons} \ref{StdPureICConsB}.
\end{Proof}

\begin{Prop}\label{PureSubquot}
Any subquotient of a pure $\ICF$-isocrystal is pure of the same slope.
\end{Prop}

\begin{Proof}
This follows from Proposition \ref{PureCritHartl} because the isocrystals $\SA{K'}$ are simple.
See also \cite[Prop.\,4.3.1]{MornevT} for a direct argument with lattices.
\end{Proof}


\begin{Prop}\label{NotNeutral}
If $K$ is separably closed, the tannakian category $\LIC{K}$ is not neutral.
\end{Prop}

\begin{Proof}
Pick any $\alpha\in\BQ$ of precise denominator $r>1$. Then by Proposition~\ref{StdPureICCons} the $\ICF$-isocrystal $\SA{K}$ has rank~$r$, and $\End(\SA{K})$ is a non-split central simple algebra of rank $r^2$ over~$\ICF$. 
Suppose now that $\LIC{K}$ possesses a fiber functor $V$ to the category of finite dimensional $\ICF$-vector spaces. 
Then $V$ is faithful and commutes with exterior powers; hence $V(\SA{K})$ is an $\ICF$-vector space of dimension~$r$. For dimension reasons the injective homomorphism of $\ICF$-algebras
$\End(\SA{K}) \into \End_{\ICF}(V(\SA{K}))$
induced by $V$ is therefore an isomorphism. But this contradicts the fact that $\End(\SA{K})$ is non-split.
\end{Proof}


\section{The slope filtration}
\label{Slopes2}

\begin{Def}\label{FiltDef}
A \emph{slope filtration} of an $\ICF$-isocrystal $M$ over $K$ is a separated exhaustive right continuous ascending filtration by $\ICF$-isocrystals $M_{\le\alpha}$ over $K$ indexed by all $\alpha\in\BQ$, such that every subquotient
$$M_\alpha := M_{\le\alpha}/M_{<\alpha}
\qquad\hbox{with}\qquad
M_{<\alpha} := \smash{\bigcup_{\beta<\alpha}M_{\le\beta}}$$
is pure of slope~$\alpha$. 
\end{Def}

Since every $\ICF$-isocrystal is an iterated extension of finitely many simple objects, saying that the slope filtration is right continuous amounts to saying that for every $\alpha$ there exists an $\varepsilon > 0$ such that $M_{\le\alpha} = M_{\le\alpha+\varepsilon}$.

\begin{Prop}\label{FiltFunctUniq}
\begin{enumerate}\StatementLabels%
\item\label{FiltFunctUniqA}
For any morphism of $\ICF$-isocrystals $f\colon M \to N$ admitting slope filtrations and every rational number $\alpha$ we have $f(M_{\le\alpha}) \subset N_{\le\alpha}$.
\item\label{FiltFunctUniqB}
In particular the slope filtration is unique if it exists.
\end{enumerate}
\end{Prop}

\begin{Proof}
The conclusion in \ref{FiltFunctUniqA} is equivalent to saying that the morphism $M_{\le \alpha}\to N/N_{\le\alpha}$ induced by $f$ is zero. Here $M_{\le \alpha}$ is a successive extension of pure $\ICF$-isocrystals of slopes $\le\alpha$, and 
$N/N_{\le\alpha}$ is a successive extension of pure $\ICF$-isocrystals of slopes $>\alpha$.
Thus \ref{FiltFunctUniqA} follows by induction from Proposition \ref{ICSlopeProp} \ref{ICSlopePropA}.

Applying \ref{FiltFunctUniqA} to the identity morphism on an $\ICF$-isocrystal with two slope filtrations implies \ref{FiltFunctUniqB}.
\end{Proof}

\begin{Prop}\label{LICBCFilt}
For any field embedding $i\colon K\into K'$ the base change functor $i^*$ preserves the slope filtration, if it exists.
\end{Prop}

\begin{Proof}
Direct consequence of Proposition \ref{LICBCPure}.
\end{Proof}


\begin{Prop}\label{FiltExist}%
Any $\ICF$-isocrystal $M$ over $K$ possesses a slope filtration.
\end{Prop}

\begin{Proof}
In the situation of Construction \ref{LICCoeffRes}, Proposition \ref{LICCoeffResPure} implies that any slope filtration of the associated $\ICFalt$-isocrystal $M'$ is a slope filtration of $M$ with the indices scaled by~$\ICd$. Thus we may without loss of generality assume that $\ICF=\BFq\ppu$ and hence $\EFK =K\ppu$.
We can then use the results of Hartl \cite[\S1.5]{HartlP}. In particular Proposition \ref{PureCritHartl} implies that an isocrystal is pure of slope $\alpha$ if and only if it is isoclinic of slope $\alpha$ in the sense of \cite[Def.\,1.5.7]{HartlP}. The claim thus holds by \cite[Prop.\,1.5.10]{HartlP}. 
\end{Proof}


\medskip
The existence of the slope filtration has several consequences:

%

\begin{Prop}\label{LICMorphStrict}
Any morphism of $\ICF$-isocrystals $f\colon M\to N$ is \emph{strict} with respect to the slope filtration, that is, for each $\alpha$ we have $f(M_{\le\alpha}) = f(M)\cap N_{\le\alpha}$.
\end{Prop}

\begin{Proof}
On the one hand the slope filtration of $M$ induces a separated exhaustive right continuous ascending filtration of $f(M)$ by $\ICF$-isocrystals $f(M_{\le\alpha})$. For this filtration every subquotient $f(M_{\le\alpha})/f(M_{<\alpha})$ is isomorphic to a quotient of $M_\alpha=M_{\le\alpha}/M_{<\alpha}$ and therefore pure of slope~$\alpha$ by Proposition \ref{PureSubquot}. By Definition \ref{FiltDef} this is thus a slope filtration of~$f(M)$.

On the other hand the slope filtration of $N$ induces a separated exhaustive right continuous ascending filtration of $f(M)$ by $\ICF$-isocrystals $f(M)\cap N_{\le\alpha}$. For this filtration every subquotient ${(f(M)\cap N_{\le\alpha})/(f(M)\cap N_{<\alpha})}$ is isomorphic to a subobject of $N_\alpha=N_{\le\alpha}/N_{<\alpha}$ and therefore pure of slope~$\alpha$ by Proposition \ref{PureSubquot}. By Definition \ref{FiltDef} this is thus also a slope filtration of~$f(M)$.

By the unicity of slope filtration these filtrations therefore coincide, as desired.
\end{Proof}

\medskip
Next the functoriality of the slope filtration from Proposition \ref{FiltFunctUniq} \ref{FiltFunctUniqA} implies that for each $\alpha$ we obtain functors 
\UseTheoremCounterForNextEquation
\begin{equation}\label{MalphaFunctors}
\LIC{K} \to  \LIC{K},\ 
M \mapsto \left\{\kern-.33em\begin{array}{l}
M_{\le\alpha}, \\
M_{<\alpha}, \\
M_\alpha.
\end{array}\right.
\end{equation}

\begin{Prop}\label{LICFuncExact}
These functors are exact.
\end{Prop}

\begin{Proof}
For any short exact sequence of $\ICF$-isocrystals $0 \to M \to N \to P \to 0$, Proposition~\ref{LICMorphStrict} implies that the induced sequence $0 \to M_{\le\alpha} \to N_{\le\alpha} \to P_{\le\alpha} \to 0$ is exact for every~$\alpha$.
Thus the functor $M \mapsto M_{\le\alpha}$ is exact.
Since $M_{<\alpha} = M_{\le\alpha-\epsilon}$ and $N_{<\alpha} = N_{\le\alpha-\epsilon}$ and $P_{<\alpha} = P_{\le\alpha-\epsilon}$ for any sufficiently small $\epsilon>0$, it follows that the functor $M \mapsto M_{<\alpha}$ is exact as well. Together this implies that the functor $M \mapsto M_\alpha = M_{\le\alpha}/M_{<\alpha}$ is exact.
\end{Proof}

\medskip
\begin{samepage}
Next we discuss splittings. 

\begin{Prop}\label{PerfectFiltSplit}
\begin{enumerate}\StatementLabels
\item\label{PerfectFiltSplitUniq}
If the slope filtration of an $\ICF$-isocrystal splits, it splits uniquely.
\item\label{PerfectFiltSplitSplit}
The slope filtration splits whenever the field $K$ is perfect. 
In this case the category $\LIC{K}$ is therefore $\BQ$-graded.
\end{enumerate}
\end{Prop}

\end{samepage}

\begin{Proof}
Assertion \ref{PerfectFiltSplitUniq} follows by induction from Proposition \ref{ICSlopeProp} \ref{ICSlopePropA}.

In the situation of \ref{PerfectFiltSplitSplit} let $M$ correspond to the $\ICFalt$-isocrystal $M'$ as in Construction \ref{LICCoeffRes}. Then the slope filtration of $M'$ splits by \cite[Prop.\,1.5.10]{HartlP}. Moreover, since that splitting is unique by part \ref{PerfectFiltSplitUniq}, it is automatically $\ICF$-equivariant. As the slope filtration of $M'$ coincides with that of $M$ up to scaling, this yields a splitting of the latter, as desired.
\end{Proof}

\begin{Rem}\label{FiltNotSplitRem}\rm
If $K$ is not perfect, there exists an $\ICF$-isocrystal over~$K$ whose slope filtration splits over no finite extension of~$K$. Namely, choose any element $\xi\in K$ which is not a $q$-th power, and let $M$ be the $\ICF$-isocrystal over $K$ with the $\EFpK$-basis $m_0,m_1$ and the action 
$$\tau m_0=m_0 \qquad\hbox{and}\qquad 
\tau m_1=(\ICu\complot 1)\, m_1+(1\complot \xi)\, m_0.$$
Denoting by $N$ the sub-isocrystal generated by~$m_0$, the slope filtration of $M$ is given by
$$M_{\le\alpha} = \left\{\kern-.33em\begin{array}{cl}
0 &\hbox{if \ $\alpha < 0$,} \\
N &\hbox{if \ $0 \le \alpha < \ICd$,} \\
M &\hbox{if \ $\alpha\ge\ICd$.}
\end{array}\right.$$
For any splitting $s\colon M \to N$ of the inclusion $N\into M$ we must have $s(m_1) = m_1-\theta\, m_0$, where $\theta \in \EFpK$ is a solution of the equation
$$ \sigma(\theta) = (\ICu\complot 1)\, \theta+(1\complot \xi).$$
However this equation has the unique solution
$$\theta = \sum_{i\ge0} \ICu^i \complot \xi^{q^{-i-1}}$$
over an algebraic closure of~$K$. The smallest field over which the slope filtration splits is therefore the infinite purely inseparable extension $K\bigl(\xi^{q^{-i-1}}|_{i\ge0}\bigr)$ of~$K$.
\end{Rem}

\medskip
By combining Propositions \ref{StdPureICCons} \ref{StdPureICConsB} and \ref{PerfectFiltSplit} \ref{PerfectFiltSplitSplit} we obtain the following form of the Dieudonn\'e--Manin classification for local isocrystals:

\begin{Prop}\label{AlgClosed}
If $K$ is algebraically closed, every $\ICF$-isocrystal is isomorphic to a direct sum of copies of $\SA{K}$ for certain $\alpha\in\BQ$. In particular the category $\LIC{K}$ is semisimple.
\end{Prop}

\begin{Prop}\label{AlgClosedPullback}
For any embedding of algebraically closed fields $i\colon K\into K'$ the base change functor $i^*\colon \LIC{K} \isoto \LIC{K'}$ is an equivalence of tannakian categories.
\end{Prop}

\begin{Proof}
By Proposition \ref{AlgClosed} and \ref{StdPureICPullback} the functor is essentially surjective on objects, and by Propositions \ref{ICSlopeProp} \ref{ICSlopePropA} and \ref{StdPureICPullback} it is fully faithful. It also commutes with the tensor product by Construction \ref{LICBC} and is therefore an equivalence of tannakian categories.
\end{Proof}

\medskip
The slope filtration interacts with the tensor structure as follows:

\begin{Prop}\label{TensorSlope}
Consider any $\ICF$-isocrystals  $M$ and $N$ over~$K$.
\begin{enumerate}\StatementLabels
\item\label{TensorSlopeB}
The slope filtration of $M\otimes N$ is given for any $\gamma\in\BQ$ by
$$\kern10pt (M \otimes N)_{\le \gamma} \, =  \sum_{\smash{\alpha + \beta \le \gamma}} \! M_{\le\alpha}\otimes N_{\le\beta}.$$
\item\label{TensorSlopeC}
For any $\gamma\in\BQ$ there is a natural isomorphism
$$(M \otimes N)_\gamma \, \cong \bigoplus_{\smash{\alpha + \beta = \gamma}} \! M_\alpha\otimes N_\beta$$
\item\label{TensorSlopeD}
For any $\alpha\in\BQ$ there is a natural isomorphism
$$(M^\vee)_{-\alpha}\ \cong\ (M_\alpha)^\vee. \kern27pt$$
\end{enumerate}
\end{Prop}

\begin{Proof}
Since $M$ and $N$ are successive extensions of their subquotients $M_\alpha$ and~$N_\beta$, the exactness of the tensor product shows that $P := M\otimes N$ is a successive extension of $M_\alpha \otimes N_\beta$ for all $\alpha$ and~$\beta$. More precisely, for any $\gamma$ consider the subobjects
$$P_{[\le\gamma]}\ := \sum_{\alpha + \beta \le \gamma} \! M_{\le\alpha}\otimes N_{\le\beta} 
\qquad\text{and}\qquad
P_{[<\gamma]}\ :=\, \bigcup_{\gamma'<\gamma}P_{[\le\gamma']}.$$
Then the epimorphisms $M_{\le\alpha}\otimes N_{\le\beta} \onto M_\alpha\otimes N_\beta$ for all $\alpha+\beta=\gamma$ induce a natural isomorphism
\UseTheoremCounterForNextEquation
\begin{equation}\label{TensorSlopePGamma}
P_{[\gamma]}\ :=\ P_{[\le\gamma]}/P_{[<\gamma]} \ \cong \bigoplus_{\alpha + \beta = \gamma} \! M_\alpha\otimes N_\beta
\end{equation}
By Proposition \ref{TensorPure} \ref{TensorPureA} this is pure of slope $\alpha+\beta$. Varying $\gamma$ this shows that the subobjects $P_{[\le\gamma]}$ form a slope filtration of~$P$. By the uniqueness of the slope filtration this proves~\ref{TensorSlopeB}, and then \eqref{TensorSlopePGamma} yields \ref{TensorSlopeC}. Assertion \ref{TensorSlopeD} is established similarly.
\end{Proof}

\medskip
Propositions \ref{LICFuncExact} and  \ref{TensorSlope} directly imply: 

\begin{Cor}\label{LICGrFunctor}
We have an exact $\ICF$-linear tensor functor
$$ \LIC{K} \to  \LIC{K},\ M\mapsto \gr M := \bigoplus_\alpha M_\alpha.$$
\end{Cor}

For any $\ICF$-isocrystal $M$ we call $\rank(M_\alpha)$ the \emph{multiplicity of $\alpha$ in~$M$}, and we call the set of $\alpha$ with $M_\alpha\not=0$ the \emph{set of slopes of~$M$.} These numerical invariants of $M$ can be neatly encoded as follows:

\begin{Cons}\label{LICNewtonPolygon}\rm
The \emph{Newton polygon of~$M$} is the unique polygon in $\BR^2$ extending to the right from the point $(0,0)$ and with a line segment of slope $\alpha$ and horizontal displacement $\rank(M_\alpha)$ for every $\alpha\in\BQ$, arranged by increasing~$\alpha$. 
Its break points have integer coordinates, because by Propositions \ref{StdPureICCons} and \ref{PureCritHartl} the precise denominator of $\alpha$ always divides $\rank(M_\alpha)$.
\end{Cons}

\begin{Prop}\label{ExtPowSlope}
For any $\ICF$-isocrystal  $M$ over~$K$ and any integer $0\le n\le\rank(M)$, the smallest slope of $\bigwedge^nM$ is the sum of the $n$ smallest slopes of $M$, counted with multiplicities.
\end{Prop}

\begin{Proof}
We show this by induction on $\rank(M)$. For $n=0$ it follows from the fact that ${\bigwedge^0M\cong\unit}$ has slope zero. For $n>0$ suppose that $M$ has the smallest slope $\alpha$ with multiplicity~$\mu$. Then $M_{\le\alpha}=M_\alpha$ and we have a short exact sequence ${0\to M_\alpha\to M\to M/M_\alpha\to0}$. From multilinear algebra we deduce that $\bigwedge^nM$ a successive extension of the $\ICF$-isocrystals $\bigwedge^\nu M_\alpha \otimes \bigwedge^{n-\nu}(M/M_\alpha)$ for all $0\le \nu\le n$. Here $\bigwedge^\nu M_\alpha$ is non-zero and pure of slope $\nu\alpha$ by Proposition \ref{TensorPure} \ref{TensorPureC}. Moreover, if $\bigwedge^{n-\nu}(M/M_\alpha)$ is non-zero, the induction hypothesis asserts that its smallest slope is the sum of the $n-\nu$ smallest slopes of $M/M_\alpha$, counted with multiplicities. By construction all those slopes are $>\alpha$.

If $n\le\mu$, this shows that the smallest slope of $\bigwedge^nM$ is that of $\bigwedge^nM_\alpha \otimes \bigwedge^0(M/M_\alpha)$ and equal to $n\alpha$. If $n>\mu$, the smallest slope of $\bigwedge^nM$ is that of ${\bigwedge^\mu M_\alpha \otimes \bigwedge^{n-\mu}(M/M_\alpha)}$ and hence equal to $\mu\alpha$ plus the sum of the $n-\mu$ smallest slopes of $M/M_\alpha$, counted with multiplicities. In either case this is the sum of the $n$ smallest slopes of~$M$, counted with multiplicities, as desired.
\end{Proof}

\begin{Prop}\label{SlopesLowBound}
For any $\alpha$ and any $\ICF$-isocrystal  $M$ over~$K$ the following are equivalent:
\begin{enumerate}\StatementLabels
\item\label{SlopesLowBoundA}
All slopes of $M$ are $\ge\alpha$.
\item\label{SlopesLowBoundB}
There exist an integer $n>0$ with $n\alpha\in\BZ$ and a lattice $L\subset M$ such that ${\tau^{n\ICd}L\subset \ICu^{n\alpha} L}$. 
\end{enumerate}%
\end{Prop}

\begin{Proof}
Suppose first that all slopes of $M$ are $\ge\alpha$. 
Take an embedding $i\colon K \into K'$ into a perfect field. Then by Proposition \ref{PerfectFiltSplit} the $\ICF$-isocrystal $i^*M$ is isomorphic to the direct sum of its pure subquotients $(i^*M)_\beta$. Choose an integer $n>0$ with $n\alpha\in\BZ$ such that, for every $\beta$ with $(i^*M)_\beta\not=0$, we have $n\beta\in\BZ$ and a lattice $L'_\beta\subset(i^*M)_\beta$ such that $\smash{\tau^{n\ICd}L'_\beta}$ generates $\smash{\ICu^{n\beta} L'_\beta}$. Since all $\beta\ge\alpha$, the direct sum of these $L'_\beta$ yields a lattice $L'\subset i^*M$ with $\tau^{n\ICd}L' \subset \ICu^{n\alpha} L'$. By Lemma \ref{CapLattice} the intersection $L := L' \cap M$ is then a lattice in~$M$, which by construction satisfies $\tau^{n\ICd}L \subset \ICu^{n\alpha} L$, proving the implication \ref{SlopesLowBoundA}$\Rightarrow$\ref{SlopesLowBoundB}.

\medskip
Conversely consider an integer $n>0$ with $n\alpha\in\BZ$ and a lattice $L\subset M$ such that ${\tau^{n\ICd}L\subset \ICu^{n\alpha} L}$. Then for each $\beta$ we obtain a lattice $((L\cap M_{\le\beta})+M_{<\beta})/M_{<\beta}$ in $M_\beta$ with the corresponding properties. To deduce (a) it therefore suffices to consider the case that $M$ is non-zero and pure of slope~$\beta$.

In that case consider the $\ICF'$-isocrystal $M'$ corresponding to $M$ as in Construction \ref{LICCoeffRes}. If the implication  \ref{SlopesLowBoundB}$\Rightarrow$\ref{SlopesLowBoundA} holds for $M'$ in place of~$M$, using the same lattice $L$ we find that all slopes of $M'$ are $\ge\smash{\frac{\alpha}{\ICd}}$. Since $M'$ is pure of slope $\smash{\frac{\beta}{\ICd}}$ by Proposition \ref{LICCoeffResPure}, this implies that $\beta\ge\alpha$, and we obtain the desired implication for~$M$. For the rest of the proof we may therefore assume that $\ICF = \BFq\ppu$, in which case $\ICd=1$ and $\smash{\EAp{K}} = K\bbu$.

Now let $r$ be the rank of~$M$. Then $\bigwedge^rM$ is an $\ICF$-isocrystal of rank~$1$, and the exterior power $\bigwedge^rL$ over $K\bbu$ is a lattice in $\bigwedge^rM$ which satisfies $\tau^{n}\bigwedge^rL \subset \ICu^{nr\alpha} \bigwedge^rL$. Since $\bigwedge^rL$ is a free module of rank $1$ over $K\bbu$, its image $\tau^n\bigwedge^rL$ therefore generates a lattice of the form $\ICu^i\,\bigwedge^rL$ for some integer $i\ge nr\alpha$. This means that $\bigwedge^rM$ is pure of slope $\ge r\alpha$. Since $\bigwedge^rM$ is pure of slope $r\beta$ by Proposition \ref{TensorPure} \ref{TensorPureC}, we deduce that $\beta\ge\alpha$, finishing the proof of the implication \ref{SlopesLowBoundB}$\Rightarrow$\ref{SlopesLowBoundA}.
\end{Proof}

\section{An explicit equivalence of categories}
\label{EquivCat}

Keeping the situation of the preceding sections, in this section we work over a fixed algebraic closure $\Kalg$ of~$K$. We let $\kalg$ denote the algebraic closure of $\BFq$ in $\Kalg$ and write $i$ for the embedding $\kalg\into \Kalg$.
Then by Proposition \ref{AlgClosedPullback} we have an equivalence of tannakian categories
$$i^*\colon \LIC{\kalg}\ \longisoto\ \LIC{\Kalg}.$$
In this section we construct an explicit quasi-inverse.

\medskip
For this consider an $\ICF$-isocrystal $M$ over~$\Kalg$ and an $\ICF$-isocrystal $P$ over~$\kalg$. Then 
\UseTheoremCounterForNextEquation
\begin{equation}\label{HomPM}
\Hom(i^*P,M)\ \cong\ \Hom_{\EFp{\kalg}[\tau]}(P,M)
\end{equation}
is a finite dimensional $\ICF$-vector space and naturally a right module over~$\End(P)$. We define
\UseTheoremCounterForNextEquation
\begin{equation}\label{UPM}
U_P(M)\ :=\ \Hom(i^*P,M)\otend{P},
\end{equation}
which is an $\ICF$-isocrystal over~$\kalg$. 
Clearly this construction is functorial in~$M$. In particular it commutes with direct sums, and so the unique decomposition $M=\bigoplus_\alpha M_\alpha$ from Proposition \ref{PerfectFiltSplit} induces a natural isomorphism
\UseTheoremCounterForNextEquation
\begin{equation}\label{UPGrGrM}
U_P(M)\ \cong\ \bigoplus_\alpha U_P(M_\alpha).
\end{equation}
Next the unique decomposition $P=\bigoplus_\beta P_\beta$ induces a decomposition $i^*P=\bigoplus_\beta i^*P_\beta$. Here each summand $i^*P_\beta$ is pure of slope $\beta$ by Proposition \ref{LICBCPure}; so by Proposition \ref{ICSlopeProp} \ref{ICSlopePropA} we have $\Hom(i^*P_\beta,M_\alpha)=0$ whenever $\beta\not=\alpha$ and hence a natural isomorphism $\Hom(i^*P,M_\alpha)\cong \Hom(i^*P_\alpha,M_\alpha)$. Together we thus obtain a natural isomorphism
\UseTheoremCounterForNextEquation
\begin{equation}\label{UPGrM}
U_P(M)\ \cong\ \bigoplus_\alpha U_{P_\alpha}(M_\alpha).
\end{equation}
Next we vary~$P$:


\begin{Prop}\label{UPthetaCons}
For any morphisms of $\ICF$-isocrystals $P
\smash{\stackrel{\smash{\raisebox{-1pt}{$\scriptstyle\iota$}}}{\longinto}} P' 
\smash{\stackrel{\smash{\raisebox{-5pt}{$\scriptstyle\pi$}}}{\longonto}} P$ such that $\pi\circ\iota=\id$,
there exists a unique morphism $\theta_{\iota,\pi}$ given on pure tensors by the following formula:
$$\xymatrix@C-22pt@R-10pt{
U_{P}(M) \ar[d]^{\theta_{\iota,\pi}}
&=& \Hom(i^*P,M) \kern-0pt \mathop{\otimes}\limits_{\raisebox{-4pt}{$\scriptstyle\End(P)$}} \kern-0pt P
\ar@<29pt>[d] & \kern31pt h\otimes p \ar@{|->}@<16pt>[d] \\
U_{P'}(M) &=& 
\kern-0pt \Hom(i^*P',M)  \kern-4pt \mathop{\otimes}\limits_{\raisebox{-4pt}{$\scriptstyle\End(P')$}} \kern-4pt P' & \kern12pt (h\circ i^*\pi) \otimes \iota(p).\\}$$
\end{Prop}

\begin{Proof}
The morphisms $\iota$ and $\pi$ induce an $\ICF$-bilinear map
\UseTheoremCounterForNextEquation
\begin{equation}\label{UPthetaConsA}
\vcenter{\xymatrix@R-20pt{
\Hom(i^*P,M) \times P \ar[r] & \Hom(i^*P',M) \otend{P'},\\  
\kern58pt (h,p) \ar@{|->}[r] & (h\circ i^*\pi) \otimes \iota(p). \kern-36pt\\}}
\end{equation}
Since $\iota$ and $\pi$ induce a direct sum decomposition $P'\cong P\oplus \ker(\pi)$, for any endomorphism $e\in\End(P)$ there exists an $e'\in\End(P')$ such that $e'\circ\iota=\iota\circ e$ and $\pi\circ e'= e\circ\pi$. The computation
$$(h\circ i^*\pi) \otimes \iota(e(p)) 
 = (h\circ i^*\pi) \otimes e'(\iota(p)) 
 = (h\circ i^*\pi\circ i^*e') \otimes \iota(p) 
 = (h\circ i^*e\circ i^*\pi) \otimes \iota(p)$$
then shows that $(h,e(p))$ and $(h\circ i^*e,p)$ have the same image under \eqref{UPthetaConsA}. The bilinear map thus corresponds to a unique $\ICF$-linear map $\theta_{\iota,\pi}$.
This is also a morphism of isocrystals, so we are done.
\end{Proof}

\begin{Prop}\label{WPMDirSysProp}
Consider any morphisms $P
\smash{\stackrel{\smash{\raisebox{-1pt}{$\scriptstyle\iota$}}}{\longinto}} P' 
\smash{\stackrel{\smash{\raisebox{-5pt}{$\scriptstyle\pi$}}}{\longonto}} P$ such that $\pi\circ\iota=\id$.
\begin{enumerate}\StatementLabels
\item\label{WPMDirSysPropFunc}
The morphism $\theta_{\iota,\pi}$ is functorial in~$M$.
\item\label{WPMDirSysPropId}
For $\iota=\pi=\id_P$ we have $\theta_{\iota,\pi} = \id_{U_P(M)}$.
\item\label{WPMDirSysPropCompos}
For any other pair of morphisms $P'
\smash{\stackrel{\smash{\raisebox{-1pt}{$\scriptstyle\iota'$}}}{\longinto}} P'' 
\smash{\stackrel{\smash{\raisebox{-5pt}{$\scriptstyle\pi'$}}}{\longonto}} P'$ such that $\pi'\circ\iota'=\id$ we have 
$$\theta_{\iota',\pi'} \circ \theta_{\iota,\pi} = \theta_{\iota'\circ\iota,\pi\circ\pi'}.$$
\item\label{WPMDirSysPropMono}
The morphism $\theta_{\iota,\pi}$ is always a monomorphism. 
\item\label{WPMDirSysPropIso}
The morphism $\theta_{\iota,\pi}$ is an isomorphism if each slope of $M$ is also a slope of~$P$.
\item\label{WPMDirSysPropIndep}
The morphism $\theta_{\iota,\pi}$ depends only on $P$, $P'$, and $M$, but not on $\iota$ or~$\pi$.
\end{enumerate}
\end{Prop}

\begin{Proof}
Assertions \ref{WPMDirSysPropFunc} through \ref{WPMDirSysPropCompos} follow directly from the construction in Proposition~\ref{UPthetaCons}. In particular \ref{WPMDirSysPropFunc} implies that everything is compatible with direct sums in~$M$. To prove the remaining assertions we may therefore without loss of generality assume that $M$ is indecomposable.
By Propositions \ref{AlgClosed} and \ref{StdPureICPullback} we are then reduced to the case that $M=i^*S$ for some $S = \SA{\kalg}$. In particular $M$ is pure of slope~$\alpha$, so by \eqref{UPGrM} we may also assume that $P$ and $P'$ are pure of slope~$\alpha$. By Proposition \ref{StdPureICCons} \ref{StdPureICConsB} it thus suffices to consider the case $P=S^{\oplus n}$ and $P'=S^{\oplus n'}$ for some $n,n'\ge0$.

If $P=0$, then \ref{WPMDirSysPropMono} through \ref{WPMDirSysPropIndep} hold trivially. Otherwise the injectivity of $\iota$ implies that $1\le n\le n'$. By the full faithfulness part of Proposition \ref{AlgClosedPullback} we then have an isomorphism
$$\Hom(S^{\oplus n},S) \isoto \Hom(i^*S^{\oplus n},i^*S) = \Hom(i^*P,M)$$
of right modules over $\End(S^{\oplus n})=\End(P)$. Plugging this into \eqref{UPM} yields an isomorphism
$$U_P(M) \ \cong\ \Hom(S^{\oplus n},S)\otend{S^{\oplus n}}.$$ 
On the other hand the evaluation map yields an evident isomorphism
$$\eval\colon \Hom(S^{\oplus n},S)\otend{S^{\oplus n}} 
\ \longisoto\ S,\quad\ h\otimes p\ \longmapsto\ h(p).$$
Repeating this with $P'$ in place of $P$ we obtain the diagram
\UseTheoremCounterForNextEquation
\begin{equation}\label{WPMDirSysPropDiag}
\vcenter{\hbox{\xymatrix@C-22pt@R-10pt{
U_{P}(M)  \ar[d]^{\theta_{\iota,\pi}}
&\cong& \kern0pt \Hom(S^{\oplus n},S)\otend{S^{\oplus n}} 
\ar@<23pt>[d] \\ 
U_{P'}(M) &\cong& \kern3pt \Hom(S^{\oplus n'},S)\otend{S^{\oplus n'}} \kern2pt\\} 
\kern-4pt\raisebox{-5pt}{$\xymatrix@C+15pt@R-20pt{
\mathstrut \ar[dr]^-{\begin{turn}{-22}$\scriptstyle\eval$\end{turn}} & \\
& \ S \\
\mathstrut \ar[ur]_-{\begin{turn}{22}$\scriptstyle\eval$\end{turn}} & \\}$}}}
\end{equation}
Here the middle vertical arrow corresponding to $\theta_{\iota,\pi}$ is given by $h\otimes p\mapsto (h\circ\pi)\otimes\iota(p)$. Since $\pi\circ\iota=\id$, it follows that the diagram commutes.
As the two evaluation morphisms in \eqref{WPMDirSysPropDiag} are isomorphisms and independent of the choice of $\iota$ and~$\pi$, the same now follows for~$\theta_{\iota,\pi}$. This proves \ref{WPMDirSysPropMono} through \ref{WPMDirSysPropIndep} in the remaining case.
\end{Proof}

\medskip
Now we endow the collection of all $\ICF$-isocrystals $P$ over~$\kalg$ with the partial order for which $P\preccurlyeq P'$ if and only if there exist morphisms $P \smash{\stackrel{\smash{\raisebox{-1pt}{$\scriptstyle\iota$}}}{\longinto}} P' \smash{\stackrel{\smash{\raisebox{-5pt}{$\scriptstyle\pi$}}}{\longonto}} P$ such that $\pi\circ\iota=\id$. Any choice of such $\iota$ and $\pi$ then induces a morphism 
\UseTheoremCounterForNextEquation
\begin{equation}\label{ThetaPP'}
\theta_{P',P} := \theta_{\iota,\pi}\colon U_P(M) \longto U_{P'}(M)
\end{equation}
that is independent of $\iota$ and~$\pi$. Proposition \ref{WPMDirSysProp} implies that the morphisms $\theta_{P',P}$ form a filtered direct system for which all morphisms far enough out are isomorphisms. 
For any $\ICF$-isocrystal $M$ over~$\Kalg$ we can therefore define
\UseTheoremCounterForNextEquation
\begin{equation}\label{UMCons}
U(M)\ :=\ \dirlim_{P}U_P(M).
\end{equation}
Proposition \ref{WPMDirSysProp} \ref{WPMDirSysPropFunc} implies that this defines a covariant $\ICF$-linear functor 
\UseTheoremCounterForNextEquation
\begin{equation}\label{UFunctorCons}
U\colon \LIC{\Kalg}\longto\LIC{\kalg}.
\end{equation}

\begin{Rem}\label{UConsRem}\rm
This functor can be realized concretely as follows. For every $M$ choose an $\ICF$-isocrystal $P_M$ over $\kalg$ whose set of slopes contains the set of slopes of~$M$. Then for any $P'$ with $P_M\preccurlyeq P'$ the morphism $\theta_{P',P_M}$ is an isomorphism; hence we can set $U(M) := U_{P_M}(M)$, avoiding all set theoretic complications. To any morphism ${f\colon M\to M'}$ associate the morphism $U(f)\colon U(M)\to U(M')$ that makes the following diagram commute:
$$\xymatrix@C-12pt{
U(M) \ar@{}[r]|-{\textstyle=} \ar[d]_{U(f)} & U_{P_M}(M) \ar[rrrrr]_-\sim^{\theta_{P_M\oplus P_{M'},P_M}} &&&&& U_{P_M\oplus P_{M'}}(M) \ar@<7pt>[d]_{U_{P_M\oplus P_{M'}}(f)}\\
U(M') \ar@{}[r]|-{\textstyle=} & U_{P_{M'}}(M') \ar[rrrrr]_-\sim^{\theta_{P_M\oplus P_{M'},P_{M'}}} &&&&& U_{P_M\oplus P_{M'}}(M') \\
}$$
The properties of the direct system imply that this yields a concrete realization of the functor~$U$. 

However, expressing $U$ as a direct limit will make it easier to keep track of its properties and will avoid having to verify some ugly compatibility conditions, especially in relation to the tensor product. 
\end{Rem}


\begin{Prop}\label{UInverseOfi*}
There is a natural isomorphism $\rho_P \colon P\isoto U(i^*P)$ that is functorial in~$P$.
\end{Prop}

\begin{Proof}
The morphism is the composite
$$\xymatrix{
P\ \ar[rr]^-{p\mapsto\id\otimes p} && \ \Hom(i^*P,i^*P)\otend{P} 
\ar@{}[r]|-{\textstyle\stackrel{\smash{\rm def}}{=}} & U_P(i^*P) \ \ar[r] &\  U(i^*P).\\}$$
Clearly this is functorial in~$P$. Also the fact that $i^*$ is an equivalence of categories implies that $\End(P) \isoto \End(i^*P) = \Hom(i^*P,i^*P)$, so the map on the left hand side is an isomorphism. As the map on the right hand side is an isomorphism by Proposition \ref{WPMDirSysProp} \ref{WPMDirSysPropIso}, the composite morphism is an isomorphism. 
\end{Proof}

\begin{Thm}\label{UMThm}
The functor $U$ is an equivalence of categories that is quasi-inverse to~$i^*$.
\end{Thm}

\begin{Proof}
Direct consequence of Proposition \ref{UInverseOfi*} and the fact that $i^*$ is an equivalence of categories.
\end{Proof}

\medskip
To turn $U$ into a tensor functor, consider $\ICF$-isocrystals $M$ and $N$ over~$\Kalg$ and $\ICF$-isocrystals $P$ and $Q$ over~$\kalg$. Then we have a natural homomorphism 
\UseTheoremCounterForNextEquation
\begin{equation}\label{HomPMQN}
\Hom(i^*P,M) \otimes_{\ICF} \Hom(i^*Q,N) \longto \Hom\bigl(i^*(P\otimes Q),M\otimes N\bigr)
\end{equation}
that is equivariant under the natural ring homomorphism 
\UseTheoremCounterForNextEquation
\begin{equation}\label{HomPQ}
\End(P)\otimes_{\ICF}\End(Q)\longto \End(P\otimes Q).
\end{equation}
This yields a natural morphism of $\ICF$-isocrystals $\eta_{P,Q}$ making the following diagram commute:
\UseTheoremCounterForNextEquation
\begin{equation}\label{EndPQ}
\raisebox{-6pt}{$\vcenter{\xymatrix@C-22pt@R-0pt{
U_P(M) \otimes U_{Q}(N) \ar[d]_{\eta_{P,Q}}
&=& \mathstrut\smash{\bigl(\Hom(i^*P,M)\otend[-4]{P}\bigr) 
\kern2pt\mathop{\otimes}
\kern4pt
\bigl(\Hom(i^*Q,N)\kern-2pt\otend[-4]{Q}\bigr)}
\ar@<1pt>[d] \\
U_{P\otimes Q}(M\otimes N) &=& 
\Hom\bigl(i^*(P\otimes Q),M\otimes N\bigr) 
\kern-10pt \mathop{\otimes}\limits_{\raisebox{-4pt}{$\scriptstyle\End(P\otimes Q)$}} \kern-8pt
(P\otimes Q). \kern82pt \\}}$}\kern-30pt
\end{equation}

\begin{Lem}\label{UMTensorCompat}
For any $P\preccurlyeq P'$ and $Q\preccurlyeq Q'$ we have $P \otimes P' \preccurlyeq Q \otimes Q'$ and the following diagram commutes:
$$\xymatrix@C+20pt{
U_P(M) \otimes U_Q(N) \ar[r]^-{\eta_{P,Q}} \ar[d]_-{\theta_{P',P}\otimes\theta_{Q',Q}}
& U_{P\otimes Q}(M\otimes N)  \ar[d]_-{\theta_{P'\otimes Q', P\otimes Q}} \\
U_{P'}(M) \otimes U_{Q'}(N) \ar[r]^-{\eta_{P',Q'}} & U_{P'\otimes Q'}(M\otimes N) \\}$$
\end{Lem}

\begin{Proof} 
Choose morphisms $P
\smash{\stackrel{\smash{\raisebox{-1pt}{$\scriptstyle\iota$}}}{\longinto}} P' 
\smash{\stackrel{\smash{\raisebox{-5pt}{$\scriptstyle\pi$}}}{\longonto}} P$ and $Q
\smash{\stackrel{\smash{\raisebox{-1pt}{$\scriptstyle\iota'$}}}{\longinto}} Q'
\smash{\stackrel{\smash{\raisebox{-5pt}{$\scriptstyle\pi'$}}}{\longonto}} Q$ such that $\pi\circ\iota=\id$ and $\pi'\circ\iota'=\id$. 
Then the morphisms
$$\begin{tikzcd}[column sep=3.32em]
P \otimes P' \rar[hook,"\iota\otimes\iota'"] & Q \otimes Q' \rar[two heads,"\pi\otimes\pi'"] & P \otimes P'
\end{tikzcd}$$
also satisfy $(\pi\otimes\pi')\circ(\iota\otimes\iota')=\id$; hence we have $P \otimes P' \preccurlyeq Q \otimes Q'$. Moreover, the two instances of the homomorphism \eqref{HomPMQN} lie in a commutative diagram
$$\vcenter{\xymatrix@R-15pt{\scriptstyle h\otimes h' \ar@{|->}[d] \\ \scriptstyle (h\circ i^*\pi)\otimes(h'\circ i^*\pi')}} 
\kern-170pt
\vcenter{\xymatrix@C-0pt@R+15pt{
\Hom(i^*P,M) \otimes_{\ICF} \Hom(i^*Q,N) \ar[r] 
\ar[d]
&  \Hom\bigl(i^*(P\otimes Q),M\otimes N\bigr) \ar[d] \\
\Hom(i^*P',M) \otimes_{\ICF} \Hom(i^*Q',N) \ar[r] 
&  \Hom\bigl(i^*(P'\otimes Q'),M\otimes N\bigr) \\}}
\kern-130pt
\vcenter{\xymatrix@R-15pt{\scriptstyle h'' \ar@{|->}[d] \\ \scriptstyle h''\circ i^*(\pi\otimes\pi')}}$$
Combining this with the two instances of the diagram \eqref{EndPQ} and the three instances of Proposition \ref{UPthetaCons} yields a three dimensional commutative diagram from which the lemma follows.
\end{Proof}

\medskip
Lemma \ref{UMTensorCompat} implies that the $\eta_{P,Q}$ form a morphism of filtered direct systems with respect to $P$ and~$Q$. They therefore induce a morphism
\UseTheoremCounterForNextEquation
\begin{equation}\label{UMVNEtaCons}
U(M)\otimes U(N) \longto U(M\otimes N).
\end{equation}

\begin{Lem}\label{UMNTensorLem}
For any $P$ and $Q$ the following diagram commutes:
$$\xymatrix@C+20pt{
P\otimes Q \ar[r]^-{\rho_{P\otimes Q}}_-\sim \ar[d]_{\rho_P\otimes\rho_Q}^\wr & U(i^*(P\otimes Q)) \ar@{=}[d]^\wr \\
U(i^*P)\otimes U(i^*Q) \ar[r]^-{\eqref{UMVNEtaCons}} & U(i^*P\otimes i^*Q)\\}$$
\end{Lem}

\begin{Proof}
By Proposition \ref{UInverseOfi*} the morphisms $\rho_P\otimes\rho_Q$ and $\rho_{P\otimes Q}$ are the composites along the left, respectively right, vertical edge of the diagram
$$\kern118pt
\raisebox{-9pt}{$\xymatrix@C-3pt@R-15pt{\scriptstyle p\otimes q \ar@{|->}[d] 
&&& \scriptstyle p\otimes q \ar@{|->}[d] \\
\scriptstyle (\id\otimes p)\otimes(\id\otimes q)
& \scriptstyle (h\otimes p)\otimes(h'\otimes q) \ar@{|->}[r] 
& \scriptstyle (h\otimes h')\otimes(p\otimes q)
& \scriptstyle \id\otimes (p\otimes q)}$}
\kern-415pt
\xymatrix@C-57pt{
P\otimes Q \ar[d] \ar@{=}[rr] && P\otimes Q \ar@<0pt>[d] \\
\mathstrut\smash{\bigl(\Hom(i^*P,i^*P)\kern-8pt\mathop{\otimes}\limits_{\End(P)}\kern-8pt P\bigr)
\otimes \bigl(\Hom(i^*Q,i^*Q)\kern-8pt\mathop{\otimes}\limits_{\End(Q)}\kern-8pt Q\bigr)} \ar@{=}[d]^-{\rm def} \ar[rr] 
&& \Hom\bigl(i^*(P\otimes Q),i^*(P\otimes Q)\bigr)\kern-17pt\mathop{\otimes}\limits_{\End(P\otimes Q)} \kern-14pt(P\otimes Q) \kern70pt \ar@<0pt>@{=}[d]^-{\rm def} \\
\ U_P(i^*P) \otimes U_Q(i^*Q)\ \ar[d]^\wr \ar[r]^-{\eta_{P,Q}} &
\kern0pt\ U_{P\otimes Q}(i^*P\otimes i^*Q))\ \ar@{=}[r]^-\sim \ar@<0pt>[d]^\wr &\ U_{P\otimes Q}(i^*(P\otimes Q)) \ar@<0pt>[d]^\wr \\
\ U(i^*P) \otimes U(i^*Q)\ \ar[r]^-{\eqref{UMVNEtaCons}} & U(i^*P\otimes i^*Q)\
\ar@{=}[r]^-\sim &\ U(i^*(P\otimes Q))\ \kern0pt \\}$$
Here the upper part evidently commutes, the middle part commutes by \eqref{EndPQ}, and the lower part commutes by the construction of the morphism~\eqref{UMVNEtaCons}. 
\end{Proof}

\begin{Prop}\label{UMNTensorProp}
The morphism \eqref{UMVNEtaCons} is an isomorphism and functorial in $M$ and~$N$.
\end{Prop}

\begin{Proof}
Since $i^*$ is essentially surjective, it suffices to prove the isomorphy in the case $M=i^*P$ and $N=i^*Q$, where it follows at once from Lemma \ref{UMNTensorLem}. The functoriality in $M$ and $N$ results directly from that of the diagram~\eqref{EndPQ}.
\end{Proof}

\begin{Thm}\label{WMTensorThm}
The functor $U$ induces an equivalence of tannakian categories that is quasi-inverse to~$i^*$.
\end{Thm}

\begin{Proof}
By Lemma \ref{UMNTensorLem} the morphisms \eqref{UMVNEtaCons} extend $U$ to a tensor functor that is a left quasi-inverse of~$i^*$. Since $i^*$ is already an equivalence of tannakian categories, that tensor functor is also compatible with the associativity, commutativity, and unit constraints.
\end{Proof}

%

\medskip
Finally, composing $U$ with the functor ``\kern1pt forgetting~$\tau$'' we obtain an exact fiber functor from $\ICF$-isocrystals over $\Kalg$ to the rigid monoidal category of free modules of finite rank over $\EFp{\kalg}$:
\UseTheoremCounterForNextEquation
\begin{equation}\label{UFiberFunctor1}
\LIC{\Kalg} \stackrel{U} \longto \LIC{\kalg} \longto \mathop{\rm Free}({\EFp{\kalg}}).
\end{equation}

\medskip
The explicit construction of the tensor functor $U$ will enable us in the next section to generalize it to the case of an arbitrary field of definition $K$ and to endow it with a representation of the Weil group of~$K$.

\section{Weil representations}
\label{WeilRep}

Keeping the situation of the preceding sections, we now return to an arbitrary field $K$ over~$\BFq$. Let $k$ denote the algebraic closure of $\BFq$ within~$K$, so that $\kalg$ is an algebraic closure of~$k$. If $K^\sep$ denotes the separable closure of $K$ within~$\Kalg$, 
the restriction of automorphisms induces natural isomorphisms
\UseTheoremCounterForNextEquation
\begin{eqnarray}\label{GammaKDef}
\GK \kern5pt\ :=\ \kern3pt \Aut(\Kalg/K) \kern3pt & \longisoto & \Gal(\Ksep/K),\\[3pt]
\UseTheoremCounterForNextEquation
\label{IKDef}
\GKgeom\ :=\ \Aut(\Kalg/K\kalg) & \longisoto & \Gal(\Ksep/K\kalg).
\end{eqnarray}
These groups lie in a natural short exact sequence 
\UseTheoremCounterForNextEquation
\begin{equation}\label{GammaExactSeq}
\xymatrix{ 1 \ar[r] & \GKgeom \ar[r] & \GK \ar[r] & \Gk := \Gal(\kalg/k) \ar[r] & 1\rlap{.}\\}
\end{equation}
Next recall that $\varGamma_{\BFq} := \Gal(\bar k/\BFq)$ is the free profinite group that is topologically generated by the \emph{arithmetic Frobenius} automorphism $\Frob_q\colon \xi\mapsto\xi^q$ relative to~$\BFq$. For any $\gamma\in\GK$ there is therefore a unique element $\deg_q(\gamma)\in\smash{\hatBZ}$ such that $\gamma|\kalg = (\Frob_q)^{\deg_q(\gamma)}$. The \emph{Weil group of~$K$} is defined as the subgroup
\UseTheoremCounterForNextEquation
\begin{equation}\label{WeilGroupDef}
\WK\ :=\ \bigl\{\gamma\in\GK \bigm| \deg_q(\gamma)\in\BZ \bigr\}.
\end{equation}
If $k$ is infinite, then $\Gk$ is a closed subgroup of infinite index of $\varGamma_{\BFq}$; hence $\WK$ is equal to the geometric Galois group~$\GKgeom$. If $k$ is finite, then $\Gk$ is a closed subgroup of index $[k/\BFq]$ of $\varGamma_{\BFq}$ and we have a  natural short exact sequence 
\UseTheoremCounterForNextEquation
\begin{equation}\label{WeilExactSeq}
\xymatrix@C+5pt{ 1 \ar[r] & \GKgeom \ar[r] & \WK \ar[r]^-{\deg_q} & \BZ\cdot[k/\BFq] \ar[r] & 0\rlap{.}\\}
\end{equation}
In either case 
we endow $\WK$ with the topology for which $\GKgeom$ is an open subgroup carrying the given profinite topology.

\medskip
The aim of this section is to attach a natural continuous representation of $\WK$ to any $\ICF$-isocrystal $M$ over~$K$. To this end, for any auxiliary $\ICF$-isocrystal $P$ over~$\kalg$ we set
\UseTheoremCounterForNextEquation
\begin{equation}\label{HomPMW}
T_P(M)\ :=\ \Hom_{\EFp{\kalg}[\tau]}(P,M\complot_K\Kalg).
\end{equation}
Consider any homomorphism $f\in T_P(M)$ and any $\gamma\in \WK$. Since $\kalg$ is a perfect field, the map $\tau_P\colon P \to P$ is bijective; hence we can form the composite map
\UseTheoremCounterForNextEquation
\begin{equation}\label{HomPMWgamma}
{}^\gamma f\ :=\ (\id\complot\gamma)\circ f\circ \tau_P^{-\deg_q(\gamma)}
\colon\ P \longto M\complot_K\Kalg.
\end{equation}
The definition of $\deg_q(\gamma)$ implies that $(\id\complot\gamma)\circ\sigma^{-\deg_q(\gamma)}$ is the identity on $\kalg$; hence the map ${}^\gamma f$ is again $\smash{\EFp{\kalg}[\tau]}$-linear. The formula \eqref{HomPMWgamma} thus defines a left action of $\WK$ on $T_P(M)$. 

Since the group \smash{$\GKgeom$} acts trivially on~$P$ and continuously on $M\complot_K\Kalg$, it also acts continuously on the finite dimensional vector space $T_P(M)$ over the local field~$\ICF$. Therefore $\WK$ acts continuously on~$T_P(M)$. By construction this action also commutes with the right action of $\End(P)$.


Next recall that $i$ denotes the embedding $\kalg\into\Kalg$, and let $j$ be the embedding $K\into\Kalg$. As in \eqref{HomPM} we thus have a natural $\End(P)$-equivariant isomorphism
\UseTheoremCounterForNextEquation
\begin{equation}\label{AltHomPMW}
T_P(M)\ \cong \ \Hom(i^*P,j^*M).
\end{equation}
With the trivial action on~$P$ this therefore induces a continuous left action of $\WK$ on 
\UseTheoremCounterForNextEquation
\begin{equation}\label{UPMW}
U_P(M)\ :=\ U_P(j^*M)\ \cong\ T_P(M)\mathop{\otimes}\limits_{\End(P)}P.
\end{equation}
By construction the morphisms $\theta_{P',P}$ from \eqref{ThetaPP'} are $\WK$-equivariant. Thus $\WK$ acts compatibly on the whole direct system; hence we obtain a continuous left action of $\WK$ on 
\UseTheoremCounterForNextEquation
\begin{equation}\label{UMW}
U(M)\ :=\ \dirlim_{P}U_P(M).
\end{equation}
As in Section \ref{EquivCat} this construction is functorial in $M$ and compatible with the tensor product.
With Theorem \ref{WMTensorThm} we thus deduce:

\begin{Thm}\label{WMBigWeilRepThm}
The construction $M \mapsto U(M)$ constitutes an exact and faithful tensor functor from \smash{$\LIC{K}$} to the tannakian category of $\ICF$-isocrystals over \smash{$\kalg$} with a continuous action of $\WK$.
\end{Thm}

We can be more precise about the image of the Weil representation in several ways. First we have:

\begin{Prop}\label{WMBigWeilRepGrad}
For any $M$ there is a natural $\WK$-equivariant isomorphism
$$U(M)\ \cong\ \bigoplus_{\alpha\in\BQ} U(M_\alpha).$$
This is functorial in $M$ and compatible with the tensor product.
\end{Prop}

\begin{Proof}
The isomorphism is induced by the unique splitting of the slope filtration of~$j^*M$. By construction the resulting decomposition \eqref{UPGrM} is preserved by the action of~$\WK$. The compatibility with the tensor product follows from Corollary \ref{LICGrFunctor}.
\end{Proof}


\medskip
Next, for any $\alpha$ we fix a choice of $\SA{\kalg}$ and recall from Proposition \ref{StdPureICCons} \ref{StdPureICConsC} that 
\UseTheoremCounterForNextEquation
\begin{equation}\label{EndAlpha}
D_\alpha\ :=\ \End(\SA{\kalg})
\end{equation}
is a central division algebra over~$\ICF$ with Hasse invariant $[-\alpha] \in \BQ/\BZ$. Then Propositions \ref{LICBCPure}, \ref{StdPureICCons} and \ref{StdPureICPullback} imply that $j^*M_\alpha$ is a direct sum of $n := \rank(M_\alpha)/\rank(\SA{\kalg})$ copies of $i^*\SA{\kalg}$.
Thus 
\UseTheoremCounterForNextEquation
\begin{equation}\label{TWDef}
T_\alpha(M)\ :=\ T_{\SA{\kalg}}(M)
\, \stackrel{\eqref{UPGrM}}{=}\, T_{\SA{\kalg}}(M_\alpha)
\end{equation}
is a free right $D_\alpha$-module of rank~$n$. 
Choosing any basis, the left action of $\WK$ on $T_\alpha(M)$ therefore corresponds to a homomorphism from $\WK$ to the multiplicative group of the matrix algebra $\Mat_{n\times n}(D_\alpha)$.
On the other hand, combining \eqref{UPMW} and \eqref{UMW} with Proposition \ref{WPMDirSysProp} \ref{WPMDirSysPropIso} and \ref{WPMDirSysPropIndep} we obtain a natural $\WK$-equivariant isomorphism
\UseTheoremCounterForNextEquation
\begin{equation}\label{WMSmallWeilRepIsom}
U(M_\alpha)\ \cong\ T_\alpha(M)\centrot_{\raisebox{-2pt}{$\scriptstyle D_\alpha$}} \SA{\kalg}.
\end{equation}
Together this shows that the representation of $\WK$ on $U(M_\alpha)$ lands in the multiplicative group of a subalgebra isomorphic to $\Mat_{n\times n}(D_\alpha)$, which is a central simple $\ICF$-algebra with Hasse invariant $[-\alpha]$. 

\medskip
In the special case of an isocrystal $M$ that is pure of slope $\alpha$, the module $T_\alpha(M)$ is exactly the Tate module from \cite[\S6]{MornevT}.

\medskip
Finally we consider the case of slope~$0$. For the simple isocrystal $\SAa{\kalg}{0}$ we can then take the unit object $\unit = \EFp{\kalg}$ with $\tau=\sigma$. From \eqref{HomPMW} we obtain a natural isomorphism
\UseTheoremCounterForNextEquation
\begin{equation}\label{WMSlope0}
T_0(M)\ \cong\ (M\complot_K \Kalg)^\tau,
\end{equation}
which by \eqref{TWDef} is a vector space of dimension $\rank(M_0)$ over $D_0 = \End(\unit) = \ICF$. Here the right hand side is the usual Tate module of~$M$, or equivalently of its slope $0$ part~$M_0$, and carries a natural action of~$\GK$ via~$\Kalg$.
The isomorphism \eqref{WMSlope0} is equivariant under $\WK$ and functorial in $M$ and compatible with the tensor product. Also, by \eqref{WMSmallWeilRepIsom} in the case $\alpha=0$ we have a natural $\WK$-equivariant isomorphism
\UseTheoremCounterForNextEquation
\begin{equation}\label{WMSmallWeilRepIsomSlope0}
U(M_0)\ \cong\ T_0(M)\centrot_{\raisebox{-2pt}{$\scriptstyle \ICF$}} \unit.
\end{equation}

%
%
%
%
%
%
%
%
%
%

\section{Local isocrystals over a finite field}
\label{LocFinField}

In this section we assume that $k$ is finite. Then
\UseTheoremCounterForNextEquation
\begin{equation}\label{LocFinFieldRing}
\EFpk\ :=\ \ICF\complot_{\BFq} k\ =\ \ICF \otimes_{\BFq} k
\end{equation}
is a finite dimensional $\ICF$-algebra, and since $k$ is perfect, the structure endomorphism $\sigma\colon\EFpk\to\EFpk$ is a ring automorphism. 

\begin{Cons}\label{LocCharPol}\rm  
For any $\ICF$-isocrystal $M$ over $k$ it follows that the map $\tau_M\colon M \to M$ is bijective.
As the Frobenius endomorphism $\Frob_k\colon \xi\mapsto\xi^{|k|}$ is the identity on~$k$, the map $\tau_M^{[k/\BFq]}\colon M \isoto M$ is $k$-linear, and hence $\EFpk$-linear by \eqref{LocFinFieldRing}. Since $M$ is a free $\EFpk$-module,
\UseTheoremCounterForNextEquation
\begin{equation}\label{CharPolFormulaLocal}
\charact_{M}(X) \ :=\ \det\nolimits_{\EFpk}\bigl(X\cdot\id_M - \tau_M^{[k/\BFq]}\bigr)
\end{equation}
is therefore a monic polynomial of degree $\rank(M)$ with coefficients in $\EFpk$. On the other hand, in view of \eqref{LocFinFieldRing} it follows from \cite[Lemma~8.1.4]{BoecklePink} that
\UseTheoremCounterForNextEquation
\begin{equation}\label{BPLemmaLocal}
\charact_{M}(X^{[k/\BFq]})\ =\ \det\nolimits_{\ICF}\bigl(X\cdot\id_M - \tau_M\bigr).
\end{equation}
Thus $\charact_{M}$ actually has coefficients in~$\ICF$. We call it the \emph{characteristic polynomial of~$M$.}
\end{Cons}

Next let $\ICFalg$ be an algebraic closure of~$\ICF$. Since $\tau_M^{[k/\BFq]}$ is an automorphism, all roots of $\charact_M$ in $\ICFalg$ are non-zero. Let $\ord_\Fp$ denote the unique valuation on $\ICFalg$ that extends the normalized valuation on $\ICF$. 

\begin{Prop}\label{SlopesRoots}
For any $\alpha\in\BQ$, the number of roots $\lambda\in\ICFalg$ of $\charact_{M}$ that satisfy 
$$\ord_\Fp(\lambda)=\alpha\cdot\frac{[k/\BFq]}{\ICd},$$
counted with multiplicities, is equal to $\rank(M_\alpha)$.
\end{Prop}

\begin{Proof}
As the assertion is invariant under extensions, it suffices to treat the case that $M$ is pure of slope~$\alpha$. Then there exist an integer $n>0$ with $n\alpha\in\BZ$ and a lattice $L\subset M$ such that $\tau^{n\ICd}L=\ICu^{n\alpha} L$. Abbreviate $e:= [k/\BFq]$. Then for any root $\lambda\in\ICFalg$ of $\charact_{M}$, the definition \eqref{CharPolFormulaLocal} implies that $\ICu^{-n\alpha e} \lambda^{n\ICd}$ is an eigenvalue of the operator $\ICu^{-n\alpha e} (\tau_M^e)^{n\ICd} = (\ICu^{-n\alpha} \tau_{\smash{M}}^{n\ICd})^e$. But since $(\ICu^{-n\alpha} \tau_{\smash{M}}^{n\ICd})^e(L)=L$, this eigenvalue is a unit at~$\Fp$. Consequently we have $\ord_\Fp(\lambda) = \frac{n\alpha e}{n\ICd} = \frac{\alpha e}{\ICd}$, and since $M$ is pure of slope~$\alpha$, the proposition follows.
\end{Proof}

\begin{Cor}\label{CharPoly}%
\begin{enumerate}\StatementLabels%
\item\label{CharPolyPositive}
The slopes of $M$ are $\ge0$ if and only if $\charact_{M}$ has coefficients in~$\ICO$. 
\item\label{CharPolySlopesUnitRoots}
The $\ICF$-isocrystal $M$ is pure of slope $0$ if and only if $\charact_{M}$ has coefficients in $\ICO$ and its constant coefficient is a unit.
\end{enumerate}
\end{Cor}

Next, since $k$ is finite, its Weil group $\Wk$ is the free discrete group generated by the arithmetic Frobenius element $\Frob_k$. Its action on $U(M)$ is determined as follows:

\begin{Prop}\label{LFFUMW}
For any $\ICF$-isocrystal $M$ over $k$ there is a natural isomorphism
$$U(M)\ \cong\ M \complot_k \kalg.$$
Under this isomorphism, the action of $\Frob_k$ on $U(M)$ corresponds to the automorphism 
$$\tau_M^{-[k/\BFq]} \complot \id_{\kalg}.$$ 
\end{Prop}

\begin{Proof}
In the notation of \S\ref{WeilRep} we have $k = K$; hence $j$ is the embedding of $k$ into its chosen algebraic closure~$\kalg$. Thus $P:=j^*M$ satisfies the condition in Proposition \ref{WPMDirSysProp}~\ref{WPMDirSysPropIso}. Combining Proposition \ref{WPMDirSysProp}~\ref{WPMDirSysPropIndep} with \eqref{UPMW} and \eqref{UMW} therefore yields a natural $W_k$-equivariant isomorphism
\UseTheoremCounterForNextEquation
\begin{equation}\label{CharPolIsoms}
\xymatrix@C+10pt{U(M)\ \cong\ U_{j^*M}(M)\ =\ \End(j^*M) 
\kern-10pt \mathop{\otimes}\limits_{\raisebox{-2pt}{$\scriptstyle\End(j^*M)$}} \kern-8pt j^*M\ 
\ar[r]^-\eval_-\sim\ & \ j^*M\ =\ M\complot_k \kalg.}
\end{equation}
According to \eqref{HomPMWgamma} the element $\Frob_k\in W_k$ acts on $\End(j^*M)$ by the formula
$$f\ \longmapsto\ (\id\complot \Frob_k) \circ f\circ \tau_{j^*M}^{-[k/\BFq]}.$$
Here $\tau_{j^*M}$ acts on $j^*M = M\complot_k \kalg$ by $\tau_M\complot\Frob_q$; hence $\tau_{j^*M}^{-[k/\BFq]}$ acts by $\tau_M^{-[k/\BFq]}\complot\Frob_k^{-1}$. Therefore $\Frob_k$ sends the identity endomorphism in $\End(j^*M)$ to the endomorphism
$$(\id\complot \Frob_k) \circ \id\circ(\tau_M^{-[k/\BFq]}\complot\Frob_k^{-1})\ =\ \tau_M^{-[k/\BFq]}\complot\id.$$
Applying the evaluation isomorphism \eqref{CharPolIsoms} thus yields the desired formula.
\end{Proof}

\medskip
Combining \eqref{CharPolFormulaLocal} and Proposition~\ref{LFFUMW}, we obtain the following formula for the characteristic polynomial of the \emph{geometric Frobenius} element $\Frob_k^{-1} \in W_k$: 

\begin{Cor}\label{CharPolGalois1}
We have 
$$\charact_{M}(X) \ =\ \det\nolimits_{\EFp{\kalg}}\bigl(X\cdot\id - \Frob^{-1}_k \bigm| U(M) \bigr).$$
\end{Cor}

Moreover, with \eqref{WMSlope0} and \eqref{WMSmallWeilRepIsomSlope0} this implies:

\begin{Cor}\label{CharPolGalois2}
If $M$ is pure of slope $0$, we have
$$\charact_{M}(X) \ =\ \det\nolimits_{\ICF}\bigl(X \cdot \id - \Frob^{-1}_k \bigm | (M\complot_k \kalg)^{\tau}\bigr).$$
\end{Cor}

\section{Local isocrystals over a ring}
\label{GoodRedLoc}

Now consider an arbitrary commutative $\BFq$-algebra~$R$. As in Section \ref{LocIsoCrys} we equip $\ICF$ with its locally compact topology and $R$ with the discrete topology and consider the completed tensor product ring
\UseTheoremCounterForNextEquation
\begin{equation}\label{FPRDef}
\EFp{R} := \ICF\complot_{\BFq}\kern-1ptR,
\end{equation}
viewed as a difference ring via the partial Frobenius endomorphism $\sigma(a\complot\xi):= a\complot\xi^q$.

\begin{Def}\label{LICRDef}
A dualizable left $\EFp{R}[\tau]$-module is called an \emph{$\ICF$-isocrystal over $R$ (with the ground field~$\BFq$).}
\end{Def}

We denote the category of $\ICF$-isocrystals over $R$ by $\LIC{R}$. By \S\ref{DiffRingIC} this category carries a rigid symmetric monoidal structure. 

\begin{Prop}\label{LICRGoodRedRank}
If $\Spec R$ is connected, then any $\ICF$-isocrystal $\CM$ over $R$ is locally free of a unique finite rank over $\EFp{R}$.
\end{Prop}

\begin{Proof}
As in \eqref{FPKDefVar} we have 
$$\EFp{R}\ \cong\ \bigl(\ICk\otimes_{\BFq}\kern-2ptR\bigr)\ppu.$$
The assumption on $R$ implies that $\Spec(\ICk\otimes_{\BFq}\kern-1ptR)$ is non-empty with finitely many connected components that are transitively permuted by~$\sigma$. Thus the same follows for $\EFp{R}$. Since $\CM$ is a finitely generated projective $\EFp{R}$-module, the result follows.
\end{Proof}


\begin{Cons}\label{LICRBC}
(Base change) \rm  
Any ring homomorphism $f\colon R \to R'$ induces a morphism of difference rings $\id\complot f\colon \EFp{R} \to \EFp{R'}$ that we again abbreviate by~$f$. From Construction \ref{BC1} we thus obtain a natural $F_\Fp$-linear monoidal base change functor
$$f^*\colon\ \LIC{R} \longto \LIC{R'},\quad \CM\mapsto f^*\CM \cong \CM\complot_{R,f}\kern-1ptR'.$$
\end{Cons}

\begin{Prop}\label{LICRKFaith}
If $R$ is a noetherian integral domain and $R'$ is non-zero, the functor $f^*$ is faithful.
\end{Prop}

\begin{Proof}
We must show that the natural map $\Hom(\CM,\CN) \to \Hom(f^*\CM,f^*\CN)$ is injective for any $\ICF$-isocrystals $\CM$ and $\CN$ over~$R$. Using the formula \eqref{Homhomtau} for the inner hom, this reduces to showing that the map 
\UseTheoremCounterForNextEquation
\begin{equation}\label{MtauDeffstar}
\CM^\tau := \{m\in \CM\mid \tau m=m\} \longto (f^*\CM)^\tau := \{m\in f^*\CM\mid \tau m=m\}
\end{equation}
is injective for any $\ICF$-isocrystal $\CM$ over~$R$. This assertion does not change on replacing $\ICF$ by the subfield $\BFq\ppu$; hence we may assume that $\ICF=\BFq\ppu$.
Then the respective difference rings are $\EFp{R} = R\ppu = R\bbu[\ICu^{-1}]$ and $\EFp{R'}=R'\ppu$. In particular $\EFp{R}$ is again a noetherian integral domain, and the kernel $\Fa$ of $\EFp{R}\to\EFp{R'}$ is a non-zero ideal.

Now the map \eqref{MtauDeffstar} is obtained by taking $\tau$-invariants for the composite map
$$\CM\ \longto\ \CM/\Fa \CM\ \longinto\ \CM\complot_{R,f}\kern-1ptR'\ =\ f^*\CM.$$
Thus it suffices to show that $\CM^\tau\cap\Fa \CM=0$. But for any element $m$ therof and every $n\ge0$ we have $m=\tau^n m \in \tau^n(\Fa \CM) \subset \Fa^{q^n}\CM.$ In particular $m$ is contained in the intersection $\bigcap_{n\ge0}\Fa^{q^n}\CM$. Since $\CM$ is a finitely generated module over the noetherian integral domain $\EFp{R}$ and $\Fa$ is a proper ideal, this intersection is zero by Krull's intersection theorem, and we are done.
\end{Proof}

\medskip

Next we consider the following difference subring of $\EFp{R}$:
\UseTheoremCounterForNextEquation
\begin{equation}\label{APRDef}
\EAp{R} := \ICO\complot_{\BFq}\kern-1ptR.
\end{equation}
As in Section \ref{Slopes1}, by a \emph{lattice} in an $\ICF$-isocrystal $\CM$ over $R$ we mean a finitely generated projective $\EAp{R}$-submodule $\CL$ such that $\CL[\ICu^{-1}] = \CM$. Combining \cite[Lemmas 8.2.1, 8.1.4]{MornevT} we obtain the following innocuous looking lemma, whose proof is actually quite subtle: 

\begin{Lem}\label{DVRLattice}
If $i$ is the embedding of a discrete valuation ring $R$ into its quotient field~$K$, then for each lattice $L \subset i^*\CM$ the intersection $L\cap\CM$ is the unique lattice in $\CM$ that generates $L$ over $\EAp{K}$.
\end{Lem}


\medskip
We also note the following easy fact:

\begin{Lem}\label{IntDomUnit}
For any integral domain $R$ we have
$$R\ppu^\times\ =\ \ICu^\BZ \cdot R\bbu^\times.$$
\end{Lem}


\medskip
Now we can generalize a fundamental theorem of Watson \cite[p.\,15, Thm.\,4.1]{Watson}:

\begin{Thm}\label{LICRKFullFaith}
If $i$ is the embedding of a noetherian normal integral domain $R$ into its quotient field~$K$, the functor $i^*$ is fully faithful. 
\end{Thm}

\begin{Proof}
By Proposition \ref{LICRKFaith} the functor is already faithful. By the same reductions as in the proof of \ref{LICRKFaith} it suffices to show that $\CM^\tau \to (i^*\CM)^\tau$ is surjective for any $\ICF$-isocrystal $\CM$ over~$R$, and we may assume that $\EFp{R} = R\ppu = R\bbu[\ICu^{-1}]$ and $\EFp{K}=K\ppu$.

%

Now let $\Fq$ run through all prime ideals of height~$1$ of~$R$. Since $R$ is a noetherian normal integral domain, we then have $R=\bigcap_\Fq R_\Fq$ within $K$ by \cite[\S11.2, Cor.\,11.4]{Eisenbud}. This implies that $R\ppu = \bigcap_\Fq R_\Fq\ppu$ within $K\ppu$. Since $\CM$ is a direct summand of a free $R\ppu$-module, writing $i_\Fq$ for the embedding $R\into R_\Fq$ it follows that $\CM = \bigcap_\Fq i_\Fq^*\CM$ within $i^*\CM$, and hence $\CM^\tau = \bigcap_\Fq(i_\Fq^*\CM)^\tau$ within $(i^*\CM)^\tau$. After replacing $R$ by $R_\Fq$ we are therefore reduced to the case that $R$ is a discrete valuation ring.

Next choose an embedding $j\colon R\into R'$ into a complete discrete valuation ring with the same uniformizer $\ICu$ and a perfect residue field, and consider the resulting embedding $j_K\colon K\into K' := R'[\ICu^{-1}]$. Then we have $R=R'\cap K$ within $K'$ and hence $R\ppu = R'\ppu\cap K\ppu$ within $K'\ppu$. Since $\CM$ is a free $R\ppu$-module, it follows that $\CM = j^*\CM \cap i^*\CM$ within $j_K^*i^*\CM$, and hence $\CM^\tau = (j^*\CM)^\tau \cap (i^*\CM)^\tau$ within $(j_K^*i^*\CM)^\tau$. Replacing $R$ by $R'$ we are therefore reduced to the case that $R$ is a complete discrete valuation ring with perfect residue field~$k$. In other words we now have $R\cong k\bbX{x}$ 
and $\EFp{R} \cong R\ppu$. 

Under these assumptions, by the formula \eqref{Homhomtau} for the inner hom it suffices to prove the bijectivity of $\Hom(\CM,\CN)\to \Hom(i^*\CM,i^*\CN)$ for all $\CM$ and~$\CN$. Furthermore, this assertion does not change on rescaling the action of $\tau$ on $\CM$ and $\CN$ by the same power of~$\ICu$. Taking a sufficiently large power we may therefore assume that all slopes of $i^*\CM$ and $i^*\CN$ are non-negative. 

By Proposition \ref{SlopesLowBound} there then exist a lattice $L\subset i^*\CM$ and an integer $n>0$ such that ${\tau^{n} L \subset L}$. After replacing $L$ by the lattice generated by $\sum_{j=0}^{n-1}\tau^j L$ we may even assume that $\tau L\subset L$. By Lemma \ref{DVRLattice} we then obtain a lattice $\CL := L\cap\CM$ in $\CM$ that generates $L$ over $K\bbu$. Since $R\bbu$ is a local ring, this lattice is a finitely generated free $R\bbu$-submodule of some rank~$r$, and by construction it satisfies $\tau \CL\subset \CL$. Choose a basis $\ell_1,\ldots,\ell_r$ and write $\tau \ell_k=\sum_j f_{jk}\ell_j$ with $f_{jk}\in R\bbu$ for all $j$ and~$k$.
Since $\CL[\ICu^{-1}]=\CM$ and $\taulin_\CM$ is an isomorphism $\sigma^* \CM \isoto \CM$, the determinant $D:=\det(f_{jk})$ is then a unit in $R\ppu$. By Lemma \ref{IntDomUnit} it is therefore a power of $\ICu$ times a unit in $R\bbu$. 
In other words $\CL$ is an \emph{equicharacteristic $\tau$-module over $(R\bbu,\sigma)$} in the sense of Watson \cite[p.\,12, Def.\,3.7]{Watson}.

The same arguments applied to $\CN$ produce a lattice in $\CN$ with the corresponding properties. We are then precisely in the situation of \cite[p.\,15, Thm.\,4.1]{Watson} and obtain the desired full faithfulness.
\end{Proof}

\begin{Cor}\label{LICRGoodRedUnique}
Let $i$ be the embedding of a noetherian normal integral domain $R$ into its quotient field~$K$, and let $M$ be an $\ICF$-isocrystal over~$K$. If there is an $\ICF$-isocrystal $\CM$ over~$R$ with an isomorphism $i^* \CM \cong M$, then this is unique up to unique isomorphism.
\end{Cor}

\begin{Proof}
For any two such models $\CM$ over~$R$, apply Theorem \ref{LICRKFullFaith} in both directions.
\end{Proof}

\begin{Rem}\label{LICRGoodRedRem}\rm
In general an $\ICF$-isocrystal $\CM$ as in Corollary \ref{LICRGoodRedUnique} does not exist, 
even after replacing $R$ by $R[\xi^{-1}]$ for some non-zero element $\xi\in R$. The reason is that the structure data of $M$ involves Laurent series with infinitely many coefficients in~$K$, which may have larger and larger denominators.
\end{Rem}

\section{Purity and the slope filtration over a ring}
\label{RedSlopFilt}

Consider an $\ICF$-isocrystal $\CM$ over an arbitrary commutative $\BFq$-algebra~$R$. Recall that by a lattice in $\CM$ we mean a finitely generated projective $\EAp{R}$-submodule $\CL$ with $\CL[\ICu^{-1}] = \CM$. The definition \ref{ICSlopeDef} of purity over a field directly generalizes as follows, in agreement with \cite[Def.\,3.4.6]{MornevT}: 

\begin{Def}\label{LocICRSlopeDef}
The $\ICF$-isocrystal $\CM$ is called \emph{pure of slope~$\alpha\in\BQ$} if there exist an integer $n>0$ with $n\alpha\in\BZ$ and a lattice $\CL\subset \CM$ such that $\tau^{n\ICd}\CL$ generates the lattice $\ICu^{n\alpha} \CL$.
\end{Def}

The definition \ref{FiltDef} of the slope filtration generalizes in the same fashion:

\begin{Def}\label{LocICRSlopeFiltDef}
A \emph{slope filtration} of $\CM$ is a separated exhaustive right continuous ascending filtration by $\ICF$-isocrystals $\CM_{\le\alpha}$ over $R$ indexed by all $\alpha\in\BQ$, such that every subquotient
$$\CM_\alpha := \CM_{\le\alpha}/\CM_{<\alpha}
\qquad\hbox{with}\qquad
\CM_{<\alpha} := \smash{\bigcup_{\beta<\alpha}\CM_{\le\beta}}$$
is an $\ICF$-isocrystal that is pure of slope~$\alpha$. 
\end{Def}

In particular this requires that each subquotient is again a finitely generated projective $\EFp{R}$-module; or again that each $\CM_{\le\alpha}$ is a direct summand $\EFp{R}$-module of~$\CM$. The definitions directly imply that these notions are preserved under base change:

\begin{Prop}\label{LICRSlopeBC}
Consider any $\BFq$-algebra homomorphism $f\colon R \to R'$.
\begin{enumerate}\StatementLabels
\item\label{LICRSlopeBCa}
If $\CM$ is pure of slope $\alpha$, then so is $f^*\CM$.
\item\label{LICRSlopeBCb}
Any slope filtration $(\CM_{\le\alpha})_{\alpha\in\BQ}$ of~$\CM$ yields a slope filtration $(f^*\CM_{\le\alpha})_{\alpha\in\BQ}$ of~$f^*\CM$.
\end{enumerate}
\end{Prop}

\medskip
For the remainder of this section we assume that $R$ is a noetherian integral domain. Let $i\colon R\into K$ denote the embedding into its field of quotients. Set $X:=\Spec R$, and for each point $x\in X$ let $i_x\colon R\to k_x$ denote the natural homomorphism to the residue field at~$x$. Then $i^*\CM$ is an $\ICF$-isocrystal over the field~$K$, and we can view $\CM$ as a good model of $i^*\CM$ over~$R$. In this section we address some questions concerning the extension of the slope filtration from $i^*\CM$ to~$\CM$. 

\medskip
First we note:

\begin{Prop}\label{FiltFunctUniqR}
\begin{enumerate}\StatementLabels%
\item\label{FiltFunctUniqRA}
For any morphism of $\ICF$-isocrystals $f\colon \CM \to \CN$ over~$R$ admitting slope filtrations and every rational number $\alpha$ we have $f(\CM_{\le\alpha}) \subset \CN_{\le\alpha}$.
\item\label{FiltFunctUniqRB}
In particular the slope filtration is unique if it exists.
\end{enumerate}
\end{Prop}

\begin{Proof}
Since $R$ is an integral domain, combining Proposition \ref{LICRSlopeBC} \ref{LICRSlopeBCb} with the fact that each $\CM_{\le\alpha}$ is a direct summand $\EFp{R}$-module of~$\CM$ implies that $\CM_{\le\alpha} = \CM\cap (i^*\CM)_{\le\alpha}$ for every~$\alpha$. As the slope filtration of $i^*\CM$ is unique by Proposition \ref{FiltFunctUniq} \ref{FiltFunctUniqB}, this proves \ref{FiltFunctUniqRB}.
Similarly, combining the formula for $\CM_{\le\alpha}$ with Proposition \ref{FiltFunctUniq} \ref{FiltFunctUniqA} proves \ref{FiltFunctUniqRA}.
\end{Proof}

\medskip
Next, by applying Proposition \ref{LICRGoodRedRank} to all subobjects $\CM_{\le\alpha}$, the existence of the slope filtration implies that the $\ICF$-isocrystals $i_x^*\CM$ over~$k_x$ have the same Newton polygon for all points $x\in X$. As the Newton polygon changes under reduction in the following example, we conclude that in general a slope filtration does not exist.

\begin{Ex}\label{LICRGoodRedEx1}\rm 
Take $F_\Fp := \BFq\ppu$ and $R := \BFq[\xi]$ with a variable~$\xi$, so that $K=\BFq(\xi)$ and $\EFp{R} = R\ppu$ and $\EFp{K} = K\ppu$. Let $\CM$ be the free $R\ppu$-module with basis $m_0,m_1$ and set $\tau m_0 = m_1$ and $\tau m_1=\ICu m_0 + \xi m_1$. This defines an $\ICF$-isocrystal of rank~$2$ over~$R$.
(In fact it arises by the constructions in Sections \ref{AMot} and \ref{Drin} 
from the Drinfeld $\BFq[t]$-module $\phi$ of rank $2$ and height $1$ over $R$ with $\phi_t = \tau^2-\xi\tau$ and $\ICu=t$.)

We claim that $i^*\CM$ has the slopes $0$ and $1$ with multiplicity $1$ each. To show this we look for an element of the form $m=\theta m_0+m_1\in i^*\CM$ with $\theta\in K\ppu$ that generates a sub-isocrystal of slope~$0$. This means that $m=\psi\cdot\tau m$ for a unit $\psi\in K\bbu^\times$. Plugging in the choice of $\tau$ and extracting coefficients this amounts to the equations
$$\theta=\psi\ICu \quad\hbox{and}\quad 1 = \psi\cdot(\sigma(\theta)+\xi).$$
Using the first equation to eliminate~$\theta$, this becomes equivalent to
$$\psi\cdot(\xi+\sigma(\psi)\cdot\ICu)=1.$$
Writing $\psi = \sum_{i\ge0} \psi_i\ICu^i$ with $\psi_i\in K$ and expanding this amounts to the recursion equations
$$\psi_0\xi=1
\quad\hbox{and}\quad \psi_\ell \xi + \kern-7pt \sum_{i+j=\ell-1} \kern-5pt \psi_i\psi_j^q = 0
\quad\hbox{for all}\ \ell\ge1.$$
By iteration these equations have a unique solution in $K$ (though not in~$R$). Thus $m$ generates a sub-isocrystal of slope $0$ of~$i^*\CM$, and an easy computation shows that the quotient isocrystal is pure of slope~$1$. This therefore constitutes the slope filtration of~$i^*\CM$.

On the other hand we have $\tau^2 m_0\equiv \xi m_0$ and $\tau^2 m_1\equiv \xi m_1$ modulo~$(\xi)$; hence the induced $\ICF$-icocrystal over the residue field $R/(\xi)$ is pure of slope $\sfrac{1}{2}$. Thus the Newton polygon changes on reduction; hence a slope filtration of $\CM$ does not exist.
\end{Ex}

\medskip
To prove the next theorem we will need some technical facts. Let $z$ be a non-zero element of the maximal ideal of~$A_\Fp$, not necessarily a uniformizer, and consider the difference subring $\BFq\bbX{z}\complot_{\BFq}\kern-1ptR \cong R\bbX{z}$ of $A_\Fp\complot_{\BFq}\kern-1ptR$.

\begin{Lem}\label{LatticeFreeLem}
Let $\CL$ be a module over $A_\Fp\complot_{\BFq}\kern-1ptR$ that is finitely generated projective over the subring $\BFq\bbX{z}\complot_{\BFq}\kern-1ptR \cong R\bbX{z}$. Then $\CL$ is finitely generated projective over $A_\Fp\complot_{\BFq}\kern-1ptR$.
\end{Lem}

\begin{Proof}
First observe that, since $A_\Fp\complot_{\BFq}\kern-1ptR$ is the completion of a ring that is finitely generated over the noetherian ring~$R$, it is again noetherian. On the other hand, since $\CL$ is finitely generated over $R\bbX{z}$, it is also finitely generated over $A_\Fp\complot_{\BFq}\kern-1ptR$. It therefore suffices to prove that $\CL$ is flat over $A_\Fp\complot_{\BFq}\kern-1ptR$.

For this assume first that $R$ is a field~$k$. Then $k\bbX{z}$ is a discrete valuation ring, and $A_\Fp\complot_{\BFq}\kern-1ptk$ is a finite product of discrete valuation rings that are finite extensions of $k\bbX{z}$. Since $\CL$ is projective over $k\bbX{z}$, it is torsion free, and hence it is a finite product of torsion free modules over the simple direct factors of $A_\Fp\complot_{\BFq}\kern-1ptk$. These modules are therefore flat, and so $\CL$ is flat, as desired. 

In the general case consider any maximal ideal $\Fm\subset R$ with residue field~$k$. Then $\CL/\Fm\CL$ is a module over the factor ring $A_\Fp\complot_{\BFq}\kern-1ptk$ which is finitely generated projective over $k\bbX{z}$. By the preceding case it is therefore flat over $A_{\Fp}\complot_{\BFq}\kern-1ptk$.
By the fiberwise criterion of flatness \cite[Lemme 11.3.10.1]{EGA4} it follows that the localizations of $\CL$ at all maximal ideals of $A_{\Fp}\complot_{\BFq}\kern-1pt R$ above $\Fm$ are flat. Since every maximal ideal of $A_\Fp\complot_{\BFq}\kern-1ptR$ lies above a maximal ideal of~$R$, varying $\Fm$ this means that $\CL$ is flat over $A_{\Fp}\complot_{\BFq}\kern-1ptR$, as desired.
\end{Proof}


\begin{Lem}\label{LocICRGroucho}
Let $\CL$ be a free $R\bbX{z}$-module of finite rank with a $\sigma^n$-linear map ${\rho\colon\CL\to\CL}$ for some integer $n\ge1$. Let $\bar\CL'$ be an $R$-submodule of $\bar\CL := \CL/z\CL$ that is generated by the image of a part of a basis of~$\CL$. Assume that $\bar\CL'$ contains $\rho\bar\CL$ and is generated by~$\rho\bar\CL'$. Then there exists a unique free $R\bbX{z}$-submodule $\CL' \subset \CL$ of rank $\rank_R(\bar\CL')$ such that ${(\CL'+z\CL)/z\CL}=\bar\CL'$ and $\rho\CL'\subset\CL'$. Moreover this $\CL'$ is a direct summand of~$\CL$.
\end{Lem}

\begin{Proof}
Choose an $R\bbX{z}$-basis $\ell_1,\ldots,\ell_r$ of $\CL$ such that the residue classes of $\ell_1,\ldots,\ell_s$ form an $R$-basis of $\bar\CL'$ for some $s\le r$. Then $\CL'$ must possess an $R\bbX{z}$-basis $\ell_1',\ldots,\ell_s'$ of length~$s$. Express each $\ell_j'$ in the form
$$\ell_j'\ =\ \sum_{i=1}^s\omega_{i,j}\ell_i+\sum_{i=s+1}^r\theta_{i,j}\ell_i$$
with $\omega_{i,j},\theta_{i,j}\in R\bbX{z}$. Then the condition $(\CL' + z\CL)/z\CL = \bar\CL_0$ means that the $s\times s$-matrix $\Omega := (\omega_{i,j})_{1\le j\le s}^{1\le i\le s}$ is invertible modulo~$(z)$ and that all $\theta_{i,j}$ lie in $R\bbX{z}z$. The first of these statements implies that $\Omega$ is invertible over $R\bbX{z}$. After replacing $\ell_1',\ldots,\ell_s'$ by another $R\bbX{z}$-basis of~$\CL'$ we can thus assume that $\Omega$ is the identity matrix and hence
\UseTheoremCounterForNextEquation
\begin{equation}\label{LocICRGroucho1}
\ell_j'\ =\ \ell_j+\sum_{i=s+1}^r\theta_{i,j}\ell_i
\end{equation}
with $\theta_{i,j}\in R\bbX{z}z$. Then $\ell_1',\ldots,\ell_s',\ell_{s+1},\ldots,\ell_r$ is an $R\bbX{z}$-basis of~$\CL$, and so $\CL'$ is a direct summand of~$\CL$. 

Next we have $\rho\CL'\subset\CL'$ if and only if for all $k$ we have
\UseTheoremCounterForNextEquation
\begin{equation}\label{LocICRGroucho2}
\rho \ell_k'\ =\ \sum_{j=1}^s\psi_{j,k}\ell_j'
\end{equation}
for elements $\psi_{j,k}\in R\bbX{z}$. It therefore remains to show that there exist unique $\theta_{i,j}$ and $\psi_{j,k}$ satisfying \eqref{LocICRGroucho1} and \eqref{LocICRGroucho2}.

To express these equations in terms of matrices we consider the $(r-s)\times s$-matrix $\Theta := (\theta_{i,j})_{1\le j\le s}^{s<i\le r}$ with coefficients in $R\bbX{z}z$ and the $s\times s$-matrix $\Psi = (\psi_{j,k})_{1\le k\le s}^{1\le j\le s}$ with coefficients in $R\bbX{z}$. Also, for each $1\le j\le r$ we write 
\UseTheoremCounterForNextEquation
\begin{equation}\label{LocICRGroucho3}
\rho \ell_j\ =\ \sum_{i=1}^r \phi_{i,j}\ell_i
\end{equation}
with elements $\phi_{i,j}\in R\bbX{z}$. We write the $r\times r$-matrix $\Phi := (\phi_{i,j})_{1\le j\le r}^{1\le i\le r}$ as a block matrix $\bigl(\begin{smallmatrix}\Phi_{11}&\Phi_{12}\\ \Phi_{21}&\Phi_{22}\end{smallmatrix}\bigr)$, where $\Phi_{11}$ has size $s\times s$ and the other blocks have the respective complementary sizes. The assumption on $\bar\CL'$ then means that $\Phi_{21}$ and $\Phi_{22}$ have coefficients in $R\bbX{z}\ICu$ and that $\Phi_{11}$ lies in $\GL_s(R\bbX{z})$.

By a direct computation that we leave to the careful reader, the equations \eqref{LocICRGroucho1} and \eqref{LocICRGroucho2} are then equivalent to the matrix equations
\UseTheoremCounterForNextEquation
\begin{eqnarray}
\label{LocICRGroucho4}
\Phi_{11} + \Phi_{12}\cdot\sigma^n(\Theta) &=& \Psi, \\[3pt]
\UseTheoremCounterForNextEquation
\label{LocICRGroucho5}
\Phi_{21} + \Phi_{22}\cdot\sigma^n(\Theta) &=& \Theta\cdot \Psi.
\end{eqnarray}
Using the first of these to eliminate~$\Psi$, and employing the fact that $\Phi_{11}$ lies in $\GL_s(R\bbX{z})$, these equations are equivalent to
\UseTheoremCounterForNextEquation
\begin{equation}\label{LocICRGroucho6}
\Theta\ =\ \Phi_{21}\cdot\Phi_{11}^{-1} + \Phi_{22}\cdot\sigma^n(\Theta)\cdot\Phi_{11}^{-1} - \Theta\cdot\Phi_{12}\cdot\sigma^n(\Theta)\cdot\Phi_{11}^{-1}.
\end{equation}

To solve this equation we observe that for any matrix $\Theta$ with coefficients in $R\bbX{z}z$, both $\Theta$ and $\sigma^n(\Theta)$ are divisible by~$z$, and by assumption $\Phi_{22}$ is also divisible by~$z$. For any $m\ge0$ it follows that the right hand side of \eqref{LocICRGroucho6} modulo $R\bbX{z}z^{m+1}$ depends only on $\Theta$ modulo $R\bbX{z}z^m$. Moreover, since $\Phi_{21}$ is divisible by~$z$, this resulting matrix again has coefficients in $R\bbX{z}z$. Thus if $\Theta$ solves the equation modulo $R\bbX{z}z^m$, the right hand side of \eqref{LocICRGroucho6} yields a solution modulo $R\bbX{z}z^{m+1}$. Moreover, if the solution modulo $R\bbX{z}z^m$ is unique, the same argument shows that the resulting solution modulo $R\bbX{z}z^{m+1}$ is unique. By induction on~$m$ we thus deduce that there exists a unique solution modulo $R\bbX{z}z^m$ for every~$m$. In the limit we therefore obtain a unique solution with coefficients in $R\bbX{z}z$.

This finishes the proof of the lemma.
\end{Proof}


\pagebreak
\begin{Thm}\label{LocICRGoodRedThm}
There exists an element $\xi\in R\setminus\{0\}$ such that $\CM$ acquires a slope filtration over $R[\xi^{-1}]$.
\end{Thm}

\begin{Proof}
For better readability we abbreviate the $\ICF$-isocrystal over $K$ obtained from $\CM$ by ${M := i^*\CM}$.
Without loss of generality we may assume that $\CM$ and hence $M$ is non-zero. Let $\alpha$ be the smallest slope of~$M$. Then the first non-trivial step in the slope filtration of $M$ is $M_{\le\alpha} = M_\alpha$. We will prove that, after replacing $R$ by $R[\xi^{-1}]$ and $\CM$ by its associated base change for some $\xi\in R\setminus\{0\}$, the $\ICF$-sub-isocrystal $M_\alpha$ extends to an $\ICF$-sub-isocrystal $\CM_\alpha$ of $\CM$ over $R$ that is pure of slope~$\alpha$, such that $\CM/\CM_\alpha$ is an $\ICF$-isocrystal.

\begin{Lem}\label{LocICRGoodRedLem2}
There exist an integer $n>0$ with $n\alpha\in\BZ$ and a lattice $L\subset M$ such that, with the operator $\rho := \ICu^{-n\alpha}\tau^{nd_\Fp}$ and the lattice $L_\alpha := L\cap M_\alpha$ in~$M_\alpha$, we have $\rho L \subset L_\alpha+\ICu L$, and $\rho L_\alpha$ generates $L_\alpha$.
\end{Lem}

\begin{Proof}
Since $M_\alpha$ is pure of slope $\alpha$, there exist an integer $n>0$ with $n\alpha\in\BZ$ and a lattice $L_\alpha\subset M_\alpha$ such that $\tau^{nd_\Fp}L_\alpha$ generates $\ICu^{n\alpha} L_\alpha$. Moreover, since all slopes of $M/M_\alpha$ are $>\alpha$, by Proposition \ref{SlopesLowBound} there exist an integer $n>0$ and a lattice $\bar L\subset M/M_\alpha$ such that $\tau^{nd_\Fp} \bar L \subset \ICu^{n\alpha+1}\bar L$. A common multiple of the two instances of $n$ does the job in both statements, and then $\rho L_\alpha$ generates $L_\alpha$ and $\rho \bar L \subset \ICu\bar L$. 
Since $M/M_\alpha$ is a projective $\EFpK$-module,  the projection $M\onto M/M_\alpha$ possesses a $\EFpK$-linear splitting $s\colon M/M_\alpha \into M$. Then there exists an integer $m$ with $\rho s(\bar L) \subset {\ICu^{-m}L_\alpha + \ICu s(\bar L)}$. Thus $L := L_\alpha + \ICu^m s(\bar L)$ is a lattice in $M$ which satisfies $\rho L \subset L_\alpha+\ICu L$ and $L\cap M_\alpha = L_\alpha$, as desired.
\end{Proof}

\medskip
For the next arguments recall that $\ICF=\ICk\ppu$ is a finite extension of $\BFq\ppu$. We are interested in the following subrings:
$$\xymatrix@R-7pt@C-10pt{
\EAp{R} \ar@{}[r]|-{\textstyle\supset} \ar@{}[d]|{\textstyle\subsetdown} & 
R\bbu \ar@{}[r]|-{\textstyle\subset} \ar@{}[d]|{\textstyle\subsetdown} & 
R\ppu \ar@{}[d]|{\textstyle\subsetdown} & 
\EFp{R} \ar@{}[l]|-{\textstyle\subset} \ar@{}[d]|{\textstyle\subsetdown} \\
\EAp{K} \ar@{}[r]|-{\textstyle\supset} & 
K\bbu \ar@{}[r]|-{\textstyle\subset} & 
K\ppu & 
\EFp{K}\rlap{.} \ar@{}[l]|-{\textstyle\subset} }$$
Since $L$ and $L_\alpha$ are finitely generated projective modules over $\EAp{K}$, they are a free modules of finite rank over the subring $K\bbu$. Moreover, by construction $L_\alpha$ is a direct summand of~$L$. Consider the $K$-vector space $\bar L := L/\ICu L$ and its subspace $\bar L_\alpha := {(L_\alpha+\ICu L)/\ICu L}$.

\begin{Lem}\label{LocICRGoodRedLem3}
There exist integers $r\ge s\ge1$ and elements $m_1,\ldots,m_r\in\CM$ such that
\begin{enumerate}\StatementLabels
\item\label{LocICRGoodRedLem3A}
The elements $m_1,\ldots,m_r$ form a basis of $L$ over $K\bbu$.
\item\label{LocICRGoodRedLem3C}
The residue classes $[m_1],\ldots,[m_s]$ form a basis of $\bar L_\alpha$ over~$K$.
\end{enumerate}
\end{Lem}

\begin{Proof}
Take a basis $\ell_1,\ldots,\ell_s$ of $L_\alpha$ over $K\bbu$ and extend it to a basis $\ell_1,\ldots,\ell_r$ of $L$ over $K\bbu$. These elements then also form a basis of $M$ over $K\ppu$.
Next observe that
the subring $K\otimes_RR\ppu$ is dense in $K\ppu$ for the $(\ICu)$-adic topology. Thus $K\otimes_R\CM$ is dense in $M$. Since $\ICu L$ is open in $M$, we can choose elements $m_1,\ldots,m_n\in K\otimes_R\CM$ such that $m_i\equiv\ell_i$ modulo $\ICu L$. These elements then again form a basis of $L$ over $K\bbu$. By construction they also satisfy the property in \ref{LocICRGoodRedLem3C}. Finally, after multiplying the $m_i$ by suitable non-zero elements of~$R$ we can achieve that they lie in~$\CM$, and are done.
\end{Proof}

\begin{Lem}\label{LocICRGoodRedLem4a}
After replacing $R$ by $R[\xi^{-1}]$ and $\CM$ by its base change $\CM[\xi^{-1}]$ for some $\xi\in R\setminus\{0\}$, the above elements $m_1,\ldots,m_r$ form a basis of $\CM$ over $R\ppu$.
\end{Lem}

\begin{Proof}
By the condition \ref{LocICRGoodRedLem3} \ref{LocICRGoodRedLem3A} the elements $m_1,\ldots,m_r$ are $R\ppu$-linearly independent. Since $\CM$ is an $R\ppu$-module of rank~$r$, every element of it is therefore a $\Quot(R\ppu)$-linear combination of $m_1,\ldots,m_r$. Since $\CM$ is finitely generated over $R\ppu$, it follows that there exists an element $\theta\in R\ppu\setminus\{0\}$ such that $\theta\CM$ is contained in the $R\ppu$-submodule generated by $m_1,\ldots,m_r$. Let $\xi$ be the lowest non-zero coefficient of~$\theta$. Then $\theta$ becomes a unit in $R[\xi^{-1}]\ppu$. Thus after replacing $R$ by $R[\xi^{-1}]$ and $\CM$ by its associated base change, the elements $m_1,\ldots,m_r$ generate~$\CM$ and therefore form a basis.
\end{Proof}

\medskip
Since $m_1,\ldots,m_r$ is a basis of $L$ over $K\bbu$ and a basis of $\CM$ over $R\ppu$, it is also a basis of the free $R\bbu$-submodule $\CL := L\cap\CM$. Since $\rho$ maps each of $L$ and $\CM$ into itself, it also maps $\CL$ to itself. Consider the $R$-submodule $\bar\CL' \subset \bar\CL := \CL/\ICu\CL$ that is generated by the residue classes $[m_1],\ldots,[m_s]$.

\begin{Lem}\label{LocICRGoodRedLem4b}
After replacing $R$ by $R[f^{-1}]$ and $\CM$ by its base change $\CM[f^{-1}]$ for some $f\in R\setminus\{0\}$, the submodule $\bar\CL'$ contains $\rho\bar\CL$ and is generated by~$\rho\bar\CL'$.
\end{Lem}

\begin{Proof}
By construction we have a natural embedding $\bar\CL\subset\bar L$ such that $\bar\CL\cap\bar L_\alpha=\bar\CL'$. Moreover, by Lemma \ref{LocICRGoodRedLem2} the image $\rho\bar L$ is equal to $\rho\bar L_\alpha$ and generates~$\bar L_\alpha$. In particular we thus have $\rho\bar L\subset\bar L_\alpha$ and hence $\rho\bar\CL\subset \bar\CL\cap\bar L_\alpha=\bar\CL'$.
Since $[m_1],\ldots,[m_s]$ is an $R$-basis of~$\bar\CL'$, for every $1\le j\le s$ we can therefore write $[\rho m_j] = \sum_{i=1}^s \psi_{ij}[m_i]$ with elements $\psi_{ij}\in R$. The fact that $\rho\bar L_\alpha$ generates~$\bar L_\alpha$ over~$K$ then means that the matrix $(\psi_{ij})$ is invertible over~$K$. After inverting its determinant, we can thus assume that it is already invertible over~$R$. Then by construction $\bar\CL'$ is generated by~$\rho\bar\CL'$, and we are done.
\end{Proof}

\begin{Lem}\label{LocICRGoodRedLem7}
The  $\ICF$-sub-isocrystal $M_\alpha$ of $M$ comes from an $\ICF$-sub-isocrystal $\CM_\alpha$ of $\CM$ over $R$ which is pure of slope~$\alpha$, such that $\CM/\CM_\alpha$ is an $\ICF$-isocrystal.
\end{Lem}

\begin{Proof}
Combining Lemmas \ref{LocICRGoodRedLem4b} and \ref{LocICRGroucho} there exists a unique direct summand $\CL'\subset\CL$ such that $(\CL'+z\CL)/z\CL=\bar\CL'$ and $\rho\CL'\subset\CL'$. On the other hand, over $K$ we already know that $L_\alpha$ is a direct summand $K\bbu$-module of $L$ such that $(L_\alpha+zL)/zL=\bar L_\alpha'$ and $\rho L_\alpha\subset L_\alpha$. Since the tensor product $K\bbu \otimes_{R\bbu}\CL' \subset L$ satisfies the same properties, by the uniqueness part of Lemma \ref{LocICRGroucho} over $K$ we conclude that $K\bbu \otimes_{R\bbu}\CL' = L_\alpha$. 

In particular this implies that $\CL'=L_\alpha\cap\CM$. Since both $L_\alpha$ and $\CM$ are $A_\Fp\complot_{\BFq}\kern-1ptR$-submodules of~$M$, the same thus follows for~$\CL'$. Being also a free $R\bbu$-module of finite rank, by Lemma \ref{LatticeFreeLem} it is therefore finitely generated projective over $A_\Fp\complot_{\BFq}\kern-1ptR = \EAp{R}$. Moreover, by construction $\CM_\alpha := M_\alpha\cap\CM$ is the $R\ppu$-submodule generated by $L_\alpha\cap\CM$ and hence by~$\CL'$. Thus $\CM_\alpha$ is finitely generated projective over $\EFp{R}$. Since $\tau$ maps each of $\CM$ and $M_\alpha$ to itself, it follows that $\tau \CM_\alpha \subset \CM_\alpha$. 

On the other hand, the fact that $\rho\bar\CL'$ generates~$\bar\CL'$ implies that $\rho\CL'$ generates $\CL'$ over $R\bbu$. By the construction of $\rho$ this means that $\tau^{nd_\Fp} \CL_\alpha = \ICu^{n\alpha}\rho \CL_\alpha$ generates $\ICu^{n\alpha}\CL_\alpha$  over $R\bbu$ and hence also over $\CE_{A_\Fp,R}$. Therefore $\tau^{nd_\Fp} \CM_\alpha$ generates $\CM_\alpha$ over $\EFp{R}$, and so $\tau \CM_\alpha$ generates $\CM_\alpha$ over $\EFp{R}$. This means that $\tau$ induces an isomorphism $\sigma^*\CM_\alpha\isoto \CM_\alpha$; hence $\CM_\alpha$ is an $\ICF$-isocrystal over~$R$ with $i^*\CM_\alpha = M_\alpha$. Moreover, the lattice $\CL'$ shows that $\CM_\alpha$ is pure of slope~$\alpha$. 

Finally, since $\CM$ and $L$ are $\CE_{A_\Fp,R}$-submodules of~$M$, the same also holds for $\CL=L\cap\CM$. Thus the factor module $\CL/\CL'$ is an $\CE_{A_\Fp,R}$-module. But by construction that is a free module of finite rank over~$R\bbu$. By Lemma \ref{LatticeFreeLem} it is therefore a finitely generated projective module over $\CE_{A_\Fp,R}$. Inverting $\ICu$ it follows that $\CM/\CM_\alpha$ is a finitely generated projective module over $\EFp{R}$. It is thus a quotient $\ICF$-isocrystal of~$\CM$, and we are done.
\end{Proof}

\medskip 
To finish the proof of Theorem~\ref{LocICRGoodRedThm}, by induction on the rank of $\CM$ we can assume that the theorem holds for the $\ICF$-isocrystal $\CM/\CM_\alpha$. After inverting another non-zero element of $R$ we may thus assume that this possesses a slope flirtation over~$R$. 
Letting $\pi\colon \CM\onto \CM/\CM_\alpha$ denote the projection, one finds that 
$$\CM_{\le\beta}\ :=\ 
\scriptstyle\biggl\{\textstyle\begin{array}{cl}
0 & \hbox{if $\beta<\alpha$,}\\[3pt]
\pi^{-1}((\CM/\CM_\alpha)_{\le\beta}) & \hbox{if $\beta\ge\alpha$}\\[3pt]
\end{array}$$
defines a slope filtration of~$\CM$, and we are done.
\end{Proof}


\begin{Prop}\label{InvertRSlopeInt}
If $R$ is a discrete valuation ring, every $\ICF$-isocrystal of rank $1$ over~$R$ is pure of an integral slope.
\end{Prop}


\begin{Proof}
{}From Proposition \ref{InvertSlopeInt} we know that $i^*\CM$ is pure of an integral slope, and so $\CM$ is pure of the same slope by \cite[Prop.\,8.2.2]{MornevT}.
\end{Proof}

\medskip
Grothendieck's specialization theorem \cite[Thm.\,2.3.1]{KatzS} for Newton polygons of $p$-adic isocrystals now has the following analogue:

\begin{Thm}\label{LocICRGoodNewtonThm}
For any $\ICF$-isocrystal $\CM$ over~$R$ and any $x\in X$ the Newton polygon of $i_x^*\CM$ lies on or above the Newton polygon of $i^*\CM$ and has the same endpoint.
\end{Thm}

\begin{Proof}
For any $x\in X$ there exists a discrete valuation of $K$ centered at~$x$. As by Proposition \ref{LICBCFilt} the Newton polygon of $i_x^*\CM$ does not change under extension of the residue field~$k_x$, it therefore suffices to consider the case that $R$ is a discrete valuation ring.

Recall from Proposition \ref{LICRGoodRedRank} that $\CM$ is locally free of a unique finite rank~$r$. By Construction \ref{LICNewtonPolygon} both Newton polygons therefore extend precisely over the interval $[0,r]$. Take any integer $0\le s\le r$. Then for each of the $\ICF$-isocrystals $i_x^*\CM$ and $i^*\CM$, the ordinate of the point on its Newton polygon with the abscissa $s$ is the sum of its $s$ smallest slopes, counted with multiplicities. By Proposition \ref{ExtPowSlope} this sum is the smallest slope of $\bigwedge^si_x^*\CM$, respectively of $\bigwedge^si^*\CM$. Since the formation of exterior powers commutes with base change, these $\ICF$-isocrystals are $i_x^*\kern1pt\bigwedge^s\CM$ and $i^*\kern1pt\bigwedge^s\CM$, respectively. 

Let $\alpha$ denote the smallest slope of $i^*\kern1pt\bigwedge^s\CM$. Then by Proposition \ref{SlopesLowBound} there exist an integer $n>0$ with $n\alpha\in\BZ$ and a lattice ${L\subset i^*\kern1pt\bigwedge^s\CM}$ such that ${\tau^{n\ICd} L \subset \ICu^{n\alpha} L}$. Since $R$ is a discrete valuation ring, Lemma \ref{DVRLattice} implies that $\CL := L \cap \bigwedge^s\CM$ is a lattice in~$\bigwedge^s\CM$. By construction this lattice satisfies ${\tau^{n\ICd} \CL \subset \ICu^{n\alpha} \CL}$. This in turn yields a lattice $i_x^*\CL \subset i_x^*\kern1pt\bigwedge^s\CM$ such that ${\tau^{n\ICd} i_x^*\CL \subset \ICu^{n\alpha} i_x^*\CL}$. Applying Proposition \ref{SlopesLowBound} again over $k_x$ thus implies that the smallest slope of $i_x^*\kern1pt\bigwedge^s\CM$ is $\ge\alpha$. Thus the Newton polygon of $i_x^*\kern1pt\bigwedge^s\CM$ lies above the Newton polygon of $i^*\kern1pt\bigwedge^s\CM$ at the abscissa~$s$, as desired.

Finally consider the special case $s=r$. Then $\bigwedge^r\CM$ is an $\ICF$-isocrystal of rank $1$ and therefore pure by Proposition \ref{InvertRSlopeInt}. Thus $i_x^*\kern1pt\bigwedge^r\CM$ and $i^*\kern1pt\bigwedge^r\CM$ are pure of the same slope. This means that the Newton polygons of $i_x^*\CM$ and $i^*\CM$ have the same endpoint, and we are done.
\end{Proof}

\medskip
Using noetherian induction, Theorems \ref{LocICRGoodRedThm} and \ref{LocICRGoodNewtonThm} imply:

\begin{Thm}\label{LocICRNewtonSemicont}
For any piecewise linear function $\ell\colon [0,\rank(\CM)] \to \BR$, the set of points $x \in X$ such that the Newton polygon of $i_x^*\CM$ lies on or above the graph of $\ell$ is Zariski-closed in~$X$.
\end{Thm}


%
%

\section{Reduction isomorphisms for Weil representations}
\label{RedUM}

In this section we consider an $\ICF$-isocrystal $\CM$ with a slope filtration, defined over a noetherian integral domain~$R$, which we now assume to be normal. Let $i\colon R\into K$ denote the embedding of $R$ into its field of quotients. Set $X:=\Spec R$, and for each point $x\in X$ let $i_x\colon R\to k_x$ denote the natural homomorphism to the residue field at~$x$. The aim of this section is to compare the objects $T_\alpha(i_x^*\CM)$ and $U(i_x^*\CM)$ for all points $x\in X$. 

\medskip
For this fix a separable closure $K^\sep$ of $K$ and let $\WK\subset \GK \cong \Gal(K^\sep/K)$ be the associated Weil and Galois group, as in Section~\ref{WeilRep}. Let $\Rsep$ denote the integral closure of~$R$ in $K^\sep$. For any point $x\in X$ pick a point $\bar x \in \Spec \Rsep$ above~$x$. Let $D_x\subset\GK$ be the corresponding decomposition subgroup with the fixed field $K_x^\hens$, and let $I_x\subset D_x$ be the inertia subgroup with the fixed field $K_x^\shens$.

Since $R$ is normal, its localization $R_x$ at $x$ is again normal. Let $\FP$ be the prime ideal of $\Rsep$ corresponding to $\bar x$. Then by Raynaud \cite[Ch.\,X \S2 Thm.\,2]{RaynaudH}, the localization of $\Rsep\cap K_x^\hens$ at the prime ideal $\FP \cap K_x^\hens$ is a henselization $\Rxhens$ of $R_x$, and the localization of $\Rsep \cap K_x^\shens$ at the prime ideal $\FP \cap K_x^\shens$ is a strict henselization $\Rxshens$ of $R_x$.
On the other hand the local ring $\Rxhens$ is normal by construction, so by \cite[Tag 0BSM]{St} the natural homomorphism $D_x/I_x \to \Gkx$ is an isomorphism.
Together we thus obtain the following commutative diagram of rings and fields with their corresponding groups:
\UseTheoremCounterForNextEquation
\begin{equation}\label{RedUMBigDiag}
\vcenter{\xymatrix@R-10pt@C-10pt{
&& K^\sep \ar@{}[d]|{\textstyle\subsetup} 
&&& 1 \ar@{}[d]|{\textstyle\subsetdown}
&&& 1 \ar@{}[d]|{\textstyle\subsetdown} \\
\ \ \kxsep\ \ar@{<<-}[r] \ar@{}[d]|{\textstyle\subsetup} 
&\ \Rxshens\ \ar@{}[r]|{\textstyle\subset}  \ar@{}[d]|{\textstyle\subsetup} 
& K_x^\shens \ar@{}[d]|{\textstyle\subsetup} 
&&& I_x \ar@{}[d]|{\textstyle\subsetdown} 
&&& I_x \ar@{}[d]|{\textstyle\subsetdown} \\
\ k_x\ \ar@{<<-}[r] \ar@{}[d]|{\textstyle\equaldown} 
&\ \Rxhens\ \ar@{}[r]|{\textstyle\subset} \ar@{}[d]|{\textstyle\subsetup}  
& K_x^\hens \ar@{}[d]|{\textstyle\subsetup} 
&& \Gkx \ar@{<<-}[r] & D_x \ar@{}[d]|{\textstyle\subsetdown} 
&& W_{k_x} \ar@{<<-}[r] & \WK\cap D_x \ar@{}[d]|{\textstyle\subsetdown} \\
\ k_x\ \ar@{<-}[r] &
\ R\  \ar@{}[r]|{\textstyle\subset} 
& K
&&& \GK 
&&& \WK \\}}
\end{equation}
Note that the above choices are unique up to conjugation by~$\GK$. 

\begin{Thm}\label{LocICRGoodFiltThmT0}
For every $x\in X$ these choices induce a $D_x$-equivariant isomorphism 
$$T_0(i^*\CM)\ \cong\ T_0(i_x^*\CM).$$
In particular $I_x$ acts trivially on $T_0(i^*\CM)$, and the action of $D_x$ factors through~$\Gkx$.
\end{Thm}

\begin{Proof}
Since $\CM$ has a slope filtration, its subquotient $\CM_0=\CM_{\le0}/\CM_{<0}$ is an $\ICF$-isycrystal over~$R$, and the passage to this pure part of slope $0$ commutes with taking the pullback under $i\colon R\into K$ and under $i_x\colon R\to k_x$. Also, by \eqref{TWDef} we have natural isomorphisms $T_0(i^*\CM) \cong T_0(i^*\CM_0)$ and $T_0(i_x^*\CM) \cong T_0(i_x^*\CM_0)$. Thus it suffices to prove the theorem in the case that $\CM$ is pure of slope~$0$.

Choosing algebraic closures $\Kalg$ of $\Ksep$ and $\bar k_x$ of~$\kxsep$, from \eqref{WMSlope0} we then obtain natural isomorphisms 
$$\begin{array}{rcl}
T_0(i^*\CM)\ \cong\!\!& (i^*\CM\complot_K \Kalg)^\tau &\!\!\cong\ 
(\CM\complot_R \Kalg)^\tau , \\[3pt]
T_0(i_x^*\CM)\ \cong\!\!& (i_x^*\CM\complot_{k_x}\kern-1pt\bar k_x)^\tau &\!\!\cong\ 
(\CM\complot_R \bar k_x)^\tau .
\end{array}$$
To prove the theorem it thus suffices to show that the natural maps
\UseTheoremCounterForNextEquation
\begin{equation}\label{LocICRGoodFiltThmT0MTauMaps}
\xymatrix@C-10pt{
(\CM\complot_R\kern-1pt \Kalg)^\tau 
& \ar@{_{ (}->}[l] (\CM\complot_R\kern-1pt \Ksep)^\tau
& \ar@{_{ (}->}[l] (\CM\complot_R\kern-1pt \Rxshens)^\tau \ar[r]
& (\CM\complot_R \kxsep)^\tau \ar@{^{ (}->}[r] 
& (\CM\complot_R \bar k_x)^\tau}
\kern-5pt
\end{equation}
are isomorphisms.

For this observe that, since $\CM$ is pure of slope~$0$, by Definition \ref{LocICRSlopeDef} there exist an integer $n>0$ and an $\ICO\complot_{\BFq}R$-lattice $\CL\subset \CM$ such that $\tau^n\CL$ generates~$\CL$. We would like to express the modules in \eqref{LocICRGoodFiltThmT0MTauMaps} in terms of $\tau$-invariants relative to~$\CL$, but cannot do that directly, because we do not know whether $\tau \CL$ is contained in~$\CL$. But we can first look at the natural maps of $\tau^n$-invariants
\UseTheoremCounterForNextEquation
\begin{equation}\label{LocICRGoodFiltThmT0LTauNMaps}
\xymatrix@C-10pt{
(\CL\complot_R\kern-1pt \Kalg)^{\tau^n} 
& \ar@{_{ (}->}[l] (\CL\complot_R\kern-1pt \Ksep)^{\tau^n}
& \ar@{_{ (}->}[l] (\CL\complot_R\kern-1pt \Rxshens)^{\tau^n} \ar[r]
& (\CL\complot_R \kxsep)^{\tau^n} \ar@{^{ (}->}[r] 
& (\CL\complot_R \bar k_x)^{\tau^n} \rlap{.}}
\end{equation}

\begin{Lem}\label{LocICRGoodFiltThmT0Lem}
These maps are isomorphisms.
\end{Lem}

\begin{Proof}
Recall that $\ICO=\ICk\bbu$ is a finitely generated free module over its subring $\BFq\bbu$. Therefore $\CL_x^\shens := \CL\complot_R\kern-1pt \Rxshens$ is a finitely generated projective module over the subalgebra ${\BFq\bbu\complot_{\BFq}\Rxshens} \cong \Rxshens\bbu$. The action of $\tau^n$ thus provides $\CL_x^\shens$ with the structure of a dualizable left $\Rxshens\bbu[\tau^n]$-module.

Next observe that $\Rxshens$, being strictly henselian, contains a subfield $\BFqn$ of order~$q^n$. Thus $\CL_x^\shens$ is a \emph{$\sigma$-module over the difference ring $(\Rxshens\bbu,\sigma^n)$} in the sense of \cite[Def.\,2.3.1]{MornevT}. By \cite[Thm. 3.3.2]{MornevT} its $\tau^n$-invariants thus form a lisse sheaf of finitely generated free $\BFqn\bbu$-modules on the \'etale site of $\Spec \Rxshens$. But as $\Rxshens$ is a strictly henselian local ring, this sheaf is constant, and so all the maps in question are isomorphisms.
\end{Proof}

\medskip
Inverting $\ICu$ in \eqref{LocICRGoodFiltThmT0LTauNMaps} we find that the natural maps
$$\xymatrix@C-10pt{
(\CM\complot_R\kern-1pt \Kalg)^{\tau^n} 
& \ar@{_{ (}->}[l] (\CM\complot_R\kern-1pt \Ksep)^{\tau^n}
& \ar@{_{ (}->}[l] (\CM\complot_R\kern-1pt \Rxshens)^{\tau^n} \ar[r]
& (\CM\complot_R \kxsep)^{\tau^n} \ar@{^{ (}->}[r] 
& (\CM\complot_R \bar k_x)^{\tau^n}}$$
are isomorphisms. By definition they are also equivariant under the induced action of~$\tau$. They therefore induce isomorphisms on the submodules of $\tau$-invariants in \eqref{LocICRGoodFiltThmT0MTauMaps}, as desired.
\end{Proof}

\medskip
Next let $\bar k$ denote the algebraic closure of $\BFq$ in~$K^\sep$. For any point $x\in X$ this is contained in $\Rxshens$, and under the reduction homomorphism it maps isomorphically to the algebraic closure of $\BFq$ in~$\kxsep$. We can therefore use the same $\bar k$ in the construction of both $U(i^*\CM)$ and $U(i_x^*\CM)$ from Section~\ref{WeilRep}.

\begin{Thm}\label{LocICRGoodFiltThmU}
For any $x\in X$ and $\alpha\in\BQ$ the above choices induce $\WK\cap D_x$-equivariant isomorphisms
$$\begin{array}{rl}
T_\alpha(i^*\CM) &\!\cong\ T_\alpha(i_x^*\CM), \\[3pt]
U(i^*\CM) &\!\cong\ U(i_x^*\CM).
\end{array}$$
In particular $I_x$ acts trivially on the left hand side, and the action of $\WK\cap D_x$ factors through~$W_{k_x}$. 
\end{Thm}

\begin{Proof}
By the slope decomposition from Proposition \ref{WMBigWeilRepGrad} and the isomorphism \eqref{WMSmallWeilRepIsom}, the statement for $U(i^*\CM)$ follows from that for $T_\alpha(i^*\CM)$ for all~$\alpha$. Moreover, by the same arguments as in the beginning of the proof of Theorem \ref{LocICRGoodFiltThmT0} with $\CM_\alpha$ in place of~$\CM_0$, we can reduce ourselves to the case that $\CM$ is pure of slope~$\alpha$.

Consider now the following inner hom $\ICF$-isocrystal over $\Rxshens$:
$$\CN\ :=\ \hom\bigl(\SA{\kalg}\complot_{\kalg} \Rxshens, \CM\complot_R\kern-1pt \Rxshens\bigr)$$
Since purity is invariant under base change, the same argument as in the proof of Proposition \ref{TensorPure} shows that $\CN$ is pure of slope~$0$. Basic properties of the Tate module from Section~\ref{WeilRep} thus yield natural isomorphisms
$$\xymatrix@R-10pt{
T_\alpha(i^*\CM)
\ar@{}[r]|-{\textstyle\cong}^{\raisebox{6pt}{$\scriptstyle\eqref{TWDef}$\kern5pt}}
& T_{\SA{\kalg}}(i^*\CM)
\ar@{}[r]|-{\textstyle\cong}^{\raisebox{6pt}{$\scriptstyle\eqref{HomPMW}$\kern41pt}}
& \kern-2pt\Hom\bigl(\SA{\kalg}\complot_{\kalg} \Kalg,\CM\complot_R\kern-1pt \Kalg\bigr)
\ar@{}[d]|{\kern15pt\textstyle\Vert\kern-3pt\wr\;\scriptstyle\eqref{Homhomtau}} \\
T_0(i^*\CN)
\ar@{}[r]|-{\textstyle\cong}^{\raisebox{6pt}{$\scriptstyle\eqref{WMSlope0}$\kern13pt}}
& (\CN\complot_{R_x^{\kern+1pt\shens}}\kern-2pt \Kalg)^\tau
\ar@{}[r]|-{\textstyle{\cong}}
& \hom\bigl(\SA{\kalg}\complot_{\kalg} \Kalg,\CM\complot_R\kern-1pt \Kalg\bigr)^\tau \rlap{.}\\}$$
Moreover, as the decomposition subgroup $D_x$ preserves the subring $\Rxshens \subset \Kalg$, the formula \eqref{HomPMWgamma} defines an $\EFp{\Rxhens}[\tau]$-linear action of $\WK\cap D_x$ on the module
$$\Hom_{\EFp{\kalg}}(\SA{\kalg},\CM\complot_R \Rxshens)\ \cong\ \CN.$$
This induces an action on $T_0(i^*\CN)$ for which the isomorphism $T_\alpha(i^*\CM) \cong T_0(i^*\CN)$ is equivariant.
On the other hand, since $\CN$ is pure of slope~$0$, as in \eqref{LocICRGoodFiltThmT0MTauMaps} we have 
natural isomorphisms
\UseTheoremCounterForNextEquation
\begin{equation}
\vcenter{\vbox{\xymatrix@C-10pt@R-15pt{
(\CN\complot_{\Rxshens}\kern-1pt \Kalg)^\tau 
& \ar@{_{ (}->}[l] (\CN\complot_{\Rxshens}\kern-1pt \Ksep)^\tau
& \ar@{_{ (}->}[l] \CN^\tau \ar[r]
& (\CN\complot_{\Rxshens} \kxsep)^\tau \ar@{^{ (}->}[r] 
& (\CN\complot_{\Rxshens} \bar k_x)^\tau\\
T_0(i^*\CN) \ar@{=}[u]&&&& T_0(i_x^*\CN).\ar@{=}[u]}}}
\end{equation}
These isomorphisms are $\WK\cap D_x$-equivariant by the definition of the action on $\CN$.
Together this therefore yields a $\WK\cap D_x$-equivariant isomorphism $T_\alpha(i^*\CM) \cong T_\alpha(i^*_x\CM)$, finishing the proof.
\end{Proof}

\newpage
\bigskip
\noindent{\huge\bf Part II: Global isocrystals}
\addtocontents{toc}{\medskip}
\addcontentsline{toc}{section}{\large Part II: Global isocrystals}

\section{Global isocrystals over a field}
\label{GlobCrys}
Now let $F$ be the function field of a smooth connected projective algebraic curve $C$ over~$\BFq$. In this section $K$ denotes an arbitrary field over~$\BFq$.

\begin{Lem}\label{GLICLem1}
The ring $F \otimes_{\BFq}\kern-1pt K$ is a non-empty finite product of integral domains that are transitively permuted by the partial Frobenius endomorphism $\sigma\colon a\otimes\xi \mapsto a\otimes\xi^q$.
\end{Lem}

\begin{Proof}
Let $\kappa$ denote the algebraic closure of $\BFq$ in~$F$. Then $\kappa$ is a finite separable extension of~$\BFq$; hence $\kappa\otimes_{\BFq}\kern-1ptK$ is a non-empty finite product of fields 
that are transitively permuted by~$\sigma$. Also $C$ is a geometrically irreducible smooth curve over~$\kappa$. Tensoring $F$ over $\kappa$ with each factor of $\kappa\otimes_{\BFq}\kern-1ptK$ thus yields an integral domain. The lemma therefore follows from the natural isomorphism $F\otimes_{\BFq}\kern-1pt K \cong F \otimes_{\kappa}\kern-1pt (\kappa\otimes_{\BFq}\kern-1pt K)$.
\end{Proof}

\medskip
We are interested in the total ring of quotients
\UseTheoremCounterForNextEquation
\begin{equation}\label{GCRingDef}
\EFK := \Quot(F \otimes_{\BFq}\kern-1pt K).
\end{equation}
Lemma \ref{GLICLem1} implies that this is a non-empty finite product of fields, and that $\sigma$ extends uniquely to an endomorphism of $\EFK$ that permutes the factors transitively.

%

\begin{Def}\label{GlobICDef}
A dualizable left $\EFK[\tau]$-module is called an \emph{$F$-isocrystal over $K$ (with the ground field~$\BFq$).}
\end{Def}

The category of $F$-isocrystals over $K$ is denoted by $\GLIC{K}$. The tensor product from \S\ref{DiffRingIC} provides this category with a rigid symmetric monoidal structure.

\medskip
Any field embedding $i\colon K \into K'$ induces a morphism of difference rings ${\EFK\into \EF{K'}}$ and therefore gives rise to a monoidal base change functor $i^*\colon \GLIC{K} \to \GLIC{K'}$, as in Construction \ref{LICBC}.

\medskip
Similarly, for any finite extension $F/F'$, giving an $F$-isocrystal $M$ over $K$ is equivalent to giving an $F'$-isocrystal $M'$ over $K$ together with an embedding $F\into\End(M')$ over~$F'$, as in Construction \ref{LICCoeffRes}.

\begin{Prop}\label{GLICTannakian}
\begin{enumerate}\StatementLabels
\item\label{GLICTannakianTannakian}
The category $\GLIC{K}$ is $F$-linear tannakian.
\item\label{GLICTannakianFiber}
The forgetful functor to $\EFK$-modules is a fiber functor.
\item\label{GLICTannakianConstRank}
Every $F$-isocrystal over $K$ is free of a unique finite rank as an $\EFK$-module.
\end{enumerate}
\end{Prop}

\begin{Proof}
This follows at once from Proposition \ref{Tannaka} as soon as we show that
the ring of invariants of $\sigma$ on $\EFK$ is the image of ${F\into\EFK}$, ${a\mapsto a\otimes1}$. 

To prove that, choose a subring $A\subset F$ such that $U:=\Spec A$ is an open chart in~$C$. Then $\EFK$ is the ring of rational functions on the curve $U_K := \Spec(A\otimes_{\BFq}\kern-2pt K)$ over~$K$. Since $\sigma$ transitively permutes the irreducible components of~$U_K$, any non-zero $\sigma$-invariant element $e\in\EFK$ does not vanish at the generic points of~$U_K$ and so has a well-defined divisor. Being $\sigma$-invariant, its pole divisor $D$ is then mapped to itself under pullback by the morphism $U_K\to U_K$ induced by~$\sigma$. By \cite[Lemma~8.9]{Pink2012} it follows that $D$  is the pullback of an effective divisor $D_0$ on~$U$. Choose a non-zero $a\in A$ whose zero divisor contains~$D_0$. Then $(a\otimes1)e$ is a regular function on $U_K$ and thus an element of $A\otimes_{\BFq}\kern-1ptK$. Being again $\sigma$-invariant, it therefore lies in $A\otimes_{\BFq}\kern-1ptK^\sigma$ with $K^\sigma=\BFq$. It is thus equal to $b\otimes1$ for some $b\in A$, and so $e=b/a\otimes1$ with $b/a\in \Quot(A)=F$, as desired.
\end{Proof}

\begin{Cor}\label{GLICTannCor}
\begin{enumerate}\StatementLabels%
\item\label{GLICTannCor1}
Every $F$-isocrystal is a finite iterated extension of simple objects.
\item\label{GLICTannCor2}
For all $F$-isocrystals $M$, $N$ the $F$-vector space
$\Hom(M,N)$ is finite-dimensional.
\end{enumerate}\end{Cor}

\begin{Proof}
This follows from Proposition \ref{GLICTannakian} \ref{GLICTannakianTannakian} by \cite[Prop.\,2.13]{DeligneT}.
\end{Proof}

\medskip
Next, let $\Fp$ be an arbitrary place of $F$ and let $\ICF$ be the $\Fp$-adic completion of~$F$.

\begin{Lem}\label{GLICBCLem}
The canonical homomorphism $F \otimes_{\BFq}\kern-1pt K \into \ICF\complot_{\BFq}\kern-1.5pt K$ 
extends uniquely to a homomorphism of difference rings
$$\ell_\Fp\colon\ \EFK\ \longto\ \EFpK = \ICF\complot_{\BFq}\kern-1.5pt K.$$
\end{Lem}

\begin{Proof}
Let $A_{(\Fp)}\subset F$ denote the valuation ring associated to~$\Fp$.
This is a localization of an $\BFq$-algebra of finite type; hence the ring $A_{(\Fp)}\otimes_{\BFq}\kern-1pt K$ is noetherian. Also $A_\Fp$ is the completion of $A_{(\Fp)}$, so the ring $A_\Fp\complot_{\BFq}\kern-1pt K$ is flat over $A_{(\Fp)}\otimes_{\BFq}\kern-1pt K$. Inverting a uniformizer of $A_{(\Fp)}$ it follows that $\EFpK$ is flat over $F\otimes_{\BFq}\kern-1pt K$. Since any flat ring homomorphism maps non-zerodivisors to non-zerodivisors, it extends uniquely to a homomorphism of the total rings of fractions. This yields the desired homomorphism $\ell_\Fp$.
This homomorphism commutes with $\sigma$ because so does the composite homomorphism
$F\otimes_{\BFq}\kern-1pt K \to \EFpK$. 
\end{Proof}

\medskip
By Construction \ref{BC1} the homomorphism $\ell_\Fp$ gives rise to a monoidal base change functor 
\UseTheoremCounterForNextEquation
\begin{equation}\label{GlobLocICCons}
\ell_\Fp\kern-4pt{}^*\colon \GLIC{K} \to \LIC{K},\ 
M\mapsto M_\Fp := \ell_\Fp\kern-4pt{}^*M.
\end{equation}
Using this and the functor $U$ from \eqref{UMW}, an $F$-isocrystal over $K$ gives rise to a whole system of Weil representations associated to all places of~$F$. 
To analyze this system of representations we first look at the special case where $K$ is finite.

\section{Global isocrystals over a finite field}
\label{FinField}

In this section we study $F$-isocrystals over a finite field $k$. In this case the ring $F\otimes_{\BFq} k$ is a finite product of fields, so from \eqref{GCRingDef} we obtain
\UseTheoremCounterForNextEquation
\begin{equation}\label{GlobFinFieldRing}
\EFk := \Quot(F \otimes_{\BFq}\kern-1pt k) = F \otimes_{\BFq}\kern-1pt k.
\end{equation}

\begin{Cons}\label{GlobCharPol}\rm  
For any $F$-isocrystal $M$ over $k$ it follows as in Construction \ref{LocCharPol} that $\tau_M\colon {M \to M}$ is bijective and that $\tau_M^{[k/\BFq]}\colon M \isoto M$ is $\EFk$-linear. Since $M$ is a free $\EFk$-module,
\UseTheoremCounterForNextEquation
\begin{equation}\label{CharPolFormulaGlobal}
\charact_{M}(X) \ :=\ \det\nolimits_{\EFk}\bigl(X\cdot\id_M - \tau_M^{[k/\BFq]}\bigr).
\end{equation}
is therefore a monic polynomial of degree $\rank(M)$ with coefficients in~$\EFk$. Also, from \cite[Lemma~8.1.4]{BoecklePink} we deduce that
\UseTheoremCounterForNextEquation
\begin{equation}\label{BPLemmaGlobal}
\charact_{M}(X^{[k/\BFq]})\ =\ \det\nolimits_{F}\bigl(X\cdot\id_M - \tau_M\bigr).
\end{equation}
Thus $\charact_{M}$ actually has coefficients in~$F$. We call it the \emph{characteristic polynomial of~$M$.}
\end{Cons}

As in Section~\ref{LocFinField}, for each place $\Fp$ of $F$ we pick an algebraic closure $\ICFalg$ of~$\ICF$ and denote by $\ord_\Fp$ the unique valuation on $\ICFalg$ that extends the normalized valuation on $\ICF$. The following theorem generalizes \cite[Prop.\,5.6.9]{GossBook} to $F$-isocrystals:

\begin{Thm}\label{CharPolThm}
For any $F$-isocrystal $M$ over $k$ and any place $\Fp$ of $F$ the associated local isocrystal $M_\Fp$ over~$k$ has the following properties:
\begin{enumerate}\StatementLabels
\item\label{CharPolThmGalois}
We have an equality of polynomials in $\EFp{\kalg}[X]${\rm:}
$$\charact_{M}(X) \ =\ 
\det\nolimits_{\EFp{\kalg}}\bigl(X\cdot\id - \Frob^{-1}_k \bigm| U(M_\Fp) \bigr).$$
In particular the polynomial on the right has coefficients in~$F$ and is independent~of~$\Fp$.
\item\label{CharPolThmNewton}
For any $\alpha\in\BQ$, the number of roots $\lambda\in\ICFalg$ of $\charact_{M}$ that satisfy 
$$\ord_\Fp(\lambda)=\alpha\cdot\frac{[k/\BFq]}{\ICd},$$
counted with multiplicities, is equal to $\rank(M_{\Fp,\alpha})$.
\item\label{CharPolThmSlopesUnitRoots}
In particular $M_\Fp$ is pure of slope~$0$ if and only if all roots of $\charact_{M}$ are units above~$\Fp$.
\end{enumerate}
\end{Thm}

\begin{Proof}
The defining formulas \eqref{CharPolFormulaLocal} and \eqref{CharPolFormulaGlobal} immediately imply that $\charact_{M} = \charact_{M_\Fp}$. Thus \ref{CharPolThmGalois} follows from Corollary \ref{CharPolGalois1}, and the rest follows from Proposition \ref{SlopesRoots}.
\end{Proof}

\section{Models of global isocrystals over a ring}
\label{GoodRedGlob}

Now let $M$ be an $F$-isocrystal over an arbitrary field~$K$. 

\medskip
To formulate our first result consider $\BFq$-subalgebras $A\subset F$ and $R\subset K$ and a non-zerodivisor $e\in A\otimes_{\BFq}\kern-1ptR$. Let $S_e$ be the multiplicative subset of $A\otimes_{\BFq}\kern-1ptR$ that is generated by the elements $\sigma^\nu(e)$ for all $\nu\ge0$. Then the localization 
\UseTheoremCounterForNextEquation
\begin{equation}\label{EAReDef}
\CE_{A,R,e}\ :=\ S_e^{-1}(A\otimes_{\BFq}\kern-1ptR)
\end{equation}
is a subring of $\EFK = \Quot(F \otimes_{\BFq}\kern-1pt K)$ that is closed under the partial Frobenius endomorphism~$\sigma$. We can therefore speak of dualizable left $\CE_{A,R,e}[\tau]$-modules as in Definition \ref{DiffRingICDef}.
Moreover, any such module yields an $F$-isocrystal over $K$ by base change via the inclusion of difference rings $\CE_{A,R,e} \into \EFK$.

\begin{Prop}\label{GlobICFpGoodRedAlg}
There exist finitely generated $\BFq$-subalgebras $A\subset F$ and $R\subset K$ and a non-zerodivisor $e\in A\otimes_{\BFq}\kern-1ptR$ such that $M$ arises by base change from some dualizable left $\CE_{A,R,e}[\tau]$-module~$\CMod$.
\end{Prop}

\begin{Proof}
From Proposition \ref{GLICTannakian} \ref{GLICTannakianConstRank} we know that $M$ is a free $\EFK$-module of some finite rank~$r$. Choosing a basis, we can therefore represent $\tau$ by an $r\times r$-matrix $\Phi$ with coefficients in~$\EFK$. By the definition of $\EFK$ there then exists a non-zerodivisor $\theta\in F\otimes_{\BFq}\kern-1ptK$ such that $\theta\Phi$ has coefficients in $F\otimes_{\BFq}\kern-1ptK$. Thus there exist finitely generated $\BFq$-subalgebras $A\subset F$ and $R\subset K$ such that $\theta$ and all coefficients of $\theta\Phi$ lie in the subring $A \otimes_{\BFq}\kern-1pt R$. 

That $\taulin_M$ is an isomorphism now implies that $\det(\Phi)$ is a unit in~$\EFK$. Thus $\theta$ and $\theta^r\kern-2pt\det(\Phi) = \det(\theta\Phi)$ are non-zerodivisors in $A\otimes_{\BFq}\kern-.5ptR$; hence so is the element $e := \theta^{r+1}\kern-2pt\det(\Phi)$. Therefore $\theta$ and $\det(\Phi)$ become units in the ring $\CE_{A,R,e}$, and so $\Phi = \theta^{-1}(\theta\Phi)$ has coefficients in $\CE_{A,R,e}$. But this means that the $\CE_{A,R,e}$-submodule $\CMod\subset M$ that is generated by the chosen basis is a dualizable left $\CE_{A,R,e}[\tau]$-module whose base change to $\EFK$ is~$M$, as desired.
\end{Proof}

\medskip
One might consider the $\CMod$ in Proposition \ref{GlobICFpGoodRedAlg} as a kind of model  of $M$ over the ring~$R$. But we avoid building a theory of such models, because the concept depends too much on the choice of~$e$.

\medskip
Proposition \ref{GlobICFpGoodRedAlg} implies that $M$ arises by base change from an $F$-isocrystal over the subfield $\Quot(R)$. For simplicity we therefore assume that $\Quot(R)=K$. For the following we abbreviate $X := \Spec R$ and let $D\subset(\Spec A)\times_{\BFq}\kern-1ptX$ denote the zero divisor of~$e$. 

\begin{Prop}\label{AllAssOnAReD}
In Proposition \ref{GlobICFpGoodRedAlg} we can achieve in addition 
\begin{enumerate}\StatementLabels
\item\label{AllAssOnAReDA}
that $\Spec A$ is a dense open chart in~$C$, and
\item\label{AllAssOnAReDeD}
that $D$ contains no irreducible component of $\Fp\times_{\BFq}\kern-1pt X$ for any point $\Fp\in\Spec A$ or of ${(\Spec A)\times_{\BFq}\kern-1pt x}$ for any point $x\in X$.
\end{enumerate}
\end{Prop}

\begin{Proof}
The finitely many generators of $A$ are regular outside some non-empty finite closed subset $S\subset C$; hence $A$ is contained in the subring $\Gamma(C\setminus S,\CO_C)$. Since $C\smallsetminus S$ is affine of finite type over~$\BFq$, this subring is finitely generated over $\BFq$ and represents the chart $C\smallsetminus S$ of~$C$. Thus we can achieve~\ref{AllAssOnAReDA}. 



Now let $\eta$ denote the generic point of~$X$. Since $\Spec A$ and $X$ are noetherian, the locus of points $x\in X$ such that $D$ contains no irreducible component of $(\Spec A)\times_{\BFq}\kern-1ptx$ is a constructible subset of~$X$. As the divisor $D$ contains no irreducible component of $(\Spec A)\times_{\BFq}\kern-1pt\eta$, 
that constructible subset is a neighborhood of~$\eta$. After inverting a suitable element of $R$ we can therefore achieve that $D$ contains no irreducible component of $(\Spec A)\times_{\BFq}\kern-1ptx$ for any $x\in X$. The same argument with $A$ in place of $R$ finishes the proof of~\ref{AllAssOnAReDeD}.
\end{Proof}

\medskip
For the following we choose $A,R,e,D$ as in Propositions \ref{GlobICFpGoodRedAlg} and~\ref{AllAssOnAReD}. As before, for each point $x\in X$ we write $i_x\colon R\to k_x$ for the natural homomorphism to the residue field. For any $f\in R$ we abbreviate $f(x) := i_x(f)$.

\begin{PropCons}\label{GlobRedCons}
For every point $x\in X$ the homomorphism of difference rings $\id\otimes i_x\colon \allowbreak {A\otimes_{\BFq}R} \to {F\otimes_{\BFq}k_x}$ extends uniquely to a homomorphism of difference rings $\CE_{A,R,e} \to \EF{k_x}$. By base change $\CMod$ thus yields an $F$-isocrystal $\CMod_x$ over $k_x$ that we call the \emph{specialization of $\CMod$ at~$x$}.
\end{PropCons}

\begin{Proof}
For every $\nu\ge0$, the zero locus of $\sigma^\nu(e)$ is the pullback $(\sigma^\nu)^*(D)\subset(\Spec A)\times_{\BFq}\kern-1ptX$. The zero locus of $(\id\otimes i_x)(\sigma^\nu(e)) \in {A\otimes_{\BFq}k_x}$ is therefore the intersection of $(\sigma^\nu)^*(D)$ with $(\Spec A)\times_{\BFq}\kern-1ptx$. Since $\sigma$ maps the subscheme $(\Spec A)\times_{\BFq}\kern-1ptx$ onto itself, it follows from \ref{AllAssOnAReD}~\ref{AllAssOnAReDeD} that this intersection contains no irreducible component of $(\Spec A)\times_{\BFq}\kern-1pt x$. Thus $(\id\otimes i_x)(\sigma^\nu(e))$ is a non-zerodivisor, so the homomorphism $\id\otimes i_x$ extends as stated.
\end{Proof}

\medskip
For every place $\Fp$ of $F$ we next want to construct an $F_\Fp$-isocrystal associated to~$\CMod$. For this we need to invert a suitable non-zero element of~$R$, as follows:

\begin{PropCons}\label{GlobICFpGoodRed}
For every place $\Fp$ of~$F$ there exists a non-zero element $\epsilon_\Fp\in R$ such that the natural homomorphism $A\otimes_{\BFq}\kern-1pt R \longto F_\Fp\complot_{\BFq}\kern-1ptR$ extends uniquely to a homomorphism of difference rings
$$\ell_\Fp\colon \CE_{A,R,e}\ \longto\ F_\Fp\complot_{\BFq}\kern-1.5ptR[\epsilon_\Fp^{-1}].$$
By base change this induces an $F_\Fp$-isocrystal $\CMod_\Fp := \ell_\Fp^{\kern+1pt*}\CMod$ over~$R[\epsilon_\Fp^{-1}]$. Moreover, for a suitable choice of $\epsilon_\Fp$ the isocrystal $\CMod_\Fp$ possesses a slope filtration. Finally, for every $x\in X$ with $\epsilon_\Fp(x)\not=0$, there is a natural isomorphism $i_x^*\CMod_\Fp \cong \CMod_{x,\Fp}$.
 \end{PropCons}

\begin{Proof}
First we claim that the inclusion $A\otimes_{\BFq}\kern-1ptR \into \ICF\complot_{\BFq}\kern-1.5pt R$ extends to a homomorphism 
\UseTheoremCounterForNextEquation
\begin{equation}\label{GlobICFpGoodRedMap}
\CE_{A,R,e}\ \into\ \ICF\complot_{\BFq}\kern-1.5pt R[\xi^{-1}]\ =\ \EFp{R[\xi^{-1}]}
\end{equation}
for some non-zero element $\xi\in R$. To show this we first note that by generic regularity \cite[Cor.\,6.12.5]{EGA4}
there exists a non-zero element $\xi\in R$ such that $R[\xi^{-1}]$ is regular. Next observe
that, since $e$ is a non-zerodivisor in $A \otimes_{\BFq}\kern-1pt R$, it is also a non-zerodivisor in $F \otimes_{\BFq}\kern-1pt K$. Its image in $\ICF\complot_{\BFq}\kern-1pt K$ is therefore a unit by Lemma \ref{GLICBCLem}. Thus $e$ is also a non-zerodivisor in the subring $\ICF\complot_{\BFq}\kern-1pt R[\xi^{-1}]$. Next, the isomorphism $\ICF\cong\ICk\ppu$ yields an isomorphism
$$\ICF\complot_{\BFq}\kern-1.5pt R[\xi^{-1}]\ \cong\ \bigl(\ICk\otimes_{\BFq}\kern-1.5pt R[\xi^{-1}]\bigr)\ppu.$$
Here the fact that $\ICk$ is a finite separable extension of $\BFq$ implies that $\ICk\otimes_{\BFq}\kern-1pt K$ is a non-empty finite product of fields~$K_j$.

As $R[\xi^{-1}]$ is regular, it follows (for instance by \cite[Exp.\,I, Th.\,9.5\,(i)]{SGA1})
that ${\ICk\otimes_{\BFq}\kern-1pt R[\xi^{-1}]}$ is a finite product of integral domains $R_j\subset K_j$. 
Thus $\ICF\complot_{\BFq}\kern-1pt R[\xi^{-1}]$ is isomorphic to a finite product of the rings $R_j\ppu$. The non-zerodivisor $e$ therefore corresponds to a tuple of non-zero elements $e_j\in R_j\ppu$. Let $\xi_j\in R_j$ denote the leading term of~$e_j$. Then $\prod_j(R_j/\xi_jR_j)$ is a quotient of $\prod_jR_j \cong \ICk\otimes_{\BFq}\kern-1pt R[\xi^{-1}]$ and hence a finitely generated $R[\xi^{-1}]$-module. Moreover, since all $\xi_j$ are non-zero, this module becomes zero on tensoring with $K$ over~$R[\xi^{-1}]$. Thus the module is annihilated by some non-zero element of $R[\xi^{-1}]$. After replacing $\xi$ by a multiple this module therefore becomes zero, which means that $\xi_j$ becomes a unit in~$R_j$. Thus each $e_j$ becomes a unit in $R_j\ppu$, and so $e$ becomes a unit in
$$\ICF\complot_{\BFq}\kern-1.5pt R[\xi^{-1}]\ \cong\ \bigl(\ICk\otimes_{\BFq}\kern-1.5pt R[\xi^{-1}]\bigr)\ppu.$$
For every $\nu\ge0$ the same therefore follows for $\sigma^\nu(e)$. The claim thus results from the construction~\eqref{EAReDef}.

By base change via \eqref{GlobICFpGoodRedMap}, the dualizable left $\CE_{A,R,e}[\tau]$-module~$\CMod$ from Proposition \ref{GlobICFpGoodRedAlg} now yields an $\ICF$-isocrystal over $R[\xi^{-1}]$. After inverting another non-zero element that acquires a slope filtration by Theorem \ref{LocICRGoodRedThm}. For a suitable non-zero multiple $\epsilon_\Fp$ of~$\xi$ the construction thus yields an $\ICF$-isocrystal $\CMod_\Fp$ over $R[\epsilon_\Fp^{-1}]$ which possesses a slope filtration.

Finally, for every $x\in X$ with $\epsilon_\Fp(x)\not=0$, we have a commutative diagram of difference rings
$$\xymatrix{\CE_{A,R,e} \ar@{^{ (}->}[r]^{\ell_\Fp} \ar[d] & 
\EFp{R[\epsilon_\Fp^{-1}]} \kern-10pt \ar[d] \\
\EF{k_x} \ar@{^{ (}->}[r]^{\ell_\Fp} & \EFp{k_x}\rlap{.}}$$
Thus the base change $i_x^*\CMod_\Fp$ is naturally isomorphic to the $\ICF$-isocrystal $\CMod_{x,\Fp}$ associated to~$\CMod_x$, and we are done.
\end{Proof}

\medskip
Next, by \emph{primes of $A$} we will mean any maximal ideals of~$A$. By the condition~\ref{AllAssOnAReD}~\ref{AllAssOnAReDA} these are in bijection with all but finitely many places of~$F$, and for any such $\Fp$ the associated $\Fp$-adic completion of $A$ is the ring of integers $\ICO$ of~$\ICF$. 

For these $\Fp$ we can be more precise about $\epsilon_\Fp$ and $\CMod_\Fp$. Let $\proj_\Fp$ denote the natural projection $A\onto k_\Fp$, and let $\Norm_\Fp$ the relative norm map $\ICk\otimes_{\BFq}\kern-1ptR \to \BFq\otimes_{\BFq}\kern-1ptR$.

\begin{Prop}\label{GlobLocICSlope0Extends}
For every prime $\Fp$ of $A$, we can choose $\epsilon_\Fp\in R\smallsetminus\{0\}$ such that
\UseTheoremCounterForNextEquation
\begin{equation}\label{epDef}
\Norm_\Fp\bigl((\proj_\Fp\otimes\id)(e)\bigr) = 1\otimes \epsilon_\Fp.
\end{equation}
Moreover, the $\ICF$-isocrystal $\CMod_\Fp$ over $R[\epsilon_\Fp^{-1}]$ is then pure of slope~$0$.
\end{Prop}

\begin{Proof}
Let $\epsilon_\Fp$ be defined by the formula \eqref{epDef}. Then the first condition in \ref{AllAssOnAReD}~\ref{AllAssOnAReDeD} implies that $(\proj_\Fp\otimes\id)(e)$ is a non-zerodivisor in $\ICk\otimes_{\BFq}\kern-1ptR$. Therefore $\epsilon_\Fp$ is non-zero.

Next recall that the residue field $\ICk$ is a finite extension of degree $\ICd$ of~$\BFq$. The relative norm map $\Norm_\Fp$ is therefore the product of the maps $\Frob_q^\mu\otimes\id$ for all integers $\mu$ modulo~$(d_\Fp)$. For every $\nu\ge0$ it follows that
\begin{eqnarray*}
1\otimes \epsilon_\Fp^{q^\nu}
\ =\ (1\otimes \epsilon_\Fp)^{q^\nu}
&=& \prod_{\mu\,{\rm mod}\,(d_\Fp)} (\Frob_q^\mu\otimes\id)(\proj_\Fp\otimes\id)(e)^{q^\nu} \\[3pt]
&=& \prod_{\mu\,{\rm mod}\,(d_\Fp)} (\Frob_q^{\mu+\nu}\otimes\id)(\proj_\Fp\otimes\id)(\id\otimes\Frob_q^\nu)(e) \\[3pt]
&=& \prod_{\mu'\,{\rm mod}\,(d_\Fp)} (\Frob_q^{\mu'}\otimes\id)(\proj_\Fp\otimes\id)(\sigma^\nu(e))
\end{eqnarray*}
Here the right hand side is $(\proj_\Fp\otimes\id)(\sigma^\nu(e))$ times some other element of $\ICk\otimes_{\BFq}\kern-1ptR$. Therefore $(\proj_\Fp\otimes\id)(\sigma^\nu(e))$ becomes a unit in $\ICk\otimes_{\BFq}\kern-1ptR[\epsilon_\Fp^{-1}]$. This implies that $\sigma^\nu(e)$ becomes a unit in ${A_\Fp\complot_{\BFq}\kern-1ptR[\epsilon_\Fp^{-1}]}$. By the construction \eqref{EAReDef} of $\CE_{A,R,e}$ it follows that the inclusion $A\otimes_{\BFq}\kern-1ptR \into A_\Fp\complot_{\BFq}\kern-1ptR[\epsilon_\Fp^{-1}]$ extends to an inclusion 
\UseTheoremCounterForNextEquation
\begin{equation}\label{GlobLocICSlope0ExtendsMap}
\CE_{A,R,e}\ \into\ A_\Fp\complot_{\BFq}\kern-1ptR[\epsilon_\Fp^{-1}]
\ =\ \EAp{R[\epsilon_\Fp^{-1}]}
\end{equation}
In particular it therefore extends to an inclusion $\CE_{A,R,e} \into \ICF\complot_{\BFq}\kern-1pt R[\epsilon_\Fp^{-1}]$. Thus the proof of Proposition \ref{GlobICFpGoodRed} works with this choice of~$\epsilon_\Fp$. 

Finally, the base change of $\CMod$ via \eqref{GlobLocICSlope0ExtendsMap} is a dualizable left $\bigl(A_\Fp\complot_{\BFq}\kern-1.5ptR[\epsilon_\Fp^{-1}]\bigr)[\tau]$-module. By Definition \ref{LocICRSlopeDef} this constitutes a lattice in $\CMod_\Fp$ that makes $\CMod_\Fp$ pure of slope~$0$.
\end{Proof}

\medskip
In particular Proposition \ref{GlobLocICSlope0Extends} implies that for all but finitely many places $\Fp$ of~$F$, the local isocrystal $M_\Fp$ over $K$ associated to $M$ is pure of slope~$0$. 

\medskip
Finally, the following result describes the points of $\Spec R[\epsilon_\Fp^{-1}]$ in terms of the divisor~$D$:

\begin{Prop}\label{pxGoodIff}
For every prime $\Fp$ of~$A$ and the element $\epsilon_\Fp$ as in \eqref{epDef}, a point $x\in\Spec R$ satisfies $\epsilon_\Fp(x)\not=0$ if and only if 
$$(\Fp\times x)\cap D\ =\ \varnothing.$$
\end{Prop}

\begin{Proof}
By definition $D$ is the zero locus of the element $e\in A\otimes_{\BFq}\kern-1ptR$. Letting $i_x\colon R\to k_x$ denote the natural homomorphism to the residue field, as before, it follows that 
$(\Fp\times x)\cap D$ is the zero locus of the element $(\proj_\Fp\otimes i_x)(e)\in \ICk\otimes_{\BFq}\kern-1ptk_x$. Thus $(\Fp\times x)\cap D = \varnothing$ if and only if $(\proj_\Fp\otimes i_x)(e)$ is a unit in $\ICk\otimes_{\BFq}\kern-1ptk_x$. This is equivalent to its norm via $\ICk\otimes_{\BFq}\kern-1ptk_x \to \BFq\otimes_{\BFq}\kern-1ptk_x \cong k_x$ being a unit, or equivalently being non-zero. But by \eqref{epDef} this norm is equal to $(\id\otimes i_x)(1\otimes \epsilon_\Fp) = 1\otimes \epsilon_\Fp(x)$. Thus $(\Fp\times x)\cap D = \varnothing$ if and only if $\epsilon_\Fp(x)\not=0$, and the lemma follows.
\end{Proof}

\section{Compatible system of Weil representations}
\label{CompSys}

We keep all notation of the preceding section. Since $K$ is finitely generated over~$\BFq$, the algebraic closure $k$ of $\BFq$ in $K$ is a finite extension of $\BFq$. As in Section \ref{WeilRep} choose a separable closure $K^\sep$ of $K$ with Galois group $\GK$ and Weil group~$\WK$, and let $\kalg$ be the algebraic closure of $\BFq$ in~$K$. We are interested in the representations of $\WK$ on $U(M_\Fp)$ for all places $\Fp$ of~$F$, viewed as free modules over $\smash{\EFp{\kalg}} = \ICF\complot_{\BFq}\kern-1pt\kalg$.

For this we now assume that $R$ is normal and consider any \emph{closed} point $x\in X=\Spec R$. As in Section~\ref{RedUM} choose a henselization and a strict henselization of $R$ at $x$, with associated inertia and decomposition group $I_x\mathrel{\triangleleft} D_x\subset \GK$. 
Since $R$ is finitely generated over~$\BFq$, the residue field $k_x$ is a finite extension of~$\BFq$. Let $\Frob_x\in\WK\cap D_x$ be any element that maps to the arithmetic Frobenius element in $W_{k_x}$. Let $\charact_{\CMod_x} \in F[X]$ be the characteristic polynomial associated to $\CMod_x$ by Construction \ref{GlobCharPol}. 

\medskip
Then for all places $\Fp$ of~$F$, the $U(M_\Fp)$ form a \emph{compatible system of representations of $\WK$ over $\smash{\EFp{\kalg}}$} in the following sense:

\begin{Thm}\label{CompatibleSystemU}
For every place $\Fp$ of $F$ and every closed point $x$ of~$X$ with $\epsilon_\Fp(x)\not=0$, the inertia group $I_x$ acts trivially on $U(M_\Fp)$, and the geometric Frobenius element $\Frob_x^{-1}$ has characteristic polynomial $\charact_{\CMod_x}$ on $U(M_\Fp)$ over $\smash{\EFp{\kalg}}$. In particular this characteristic polynomial has coefficients in~$F$ and is independent~of~$\Fp$.
\end{Thm}

\begin{Proof}
By Proposition \ref{GlobICFpGoodRed} the local isocrystal $M_\Fp$ associated to $M$ comes from an $\ICF$-isocrystal $\CMod_\Fp$ over $R[\epsilon_\Fp^{-1}]$ which possesses a slope filtration, and since $\epsilon_\Fp(x)\not=0$, there is a natural isomorphism $i_x^*\CMod_\Fp \cong \CMod_{x,\Fp}$. From Theorem \ref{LocICRGoodFiltThmU} we thus obtain $\WK\cap D_x$-equivariant isomorphisms of $\smash{\EFp{\kalg}}$-modules
$$U(M_\Fp)\ \cong\ U(i_x^*\CMod_\Fp)\ \cong\ U(\CMod_{x,\Fp}).$$
Since the inertia group $I_x$ acts trivially on the right hand side, it thus also acts trivially on the left. The assertion about the characteristic polynomial of $\Frob_x^{-1}$ is now a direct consequence of Corollary \ref{CharPolGalois1}.
\end{Proof}

\medskip
Likewise, for almost all places $\Fp$ of~$F$, the $\Fp$-adic Tate modules $T_0(M_\Fp)$ associated to $M$ by \eqref{WMSlope0} form a compatible system of representations of $\GK$ over $\ICF$ in the following sense:

\begin{Thm}\label{CompatibleSystemT}
For every prime $\Fp$ of~$A$ and every closed point $x$ of~$X$ with $(\Fp\times x)\cap D = \varnothing$, the inertia group $I_x$ acts trivially on $T_0(M_\Fp)$, and the geometric Frobenius element $\Frob_x^{-1}$ has characteristic polynomial $\charact_{\CMod_x}$ on $T_0(M_\Fp)$ over~$\ICF$. In particular this characteristic polynomial has coefficients in~$F$ and is independent~of~$\Fp$.
\end{Thm}

\begin{Proof}
Combine Theorem \ref{CompatibleSystemU} with Proposition \ref{pxGoodIff} and the isomorphism \eqref{WMSmallWeilRepIsomSlope0}.
\end{Proof}

\section{Application to \texorpdfstring{$A$}{A}-motives} 
\label{AMot}

Recall that $F$ is the function field of a smooth connected projective algebraic curve $C$ over~$\BFq$. In this section we let $A\subset F$ denote the subring of all functions that are regular outside a fixed closed point $\infty\in C$. Then the ideal class group of $A$ is finite, so for every maximal ideal $\Fp\subset A$ we can choose an element $a_\Fp\in A$ with $(a_\Fp)=\Fp^n$ for some $n\ge1$. This will help our constructions when $\Fp$ is not a principal ideal.

Consider an $\BFq$-algebra $R$ and view $A\otimes_{\BFq}\kern-1ptR$ as a difference ring via the partial Frobenius endomorphism $\sigma(a\otimes\xi):= a\otimes\xi^q$. Fix a morphism $c\colon\Spec R\to\Spec A \subset C$ corresponding to an $\BFq$-algebra homomorphism $c^*\colon A\to R$, and let $J_{c,R}$ denote the kernel of the ring homomorphism $A\otimes_{\BFq}\kern-1ptR\to R$, ${a\otimes\xi\mapsto c^*(a)\xi}$. In analogy to Hartl--Juschka \cite[Def.\,2.3.1]{HartlJuschka} we say:

\begin{Def}\label{AMotDef}
An \emph{(effective) $A$-motive $\CM$ of characteristic~$c^*$ over~$R$} is a left module over $(A\otimes_{\BFq}\kern-1ptR)[\tau]$ that is finitely generated projective over $A\otimes_{\BFq}\kern-1ptR$, such that $\coker(\taulin_\CM)$ is annihilated by a power of~$J_{c,R}$. 
\end{Def}


First consider an $A$-motive $M$ of characteristic $c^*$ over a field~$K$. Then by base change via the embedding $A\otimes_{\BFq}\kern-1ptK \into \Quot(F \otimes_{\BFq}\kern-1pt K) = \EFK$, we obtain a left $\EFK[\tau]$-module $M_F$ that is finitely generated projective over $\EFK$.

\begin{Prop}\label{AMotIsoCrys}
The module $M_F$ is an $F$-isocrystal over $K$.
\end{Prop}

\begin{Proof}
The structure morphism of $M_F$ is surjective because $J_{c,K}$ is a maximal ideal of $A\otimes_{\BFq}\kern-1pt K$ and so generates the ideal $(1)$ of the ring $\Quot(A\otimes_{\BFq}\kern-1.5ptK) = \EFK$. Next, pick a transcendental element $a \in A$. Then $M_F$ and $\sigma^*(M_F)$ are vector spaces of the same dimension over the field $K(a)$, and so the surjection $\taulin_{M_F}$ is an isomorphism.
\end{Proof}

\begin{Prop}\label{AMotCharPolA}
If $K$ is finite, the characteristic polynomial $\charact_{M_F}$ has coefficients~in~$A$. 
\end{Prop}

\begin{Proof}
In this case $M$ is a finitely generated projective $A$-submodule of $M_F$ such that ${F\otimes_AM =M_F}$. Since it is also preserved by the endomorphism~$\tau_{M_F}$, the claim follows from the formula \eqref{BPLemmaGlobal}.
\end{Proof}

\medskip
Next consider the prime ideal $\wp := \ker(c^*)$ of~$A$, which may be zero or maximal. 

\begin{Prop}\label{AMotFpSlope}
For every place $\Fp$ of $F$ the associated $\ICF$-isocrystal $M_{F,\Fp}$ over~$K$
\begin{enumerate}\StatementLabels
\item\label{AMotFpSlopeA}
... is pure of slope $0$ if $\Fp\not=\wp,\infty$.
\item\label{AMotFpSlopeB}
... has all slopes $\ge0$ if $\Fp=\wp$.
\end{enumerate}
\end{Prop}

\begin{Proof}
For $\Fp\not=\infty$ the ring $\ICO$ is simply the $\Fp$-adic completion of~$A$. By base change via $A\otimes_{\BFq}\kern-1.5ptR \into \ICO\complot_{\BFq}\kern-1.5ptK = \EAp{K}$, from $M$ we therefore obtain a lattice $L_\Fp\subset M_{F,\Fp}$ that satisfies $\tau L_\Fp\subset L_\Fp$. Thus \ref{AMotFpSlopeB} follows from Proposition \ref{SlopesLowBound}. In the case \ref{AMotFpSlopeA} the maximal ideal $J_{c,K}$ does not lie above the ideal $\Fp\subset A$ and thus generates the ideal $(1)$ in $A_\Fp\complot_{\BFq}\kern-1.5ptK$. Consequently $\coker(\taulin_{L_\Fp}) = 0$, and so $M_{F,\Fp}$ is pure of slope~$0$.
\end{Proof}

%

\begin{Rem}\label{AMotPure}\rm
As a direct consequence of Definition \ref{ICSlopeDef}, the $F_\infty$-isocrystal $M_{F,\infty}$ is pure of slope $\alpha$ if and only if $M$ is pure of weight $-\alpha$ in the sense of \cite[\S1.10]{AndersonT} and \cite[Def.\,2.3.8]{HartlJuschka}. In that case \eqref{WMSmallWeilRepIsom} yields an isomorphism
$$U(M_{F,\infty})\ \cong\ 
\smash{T_\alpha(M_{F,\infty})\centrot_{\raisebox{-2pt}{$\scriptstyle D_\alpha$}} \SA{\kalg}.}$$
\end{Rem}

\medskip
Next, for any place $\Fp\not=\wp,\infty$ of $F$, combining Proposition \ref{AMotFpSlope}~\ref{AMotFpSlopeA} with \eqref{WMSlope0} yields the usual \emph{rational $\Fp$-adic Tate module} of $M$ (compare \cite[\S2.3.5, p.70]{HartlJuschka}):
\UseTheoremCounterForNextEquation
\begin{equation}\label{T0MAMotDef}
T_0(M_{F,\Fp})
\ \cong\ \bigl(\ICF\complot_A (M\otimes_K\kern-.5pt\Kalg)\bigr)^\tau.
\end{equation}
Likewise, for every place $\Fp$ of $F$ we now also have the associated module $U(M_{F,\Fp})$. In the remainder of this section we want to make more precise how these form compatible systems of representations following Theorems \ref{CompatibleSystemU} and \ref{CompatibleSystemT}.

\medskip
For this let us return to $A$-motives over an arbitrary ring~$R$. As a direct consequence of Definition \ref{AMotDef}, for any $\BFq$-algebra homomorphism $f\colon R\to R'$, by base change any $A$-motive $\CM$ of characteristic~$c^*$ over~$R$ yields an $A$-motive $f^*\CM$ of characteristic $f\circ c^*$ over~$R'$. Conversely we have:

\begin{Prop}\label{AMotKR}
For any $A$-motive $M$ of characteristic~$c^*$ over a field $K$, there exists a finitely generated $\BFq$-subalgebra $R\subset K$ such that $c$ extends to a morphism ${\Spec R\to\Spec A}$, again denoted~$c$, and $M$ arises by base change from an $A$-motive of characteristic $c^*$ over~$R$.
\end{Prop}

\begin{Proof}
Since $M$ is a finitely generated projective $A\otimes_{\BFq}\kern-1ptK$-module, it is a direct summand of $(A\otimes_{\BFq}\kern-1ptK)^{\oplus r}$ for some integer $r\ge0$. It is thus the image of an idempotent matrix $I\in\Mat_{r\times r}(A\otimes_{\BFq}\kern-1ptK)$. For some finitely generated $\BFq$-subalgebra $R\subset K$ this matrix has coefficients in $A\otimes_{\BFq}\kern-1ptR$, and the image $\CM$ of $(A\otimes_{\BFq}\kern-1ptR)^{\oplus r}$ under $I$ is then a finitely generated projective $A\otimes_{\BFq}\kern-1ptR$-submodule of~$M$, whose base change to $K$ yields back~$M$.

Next, after enlarging $R$ we may assume that the images under $\tau$ of finitely many generators of $\CM$ again lie in~$\CM$. Then $\CM$ is a left $(A\otimes_{\BFq}\kern-1ptR)[\tau]$-module. After enlarging $R$ again we may also assume that $c^*(A)\subset R$, so that $c$ extends to a morphism ${\Spec R\to\Spec A}$. Then the cokernel of $\taulin_\CM$ is a finitely generated $A\otimes_{\BFq}\kern-1ptR$-module whose tensor product with $K$ over $R$ is annihilated by some power $J_{c,R}^n$. After inverting finitely many elements of $R$ we can thus achieve that $\coker(\taulin_\CM)$ itself is annihilated by  $J_{c,R}^n$, and then $\CM$ has the desired properties.
\end{Proof}

\medskip
For the following we thus fix an $A$-motive $\CM$ of characteristic $c^*$ defined over an integral domain $R$ which is finitely generated over~$\BFq$. As before we let $i\colon R\into K$ denote the embedding into the quotient field. We abbreviate $X := \Spec R$, and for each point $x\in X$ we let $i_x\colon R\to k_x$ denote the natural homomorphism to the residue field at~$x$. Again we consider the prime ideal $\wp := \ker(c^*)$ of~$A$, which may be zero or maximal. We can then make the results of Section~\ref{GoodRedGlob} more precise, as follows:

\begin{Prop}\label{AMotRGoodP0Max}
If $\wp$ is a maximal ideal, let $A'\subset F$ denote the subring of all functions that are regular outside $\infty$ and~$\wp$ and set $e:=1$. Then by base change $\CM$ yields a dualizable left $\CE_{A',R,1}[\tau]$-module $\CMod$ that satisfies the conditions in Propositions \ref{GlobICFpGoodRedAlg} and \ref{AllAssOnAReD}.
\end{Prop}

\begin{Proof}
By assumption $\wp$ is the only point in the image of $c\colon X\to\Spec A$; hence the inverse image under $c$ of the open subscheme $\Spec A' \subset\Spec A$ is empty. Let $\CMod$ denote the left $(A'\otimes_{\BFq}\kern-1ptR)[\tau]$-module that is obtained from $\CM$ by base change. Then the condition on $\coker(\taulin_\CM)$ in Definition \ref{AMotDef} implies that $\taulin_{\CMod}$ is surjective. Since $\CMod$ is a finitely generated projective $A'\otimes_{\BFq}\kern-1ptR$-module, and $\taulin_{\CMod}$ is generically an isomorphism by Proposition \ref{AMotIsoCrys}, it follows that $\taulin_{\CMod}$ is an isomorphism. With the non-zerodivisor $e=1$ we then have $A'\otimes_{\BFq}\kern-1ptR = \CE_{A',R,1}$ in~\eqref{EAReDef}, and so $\CMod$ is a dualizable left $\CE_{A',R,1}[\tau]$-module. Thus $\CMod$ satisfies the condition in Proposition \ref{GlobICFpGoodRedAlg}, and with the zero divisor $D:=\varnothing$ of $e=1$ it satisfies the conditions in Proposition \ref{AllAssOnAReD}.
\end{Proof}

\begin{Prop}\label{AMotRFpICP0Max}
If $\wp$ is a maximal ideal, for every place $\Fp$ of~$F$ the $A$-motive $\CM$ gives rise to an $\ICF$-isocrystal $\CM_{F,\Fp}$ over~$R$, and that is pure of slope~$0$ if $\Fp\not\in\{\wp,\infty\}$.
\end{Prop}

\begin{Proof}
Let $A'$ and $\CMod$ be as in Proposition \ref{AMotRGoodP0Max}. Then for any place $\Fp$ of~$F$ the embedding $A'\into\ICF$ induces an embedding of difference rings
$$\CE_{A',R,1}\ =\ A'\otimes_{\BFq}\kern-1ptR\ \longinto\ \ICF\complot_{\BFq}\kern-1.5ptR\ =\ \EFp{R}.$$
By base change $\CMod$ therefore yields the desired $\ICF$-isocrystal $\CM_{F,\Fp}$ over~$R$. If in addition $\Fp\not\in\{\wp,\infty\}$, the embedding $A'\into\ICO$ induces an embedding of difference rings
$$\CE_{A',R,1}\ =\ A'\otimes_{\BFq}\kern-1ptR\ \longinto\ \ICO\complot_{\BFq}\kern-1.5ptR\ =\ \EAp{R}.$$
By base change $\CMod$ therefore yields a dualizable left $\EAp{R}[\tau]$-module and hence a lattice that makes $\CM_{F,\Fp}$ pure of slope~$0$.
\end{Proof}


\begin{Prop}\label{AMotRGoodP0Zero}
If $\wp$ is the zero ideal, take any element $a\in A$ that is transcendental over~$\BFq$. Then $e := a\otimes1-1\otimes c^*(a)$ is a non-zerodivisor in $A\otimes_{\BFq}\kern-1ptR$, and by base change $\CM$ yields a dualizable left $\CE_{A,R,e}[\tau]$-module $\CMod$ that satisfies the conditions in Propositions \ref{GlobICFpGoodRedAlg} and~\ref{AllAssOnAReD}.
\end{Prop}

\begin{Proof}
As $a$ is transcendental over~$\BFq$, the ring $A$ is a finitely generated projective module over the subring $\BFq[a]$, which is isomorphic to the polynomial ring $\BFq[t]$. Thus $A\otimes_{\BFq}\kern-1ptR$ is a finitely generated projective module over $\BFq[a]\otimes_{\BFq}\kern-1ptR \cong R[t]$. Under this isomorphism the element $e = a\otimes1-1\otimes c^*(a)$ corresponds to the polynomial $t-c^*(a) \in R[t]$. As that is a non-zerodivisor in $R[t]$, it is also a non-zerodivisor in the flat $R[t]$-algebra $A\otimes_{\BFq}\kern-1ptR$, proving the first statement. Note also that by the construction of $J_{c,R}$ we have $e \in J_{c,R}$. 

Now let $\CMod$ denote the $\CE_{A,R,e}$-module obtained by base change from $\CM$ via the embedding $A\otimes_{\BFq}\kern-1.5ptR \into\CE_{A,R,e}$. This module inherits an action of $\tau$ for which $\coker(\taulin_{\CMod})$ is annihilated by a power of $e$ and is therefore zero. Since $\CMod$ is a finitely generated projective $\CE_{A,R,e}$-module, and $\taulin_{\CMod}$ is generically an isomorphism by Proposition \ref{AMotIsoCrys}, it follows that $\taulin_{\CMod}$ is an isomorphism. Thus $\CMod$ is a dualizable left $\CE_{A,R,e}[\tau]$-module, as required in Proposition \ref{GlobICFpGoodRedAlg}.

Next, the condition \ref{AllAssOnAReD} \ref{AllAssOnAReDA} holds by assumption. Finally, since $a$ is transcendental over~$\BFq$ and $c$ is injective, the image $c^*(a) \in R$ is again transcendental over~$\BFq$. This element therefore defines a non-constant morphism $c^*(a)\colon X\to\BA^{\kern-1pt1}$. Taking the transpose of its graph in $\BA^{\kern-1pt1}\times_{\BFq}\kern-1pt X$ and pulling it back under the finite projection
$$(\Spec A)\times_{\BFq}\kern-1pt X\ \longonto\ 
\bigl(\Spec\BFq[a]\bigr)\times_{\BFq}\kern-1pt X\ \cong\ \BA^{\kern-1pt1}\times_{\BFq}\kern-1.5pt X$$
yields precisely the zero divisor $D$ of~$e$. The fact that this is the inverse image of a transpose graph then implies that $D$ contains no irreducible component of ${(\Spec A)\times_{\BFq}\kern-1pt x}$ for any point $x\in X$. Moreover, the fact that the morphism $c^*(a)$ is non-constant implies that $D$ contains no irreducible component of $\Fp\times_{\BFq}\kern-1pt X$ for any point $\Fp\in\Spec A$. Thus $D$ satisfies the condition \ref{AllAssOnAReD} \ref{AllAssOnAReDeD}, and we are done.
\end{Proof}

\begin{Prop}\label{AMotRFpICP0Zero}
If $\wp$ is the zero ideal, consider any place $\Fp$ of~$F$.
\begin{enumerate}\StatementLabels
\item\label{AMotRFpICP0ZeroA}
If $\Fp\not=\infty$, the $A$-motive $\CM$ gives rise to an $\ICF$-isocrystal $\CM_{F,\Fp}$ over $R[c^*(a_\Fp)^{-1}]$ that is pure of slope~$0$.
\item\label{AMotRFpICP0ZeroB}
The $A$-motive $\CM$ gives rise to an $F_\infty$-isocrystal $\CM_{F,\infty}$ over~$R$.
\end{enumerate}
\end{Prop}

\begin{Proof}
In \ref{AMotRFpICP0ZeroA} recall that $a_\Fp\in A$ satisfies $(a_\Fp)=\Fp^n$ for some $n\ge1$; hence $a_\Fp$ is transcendental over~$\BFq$. Set $e := a_\Fp\otimes1-1\otimes c^*(a_\Fp)$ and let $\CMod$ be as in Proposition \ref{AMotRGoodP0Zero} for $a=a_\Fp$. Since $a_\Fp$ lies in $\Fp$, and $c^*(a_\Fp)$ is a unit in $R[c^*(a_\Fp)^{-1}]$, for any $\nu\ge0$ the element
$$\sigma^\nu(e)\ =\ a_\Fp\otimes1-1\otimes c^*(a_\Fp)^{q^\nu}\ \in\ \ICO\complot_{\BFq}\kern-1ptR[c^*(a_\Fp)^{-1}]$$
is a unit. By the definition of $\CE_{A,R,e}$ we therefore obtain an embedding of difference rings
$$\CE_{A,R,e}\ =\ S_e^{-1}(A\otimes_{\BFq}\kern-1ptR)\ \longinto\ \ICO\complot_{\BFq}\kern-1ptR[c^*(a_\Fp)^{-1}]\ =\ \EAp{R[c^*(a_\Fp)^{-1}]}.$$
By base change $\CMod$ therefore yields a dualizable left $\EAp{R[c^*(a_\Fp)^{-1}]}[\tau]$-module, or equivalently a lattice in an $\ICF$-isocrystal over $R[c^*(a_\Fp)^{-1}]$ that is pure of slope~$0$, as desired.

In \ref{AMotRFpICP0ZeroB} take any element $a\in A$ that is transcendental over~$\BFq$ and let $e := a\otimes1-1\otimes c^*(a)$ and $\CMod$ be as in Proposition \ref{AMotRGoodP0Zero}. For any $\nu\ge0$ consider the element
$$\sigma^\nu(e)\ =\ a\otimes1-1\otimes c^*(a)^{q^\nu}
\ =\ (a\otimes1)\bigl(1\otimes1-a^{-1}\otimes c^*(a)^{q^\nu}\bigr)
\ \in\ F_\infty\complot_{\BFq}\kern-1ptR.$$
Here $a\otimes1$ is a unit in $F_\infty\complot_{\BFq}\kern-1ptR$, and as $a^{-1}$ lies in the maximal ideal of $A_\infty$, the element $1\otimes1-a^{-1}\otimes c^*(a)^{q^\nu}$ is also a unit in $F_\infty\complot_{\BFq}\kern-1ptR$. We conclude that $\sigma^\nu(e)$ is a unit in $F_\infty\complot_{\BFq}\kern-1ptR$. By the definition of $\CE_{A,R,e}$ we therefore obtain an embedding of difference rings
$$\CE_{A,R,e}\ =\ S_e^{-1}(A\otimes_{\BFq}\kern-1ptR)\ \longinto\ F_\infty\complot_{\BFq}\kern-1ptR\ =\ \CE_{F_\infty,R}.$$
By base change $\CMod$ therefore yields the desired $F_\infty$-isocrystal over~$R$.
\end{Proof}

\medskip
Now we assume in addition that $R$ is normal. We can then make Theorems \ref{CompatibleSystemT} and \ref{CompatibleSystemU} more precise, as follows. Abbreviate $M := i^*\CM$, and for every $x\in X$ let $\CM_x := i_x^*\CM$ denote the $A$-motive over $k_x$ obtained by base change. The next theorem generalizes Gardeyn \cite[Prop.\,3.3]{GardeynT} to include the case of special characteristic:

\begin{Thm}\label{CompatibleSystemTAMot}
For every prime $\Fp\not=\wp$ of~$A$ and every closed point $x$ of~$X$ with $c(x)\not=\Fp$, the inertia group $I_x$ acts trivially on $T_0(M_{F,\Fp})$, and the geometric Frobenius element $\Frob_x^{-1}$ has characteristic polynomial $\charact_{\CM_{x,F}}$ on $T_0(M_{F,\Fp})$ over~$\ICF$. In particular this characteristic polynomial has coefficients in~$A$ and is independent~of~$\Fp$.
\end{Thm}


\begin{Proof}
If $\wp$ is a maximal ideal, by Proposition \ref{AMotRGoodP0Max} the conditions in Propositions \ref{GlobICFpGoodRedAlg} and \ref{AllAssOnAReD} are satisfied with $A'$ in place of $A$ and with $e=1$ and $D=\varnothing$, and by construction we have $\CMod_x = \CM_{x,F}$. Since the primes of $A'$ correspond precisely to the primes $\Fp\not=\wp$ of~$A$, by Theorem \ref{CompatibleSystemT} the desired statements thus hold for all those and all closed points $x$ of~$X$.

\medskip
If $\wp$ is the zero ideal, fix a prime $\Fp$ of~$A$ and a closed point $x\in X$ with $c(x)\not=\Fp$. Letting $\Fm_x\subset R$ denote the maximal ideal associated to~$x$, the point $c(x)$ then corresponds to the maximal ideal $(c^*)^{-1}(\Fm_x)\subset A$. As that is different from the maximal ideal $\Fp\subset A$, there exists an element $a\in(c^*)^{-1}(\Fm_x)$ with $a\equiv1\bmod\Fp$. This element is then transcendental over~$\BFq$. Also, the congruences $a\equiv1\bmod\Fp$ and  $c^*(a)\equiv0\bmod\Fm_x$ imply that 
$$e\ :=\ a\otimes1-1\otimes c^*(a)\ \equiv\ 1\otimes1 \ \bmod\ (\Fp\otimes_{\BFq}\kern-1ptR+A\otimes_{\BFq}\kern-1pt\Fm_x).$$
Thus $e$ maps to the identity element in the tensor product of the residue fields $k_\Fp\otimes_{\BFq}\kern-1ptk_x$. Letting $D\subset(\Spec A)\times_{\BFq}X$ denote the zero divisor of~$e$, this means that $(\Fp\times x)\cap D = \varnothing$. Also by construction we then have $\CMod_x = \CM_{x,F}$. By Proposition \ref{AMotRGoodP0Zero} and Theorem \ref{CompatibleSystemT} the desired statements thus hold for $\Fp$ and~$x$, as desired. 
 
Finally, the fact that $\charact_{\CM_{x,F}}$ has coefficients in~$A$ follows from
Proposition \ref{AMotCharPolA}. 
\end{Proof}


\begin{Thm}\label{CompatibleSystemUAMot}
For every place $\Fp$ of $F$ there exists a divisor $D_\Fp\subset X$, such that $D_\Fp = c^{-1}(\Fp)$ whenever $\Fp\not\in\{\wp,\infty\}$, and that for every closed point $x$ of~$X\smallsetminus D_\Fp$, the inertia group $I_x$ acts trivially on $U(M_{F,\Fp})$, and the geometric Frobenius element $\Frob_x^{-1}$ has characteristic polynomial $\charact_{\CM_{x,F}}$ on $U(M_{F,\Fp})$ over $\smash{\EFp{\kalg}}$. In particular this characteristic polynomial has coefficients in~$A$ and is independent~of~$\Fp$.
\end{Thm}

\begin{Proof}
For $\Fp\not\in\{\wp,\infty\}$ this follows from Theorem \ref{CompatibleSystemTAMot} and the isomorphism \eqref{WMSmallWeilRepIsomSlope0}. For $\Fp\in\{\wp,\infty\}$ it follows from Theorem \ref{CompatibleSystemU} applied to the $\CMod$ from Proposition \ref{AMotRGoodP0Max} or \ref{AMotRGoodP0Zero}.
\end{Proof}

\section{Application to Drinfeld modules}
\label{Drin}

As in the preceding section we let $A\subset F$ denote the subring of all functions that are regular outside a fixed closed point $\infty\in C$. For any non-zero element $a\in A$ we set 
\UseTheoremCounterForNextEquation
\begin{equation}\label{DegADef}
\deg_A(a)\ :=\ \dim_{\BFq}(A/Aa)\ \in\ \BZ^{\ge0}.
\end{equation}
Consider an $\BFq$-algebra $R$. Then a polynomial $\psi=\sum_{i=0}^d \psi_i\tau^i \in R[\tau]$ with all $\psi_i\in R$, whose highest coefficient $\psi_d$ is a unit in~$R$, is called \emph{elliptic of degree~$d$}.

\medskip
Consider an integer $r\ge1$. By a \emph{(standard) Drinfeld $A$-module of rank $r$ over~$R$} we mean an $\BFq$-algebra homomorphism $\phi\colon A\to R[\tau]$ such that the image $\phi_a$ of each non-zero $a\in A$ is elliptic of degree $r\deg_A(a)$. 
(There is a more general notion of Drinfeld $A$-modules over a base, but Zariski locally any such Drinfeld $A$-module is isomorphic to one of the special form here.)
For any such~$\phi$, the $\BFq$-algebra homomorphism 
\UseTheoremCounterForNextEquation
\begin{equation}\label{CharPhiDef}
A\longto R,\ \ a\mapsto c^*(a) := d\phi_a
\end{equation}
is called the \emph{characteristic of~$\phi$.}

\medskip
{}From now on we assume that $R$ is a noetherian integral domain with quotient field~$K$. Then the kernel $\wp$ of $c^*$ is either zero or a maximal ideal of~$A$. In the latter case recall that $a_\wp\in A$ satisfies $(a_\wp)=\wp^n$ for some $n\ge1$. Thus there exists a unique integer $1\le h\le r$ such that $\phi_{a_\wp} = \epsilon_\wp\tau^{h\deg_A(a_\wp)}+(\hbox{higher terms})$ with $\epsilon_\wp\in R\setminus\{0\}$. The integer~$h$ is independent of the choice of $a_\wp$ and called the \emph{height of~$\phi$ over~$K$.} If the coefficient $\epsilon_\wp$ is a unit in~$R$, we say that \emph{$\phi$ has constant height~$h$ over~$R$}. In particular $\phi$ acquires constant height $h$ over $R[\epsilon_\wp^{-1}]$.

\medskip
For any $\phi$ as above we endow $\CM:=R[\tau]$ with the structure of an $A\otimes_{\BFq}\kern-1ptR$-module by setting $(a\otimes\xi)m := \xi m\phi_a$. Left multiplication by $\tau$ turns this into a left module over $(A\otimes_{\BFq}\kern-1ptR)[\tau]$. 

\begin{Prop}[Drinfeld]\label{DrinAMot}
This is an $A$-motive of rank $r$ and characteristic $c^*$ over~$R$.
\end{Prop}

\begin{Proof}
The module $\CM$ defines a quasicoherent sheaf $\CF$ on the scheme ${\Spec(A \otimes_{\BFq}\kern-1pt R)}$, and $\taulin_\CM\colon\CM\to\CM$ corresponds to a homomorphism $\taulin_\CF\colon \sigma^*\CF\to\CF$. Drinfeld \cite[Prop.\,3]{DrinfeldComm} proves that this sheaf is locally free of rank ~$r$ and that $\coker(\taulin_\CF)$ is supported on the graph of~$c$. Thus $\CM$ satisfies the conditions in Definition~\ref{AMotDef}.
\end{Proof}

\medskip
For every place $\Fp$ of $F$ this $A$-motive yields an $\ICF$-isocrystal as follows. 

\begin{Prop}\label{DrinFpIC}
For any place $\Fp\not=\wp,\infty$ of $F$, the $A$-motive $\CM$ gives rise to an $\ICF$-isocrystal $\CM_{F,\Fp}$ over $R[c^*(a_\Fp)^{-1}]$ that is pure of slope~$0$.
\end{Prop}

\begin{Proof}
Direct consequence of Propositions \ref{AMotRFpICP0Max} and~\ref{AMotRFpICP0Zero}.
\end{Proof}


\begin{Prop}\label{DrinFinftyIC}
For $\Fp=\infty$ the $A$-motive $\CM$ gives rise to an $F_\infty$-isocrystal $\CM_{F,\infty}$ over~$R$ that is pure of slope~$\minusoneoverr$.
\end{Prop}

\begin{Proof}
By construction $\infty$ represents a closed point on~$C$ such that $\Spec A = C\setminus\{\infty\}$. Abbreviating $X:=\Spec R$, the scheme ${\Spec(A \otimes_{\BFq}\kern-1pt R)}$ is thus the complement of the divisor $\infty_X := {\infty\times_{\BFq}\kern-1ptX}$ within $C_X := {C\times_{\BFq}\kern-1ptX}$. Moreover $\sigma\colon {A \otimes_{\BFq}\kern-1pt R}\to {A \otimes_{\BFq}\kern-1pt R}$ corresponds to a morphism $\sigma\colon C_X\to C_X$ which is the product of the identity endomorphism of $C$ and of the absolute $q$-Frobenius endomorphism of~$X$. Let $d_\infty$ denote the degree of $\infty$ over~$\BFq$.


The module $\CM$ corresponds to a locally free sheaf $\CF$ of rank $r$ on $C_X\smallsetminus\infty_X$ endowed with a homomorphism $\taulin_\CF\colon \sigma^*\CF\to\CF$. Drinfeld \cite[Prop.3]{DrinfeldComm} constructs an increasing family $(\CF_n)_{n\in\BZ}$ of locally free sheaves of rank $r$ on $C_X$ extending~$\CF$, and a compatible system of morphisms $\sigma^*\CF_n \to \CF_{n+1}$ extending~$\taulin_\CF$, such that for every $n$ we have $\CF_{n + rd_\infty} = \CF_n(\infty_X)$ and the induced morphism $\sigma^*(\CF_n/\CF_{n-1}) \longto \CF_{n+1}/\CF_n$ is an isomorphism.

Recall that $A_\infty$ denotes the valuation ring in $F_\infty$ with a uniformizer $z_\infty \in A_\infty$. The pullbacks of the $\CF_n$ under the natural morphism $\Spec(A_\infty\complot_{\BFq}\kern-1ptR) \to C_X$ then correspond to finitely generated projective $\CE_{A_\infty,R}$-modules $\CL_n$ and hence to lattices in $\CM_{F,\infty} = F_\infty\complot_{A}\CM$. These form an increasing family satisfying $\tau\CL_n\subset\CL_{n+1}$ and $\CL_{n + rd_\infty} = z_\infty^{-1}\CL_n$, such that $\taulin_\CM$ induces isomorphisms $\sigma^*(\CL_n/\CL_{n-1}) \isoto \CL_{n+1}/\CL_n$.

By induction $\taulin_\CM$ therefore gives rise to isomorphisms $\sigma^*(\CL_n/\CL_{n-m}) \isoto \CL_{n+1}/\CL_{n+1-m}$ for all $m\ge1$.
In the special case $m=rd_\infty$ this means that $\taulin_\CM$ thus induces isomorphisms $\sigma^*(\CL_n/z_\infty\CL_n) \isoto \CL_{n+1}/z_\infty\CL_{n+1}$. Since $\CL_n$ and $\CL_{n+1}$ are finitely generated projective modules over $\CE_{A_\infty,R}$, it follows that $\taulin_\CM$ induces isomorphisms $\sigma^*\CL_n \isoto \CL_{n+1}$.
 
Inverting $z_\infty$ this shows that the structure morphism of $\CM_{F,\infty}$ is an isomorphism; hence $\CM_{F,\infty}$ is an $F_\infty$-isocrystal over~$R$ (as we already knew from Propositions \ref{AMotRFpICP0Max} and \ref{AMotRFpICP0Zero}~\ref{AMotRFpICP0ZeroB}). On the other hand it follows that $\tau^{rd_\infty}$ induces isomorphisms $(\sigma^{rd_\infty})^*\CL_n \isoto \CL_{n+rd_\infty} = z_\infty^{-1}\CL_n$. Thus $\tau^{rd_\infty}\CL_n$ generates the lattice $z_\infty^{-1}\CL_n$, making $\CM_{F,\infty}$ pure of slope~$\minusoneoverr$, as desired.
\end{Proof}


\begin{Prop}\label{DrinFp0IC}
If $\wp$ is a maximal ideal, then $\CM$ gives rise to an $F_\wp$-isocrystal $\CM_{F,\wp}$ over $R[\epsilon_\wp^{-1}]$ that possesses a slope filtration and has slope~$0$ with multiplicity $r-h$ and slope~$\oneoverh$ with multiplicity~$h$.
\end{Prop}

\begin{Proof}
Without loss of generality we can assume that $\epsilon_\wp$ is a unit in~$R$, so that $R[\epsilon_\wp^{-1}]=R$. Consider the $\wp$-adic completion $\CL := A_\wp\complot_{A}\CM$ of~$\CM$. Then $\CM_{F,\wp} := F_\wp\otimes_{A_\wp} \CL$ is an $F_\wp$-isocrystal over $R$ by Proposition \ref{AMotRFpICP0Max}. Let $d_\wp$ denote the degree of the residue field $k_\wp$ over~$\BFq$. We shall construct an $\CE_{A_\wp,R}$-submodule $\CL_0 \subset \CL$
with the following properties:
\begin{enumerate}\StatementLabels
\item\label{P0FiltSub}
The submodule $\CL_0$ is finitely generated and projective.
\item\label{P0FiltQuot}
The quotient $\CL/\CL_0$ is finitely generated and projective.
\item\label{P0FiltSubGen}
The module $\CL_0$ is generated by $\tau\CL_0$.
\item\label{P0FiltQuotGen}
The module $\wp^{n}(\CL/\CL_0)$ is generated by $\tau^{h n d_\wp}(\CL/\CL_0)$.
\end{enumerate}
For this we note that by assumption $a_\wp$ is transcendental over~$\BFq$. We abbreviate $z:=a_\wp$ and observe that $A$ is a free module of rank $d := {\deg_A(z) = n d_\wp}$ over the subring $\BFq[z]$. Also $\CL$ then coincides with the $(z)$-adic completion of the $R[z]$-module~$\CM$. We will first construct $\CL_0\subset\CL$ as an $R\bbX{z}$-submodule.

\medskip
Since $\phi_z\in R[\tau]$ is elliptic of degree~$rd$, the $R[z]$-module $\CM = R[\tau]$ is free of rank $rd$ with the basis $m_1,\ldots,m_{rd}$, where $m_i := \tau^{rd-i}$ for all~$i$. By assumption we can write
$$\varphi_z = - \psi_{hd} \tau^{hd} - \psi_{hd+1} \tau^{hd+1} - \dotsc - \psi_{rd-1} \tau^{rd-1} + \psi_{rd} \tau^{rd}$$
with elements $\psi_i \in R$ such that $-\psi_{hd}=\epsilon_\wp$. Setting $s:=r-h$, the action of $\tau$ on $\CM$ with respect to the basis $m_1,\dotsc,m_{rd}$ is then given by the matrix
\UseTheoremCounterForNextEquation
\begin{equation}\label{MotAction}
\begin{tikzpicture}[baseline=(current  bounding  box.center)]
\matrix (mat) [matrix of math nodes, left delimiter=(, right delimiter=), outer sep=-5pt, column sep=.72em]
{
\frac{\psi_{rd-1}}{\psi_{rd}} & 1 \\
\phantom{0} & 0 \\
\phantom{0} & & & \phantom{0} & 1 \\
\frac{\psi_{hd}}{\psi_{rd}} & & & & 0 & 1 \\
0 & & & & & 0 & 1\\
\phantom{0}\\
0 & & & & & & & \phantom{0} & 1\\
\frac{z}{\psi_{rd}} & & & & & & & & 0\\
};
\draw[loosely dotted,line width=1pt] (-2.88,1.5) -- (-2.88,0.8);
\draw[loosely dotted,line width=1pt] (-2.88,-0.76) -- (-2.88,-1.21);
\draw[loosely dotted,line width=1pt] (-1.5,2.0) -- (-0.15,1.0);
\draw[loosely dotted,line width=1pt] (-1.5,1.4) -- (-0.15,0.4);
\draw[loosely dotted,line width=1pt] (1.9,-0.7) -- (3,-1.5);
\draw[loosely dotted,line width=1pt] (1.2,-0.7) -- (3,-2.0);
\draw[dotted,line width=.75pt] (0.4,-2.5) -- (0.4,2.5);
\draw[dotted,line width=.75pt] (-3.2,-0.14) -- (0.33,-0.14);
\draw[dotted,line width=.75pt] (0.53,-0.14) -- (3.2,-0.14);
\draw[snake=brace,line width=.66pt] (-4.2,0.12) -- (-4.2,2.4) node [left=2pt,midway] {\ensuremath{sd}} ;
\draw[snake=brace,line width=.66pt] (-4.2,-2.4) -- (-4.2,-0.27) node [left=2pt,midway] {\ensuremath{hd}} ;
\draw[snake=brace,line width=.66pt] (-3.24,3) -- (0.24,3) node [above=2pt,midway] {\ensuremath{sd}} ;
\draw[snake=brace,line width=.66pt] (0.55,3) -- (3.24,3) node [above=2pt,midway] {\ensuremath{hd}} ;
\end{tikzpicture}
\end{equation}
with all blank entries being zero. 
Let $\bar\CL_0$ be the free $R$-submodule of $\bar\CL := \CL/z\CL$ that is generated by the residue classes of $m_1,\dotsc,m_{sd}$.

\begin{Lem}\label{DrinFp0ICLem1}
The submodule $\bar\CL_0$ contains $\tau^{hd}\bar\CL$ and is generated by $\tau\bar\CL_0$ over $R$.
\end{Lem}

\begin{Proof}
Since the lower left block of the matrix \eqref{MotAction} is zero modulo $(z)$, we have ${\tau\bar\CL_0\subset\bar\CL_0}$. Also the upper left block has determinant ${\pm\frac{\psi_{hd}}{\psi_{rd}} = \mp \frac{\epsilon_\wp}{\psi_{rd}}}$,
which by assumption is a unit in~$R$. Together this implies that $\tau\bar\CL_0$ generates $\bar\CL_0$. 
On the other hand the lower right block of \eqref{MotAction} is a nilpotent $hd\times hd$-matrix; hence its $hd$-th power is zero. As the lower left block is zero modulo~$(z)$, it follows that $\tau^{hd}\bar\CL \subset \bar\CL_0$.
\end{Proof}

\medskip
The second part of Lemma \ref{DrinFp0ICLem1} implies that $\bar\CL_0$ is also generated by~$\tau^{hd}\bar\CL_0$ over~$R$. Applying Lemma \ref{LocICRGroucho} with $\rho=\tau^{hd}$ and $\bar\CL' = \bar\CL_0$, there therefore exists a unique $R\bbX{z}$-module direct summand $\CL_0\subset\CL$ such that ${(\CL_0+z\CL)/z\CL}=\bar\CL_0$ and $\tau^{hd}\CL_0\subset\CL_0$. 
Since $\bar\CL_0$ is free of rank $sd$ over~$R$, the module $\CL_0$ is free of rank $sd$ over~$R\bbX{z}$.

\begin{Lem}\label{DrinFp0ICLem2}
The submodule $\CL_0$ is generated by $\tau\CL_0$ over $R\bbX{z}$.
\end{Lem}

\begin{Proof}
Let $\CL_1$ be the $R\bbX{z}$-submodule of $\CL$ generated by $\tau\CL_0$. As the structure morphism $\taulin_\CL\colon \sigma^*\CL \to \CL$ is injective, this is again a free $R\bbX{z}$-module of rank~$sd$. Since $\bar\CL_0$ is generated by $\tau\bar\CL_0$, it also follows that ${(\CL_1+z\CL)/z\CL}=\bar\CL_0$. On the other hand the inclusion $\tau^{hd}\CL_0 \subset \CL_0$ implies that $\tau^{hd}\CL_1\subset\CL_1$. The uniqueness part of Lemma \ref{LocICRGroucho} thus shows that $\CL_1 = \CL_0$.
\end{Proof}

\medskip
Next, the quotient $\CL/\CL_0$ is free of rank $hd$ over $R\bbX{z}$ and inherits an action of~$\tau^{hd}$.

\begin{Lem}\label{DrinFp0ICLem3}
The $R\bbX{z}$-module $z(\CL/\CL_0)$ is generated by $\tau^{hd}(\CL/\CL_0)$.
\end{Lem}

\begin{Proof}
The matrix \eqref{MotAction} has determinant $\pm\frac{z}{\psi_{rd}}$, which is $z$ times a unit in~$R$. Thus the morphism induced by $\tau$ on the highest exterior power of $\CL$
$$\det\taulin_\CL\colon \sigma^*(\det\CL)\ \longto\ \det\CL$$ 
is injective with the image $z(\det\CL)$. It follows that the morphism
$$\det(\tau^{hd}_\CL)^{\textup{lin}}\colon (\sigma^{hd})^*(\det\CL)\ \longto\ \det\CL$$
is injective with the image $z^{hd}(\det\CL)$. Since the morphism
$$\det(\tau^{hd}_{\CL_0})^{\textup{lin}}\colon (\sigma^{hd})^*(\det\CL_0)\ \longisoto\ \det\CL_0$$
is an isomorphism, we deduce that the morphism
$$\det(\tau^{hd}_{\CL/\CL_0})^{\textup{lin}}\colon (\sigma^{hd})^*\bigl(\det(\CL/\CL_0)\bigr)\ \longto\ \det(\CL/\CL_0)$$
is injective with the image $z^{hd}(\det(\CL/\CL_0))$. Now the fact that $\tau^{hd}\bar\CL\subset\bar\CL_0$ implies that $\tau^{hd}(\CL/\CL_0)\subset z(\CL/\CL_0)$. Thus we can write $(\tau^{hd}_{\CL/\CL_0})^{\textup{lin}} = z{\cdot}\rho$ for some linear map 
$$\rho\colon (\sigma^{hd})^*(\CL/\CL_0)\ \longto\ \CL/\CL_0.$$
The stated property of $\det(\tau^{hd}_{\CL/\CL_0})^{\textup{lin}}$ then implies that 
$\det\rho$ is surjective. Thus $\rho$ is surjective, and the lemma follows.
\end{Proof}

\medskip
Since $Az=\wp^n$, the facts collected so far mean that $\CL_0$ satisfies the analogues of properties \ref{P0FiltSub} through \ref{P0FiltQuotGen} with respect to the ring $R\bbX{z}$ in place of $\CE_{A_\wp,R}$. Next we show:

\begin{Lem}\label{DrinFp0ICLem4}
The submodule $\CL_0$ is an $\CE_{A_\wp,R}$-submodule of~$\CL$.
\end{Lem}

\begin{Proof}
Lemma \ref{DrinFp0ICLem1} implies that $\bar\CL_0$ is generated by $\tau^{hd}\bar\CL_0$ 
and hence also by $\tau^{hd}\bar\CL$ over~$R$. Next recall that $\CL$ and $\bar\CL$ are $A_\wp$-modules. For any unit $a\in A_\wp^\times$, the equality $a\bar\CL=\bar\CL$ thus implies that $a\bar\CL_0=\bar\CL_0$. It follows that $a\CL_0\subset\CL$ is again a free $R\bbX{z}$-submodule of rank $sd$ such that ${(a\CL_0+z\CL)/z\CL}=\bar\CL_0$ and $\tau^{hd}a\CL_0\subset a\CL_0$. By the uniqueness part of Lemma \ref{LocICRGroucho} we thus have $a\CL_0 = \CL_0$. 

Now observe that since $A_\wp$ is a local ring, for any non-unit $a\in A_\wp$ the element $1+a$ is a unit in $A_\wp^\times$. By the argument above we thus have $(1+a)\CL_0=\CL_0$. Together we therefore conclude that $A_\wp\CL_0\subset\CL_0$. 

Finally, since $A_\wp$ generates $\CE_{A_\wp,R} = A_\wp \complot_{\BFq}\kern-1ptR$ over $R\bbX{z}$, the lemma follows.
\end{Proof}

\begin{Lem}\label{DrinFp0ICLem5}
The submodule $\CL_0$ satisfies \ref{P0FiltSub} through \ref{P0FiltQuotGen} with respect to $\CE_{A_\wp,R}$. 
\end{Lem}

\begin{Proof}
Since $\CL_0$ is a free module of finite rank over~$R\bbX{z}$ and an $\CE_{A_\wp,R}$-module by Lemma \ref{DrinFp0ICLem4}, by Lemma \ref{LatticeFreeLem} it is therefore finitely generated projective over $\CE_{A_\wp,R}$, proving~\ref{P0FiltSub}.
Likewise, since $\CL/\CL_0$ is a free module of finite rank over~$R\bbX{z}$ and an $\CE_{A_\wp,R}$-module by Lemma \ref{DrinFp0ICLem4}, by Lemma \ref{LatticeFreeLem} it is therefore finitely generated projective over $\CE_{A_\wp,R}$, proving~\ref{P0FiltQuot}. Finally, since \ref{P0FiltSubGen} and \ref{P0FiltQuotGen} already hold with respect to the subring $R\bbX{z}$, they hold a fortiori with respect to $\CE_{A_\wp,R}$.
\end{Proof}

\medskip
The properties \ref{P0FiltSub} and \ref{P0FiltSubGen} imply that $\CM_{F,\wp,0} := F_\wp \otimes_{A_\wp}\kern-1pt\CL_0$ is an $F_\wp$-sub-isocrystal of $\CM_{F,\wp}$, and the lattice $\CL_0$ shows that it is pure of slope~$0$. Likewise, the properties \ref{P0FiltQuot} and \ref{P0FiltQuotGen} imply that $\CM_{F,\wp}/\CM_{F,\wp,0}$ is an $F_\wp$-isocrystal, and the lattice $\CL/\CL_0$ shows that it is pure of slope~$\oneoverh$. Thus $\CM_{F,\wp}$ carries a slope filtration over $R$ with the desired slopes.

Finally, the fact that $\CL_0$ is free of rank $sd$ over~$R\bbX{z}$ implies that $\CM_{F,\wp,0} = \CL_0[z^{-1}]$ is free of rank $sd$ over $R\ppX{z}$. On the other hand $\CM_{F,\wp,0}$ is locally free of constant rank over $\CE_{F_\wp,R} = F_\wp\complot_{\BFq}\kern-1ptR$ by Proposition \ref{LICRGoodRedRank}. Since $F_\wp$ has dimension $d$ over $\BFq\ppX{z}$, it follows that $\CM_{F,\wp,0}$ is locally free of rank $s=r-h$ over $\CE_{F_\wp,R}$. The same argument with $\CL/\CL_0$ shows that $\CM_{F,\wp}/\CM_{F,\wp,0}$ is locally free of rank $h$ over $\CE_{F_\wp,R}$. The proof of the proposition is thus complete.
\end{Proof}

\medskip
Now we assume that $\phi$ is defined over a field~$K$. For any place $\Fp\not=\infty$ of $F$ the rational Tate module $V_\Fp(\phi)$ is an $\ICF$-vector space of dimension $r$ if $\Fp\not=\wp$, respectively of dimension $r-h$ if $\Fp=\wp$. In either case it comes with a natural continuous representation of~$\GK$. 
Let $M$ denote the $A$-motive over~$K$ that is associated to~$\phi$, and by $M_F$ the $F$-isocrystal associated to~$M$. Recall that $\Omega_A$ is a projective $A$-module of rank~$1$.

\begin{Prop}\label{DrinFpTate}
For any place $\Fp\not=\infty$ of~$F$, 
there is a natural $\GK$-equivariant isomorphism of $F_\Fp$-vector spaces
$$V_\Fp(\phi)\ \cong\ T_0(M_{F,\Fp})^\vee \otimes_A\Omega_A.$$
\end{Prop}

\begin{Proof}
Anderson \cite[Prop.\,1.8.3]{AndersonT} proved this when $A=\BFq[t]$ and $\wp$ is the zero ideal, and we will adapt his proof to the general case.

First observe that both sides are defined using the same algebraic closure $\Kalg$ of~$K$. Any natural isomorphism will thus automatically be $\GK$-equivariant. After replacing $K$ by $\Kalg$ we may therefore assume that $K$ is algebraically closed. 

We view $K$ as a left module over $K[\tau]$ by the action $(\sum_i\psi_i\tau^i,x) \mapsto \sum_{i\ge0}\psi_i x^{q^i}$. Consider any non-zero element $a\in A$ and abbreviate $\phi[a] := \{\xi\in K\mid \phi_a(\xi)=0\}.$ 

\begin{Lem}\label{DrinFpTatePairing1}
There is a natural perfect pairing of finite dimensional $\BFq$-vector spaces
$$B\colon\ (M/aM)^\tau \times \phi[a] \longto \BFq,\ \ ([m],\xi)\mapsto m(\xi).$$
\end{Lem}

\begin{Proof}
The definition of the $A\otimes_{\BFq}\kern-1pt K[\tau]$-module structure on $M=K[\tau]$ implies that $M/aM = K[\tau]/K[\tau]\phi_a$. For each $\xi\in\phi[a]$, the value $m(\xi)\in K$ for $m\in M$ thus depends only on the residue class $[m] \in M/aM$. Moreover, since $\tau$ acts on $K$ through the $q$-Frobenius map~$\sigma$, for any $\tau$-invariant element $[m] \in (M/aM)^\tau$ we have $m(x)\in K^\sigma=\BFq$. The map is therefore well-defined. Clearly it is $\BFq$-bilinear. Also, for any $x\in\phi[a]$ with $m(\xi)=0$ for all $[m] \in (M/aM)^\tau$, taking $m=1$ we deduce that $\xi=0$. Thus the pairing is non-degenerate on the right hand side.

Next, the assumption $a\not=0$ implies that $\phi_a$ is non-zero as well. Since $K$ is perfect, we can therefore write $\phi_a$ uniquely in the form $\tau^n\psi$ for an integer $n\ge0$ and an element $\psi\in K[\tau]$ with non-zero constant term. By definition $\phi[a]$ is then simply the zero set of~$\psi$. Let $d\ge0$ be the degree of $\psi$ with respect to~$\tau$. Since $\psi(X)$ is a separable polynomial in~$X$, it follows that $\phi[a]$ is an $\BFq$-vector space of finite dimension~$d$. As the pairing is already non-degenerate on the right hand side, to prove that it is perfect it thus suffices to show that $\dim_{\BFq}(M/aM)^\tau$ is equal to $d$ as well.

For this consider the $K$-linear map
$$K[\tau] \longto K[\tau]/K[\tau]\psi \times K[\tau]/K[\tau]\tau^n,\ \ m\mapsto ([m],[m]).$$
Its kernel is the intersection $K[\tau]\psi \cap K[\tau]\tau^n$. Since $K$ is a field and the constant term of $\psi$ is non-zero, for any $\theta\in K[\tau]$ we have $\theta\psi\in K[\tau]\tau^n$ if and only if $\theta$ itself lies in $K[\tau]\tau^n$. Thus the kernel is $K[\tau]\tau^n\psi = K[\tau]\phi_a$, and so the map induces an injective homomorphism
$$M/aM\ =\ K[\tau]/K[\tau]\phi_a\ \longinto\ K[\tau]/K[\tau]\psi \times K[\tau]/K[\tau]\tau^n.$$
Here the two $K$-vector spaces on the right hand side have respective dimensions $d$ and~$n$, and since $\phi_a=\tau^n\psi$ the one on the left has dimension $n+d$. Thus the map is already an isomorphism. Now the fact that $\tau$ is nilpotent on $K[\tau]/K[\tau]\tau^n$ implies that its subspace of $\tau$-invariants is zero. On the other hand $\tau$ is surjective on $K[\tau]/K[\tau]\psi$; hence by Lang's theorem \cite[Exp.\,XXII Prop.\,1.1]{SGA7}, its subspace of $\tau$-invariants has dimension~$d$ over~$\BFq$. 
Together this proves that $(M/aM)^\tau$ has dimension $d$ over~$\BFq$, as desired.
\end{Proof}

\medskip
Next, the pairing $B$ from Lemma \ref{DrinFpTatePairing1} satisfies $B([bm],\xi) = B([m\phi_b],\xi) = m\phi_b(\xi) = B([m],\phi_b(\xi))$ for all $b\in A$. Endowing $\Hom_{\BFq}(A/(a),\BFq)$ with the structure of an $A$-module by $(bh)([b']) = h(b'b)$ for all $b,b'\in A$, we deduce a natural perfect $A$-bilinear pairing
\UseTheoremCounterForNextEquation
\begin{equation}\label{DrinFpTatePairing2}
(M/aM)^\tau \times \phi[a] \longto \Hom_{\BFq}(A/(a),\BFq),
\ \ ([m],\xi)\mapsto \bigl(b\mapsto m(\phi_b(\xi))\bigr).
\end{equation}
Now recall that there is also a natural perfect pairing of $\BFq$-vector spaces
$$A/(a) \times a^{-1}\Omega_A/\Omega_A \longto \BFq,
\ \ ([b],[\omega]) \mapsto \trace_{k_\infty/\BFq}(\mathop{\rm res}\nolimits_\infty(b\omega)).$$
This corresponds to a natural isomorphism of $A$-modules
$$\Hom_{\BFq}(A/(a),\BFq)\ \cong\ a^{-1}\Omega_A/\Omega_A\ \cong\ (a^{-1}A/A) \otimes_A \Omega_A.$$
The pairing \eqref{DrinFpTatePairing2} thus yields a natural perfect $A$-bilinear pairing
$$(M/aM)^\tau \times \phi[a] \longto (a^{-1}A/A) \otimes_A \Omega_A.$$
Thus in turn corresponds to a natural isomorphism of $A$-modules
\UseTheoremCounterForNextEquation
\begin{equation}\label{DrinFpTatePairing4}
\Hom_A\bigl(a^{-1}A/A,\phi[a]\bigr)\ \cong\ \Hom_A\bigl((M/aM)^\tau,A/(a)\bigr) \otimes_A \Omega_A.
\end{equation}

All this holds for an arbitrary non-zero element $a\in A$, but now we specialize it to the case $a=a_\Fp^i$ for any integer $i\ge1$. Since $(a_\Fp)=\Fp^n$ for some $n\ge1$, the isomorphism \eqref{DrinFpTatePairing4} yields a natural isomorphism of $A$-modules
$$\Hom_A\bigl(\Fp^{-ni}/A,\phi[\Fp^{ni}]\bigr)\ \cong\ \Hom_A\bigl((M/\Fp^{ni}M)^\tau,A/\Fp^{ni}\bigr) \otimes_A \Omega_A.$$
Also, one checks that these isomorphisms are compatible with the natural transition maps from $\Fp^{ni}$ to $\Fp^{n(i+1)}$. Taking the inverse limit over $i$ we thus obtain a natural isomorphism of $A_\Fp$-modules
$$T_\Fp(\phi)\ :=\ \textstyle \Hom_A\bigl(F_\Fp/A_\Fp,\bigcup_i\phi[\Fp^{ni}]\bigr)\ \cong\ \Hom_{A_\Fp}\bigl((M\otimes_AA_\Fp)^\tau,A_\Fp\bigr) \otimes_A \Omega_A.$$
By construction that is independent of the choice of~$a_\Fp$. Tensoring it with $F_\Fp$ over $A_\Fp$ finally yields a natural isomorphism of $F_\Fp$-vector spaces
$$V_\Fp(\phi)\ :=\ T_\Fp(\phi)\otimes_{A_\Fp}F_\Fp\ \cong\ \Hom_{F_\Fp}\bigl((M\otimes_AF_\Fp)^\tau,F_\Fp\bigr) \otimes_A \Omega_A
\ =\ T_0(M_{F,\Fp})^\vee \otimes_A\Omega_A.$$
This finishes the proof of Proposition \ref{DrinFpTate}.
\end{Proof}

\medskip
In view of Proposition \ref{DrinFpTate}, we define the following analogues of $V_\Fp(\phi)$. 

\begin{Cons}\label{DrinVinftyDef}\rm
For $\Fp=\infty$ we set
$$V_{\infty,\oneoverr}(\phi)\ :=\ T_{\minusoneoverr}(M_{F,\infty})^\vee \otimes_A\Omega_A.$$
Recall from Proposition \ref{DrinFinftyIC} that $\CM_{F,\infty}$ is an $F_\infty$-isocrystal of rank $r$ that is pure of slope~$\minusoneoverr$. Thus \eqref{TWDef} implies that $T_{\minusoneoverr}(M_{F,\infty})$ is a one-dimensional right vector space over a central division algebra $D_\infty := \End(\SAa{\kalg}{\minusoneoverr})$ with Hasse invariant $[\oneoverr]$ over~$F_\infty$. Its dual is therefore naturally a one-dimensional left vector space over~$D_\infty$. Since $\Omega_A$ is an invertible $A$-module, the same thus also holds for $V_{\infty,\oneoverr}(\phi)$.

For any basis element $v\in V_{\infty,\oneoverr}(\phi)$ the natural left action of $\gamma\in\WK$ is then given by the formula ${}^\gamma v = \rho(\gamma)v$ for a unique continuous homomorphism $\rho\colon \WK \to (D_\infty^\opp)^\times$. By \cite[Thm.\,9.2.2]{MornevT}, up to identification of the division algebra $D_\infty^\opp$ (and in particular up to conjugation) this homomorphism coincides with that introduced by J.-K. Yu in \cite{Yu2003}. 
\end{Cons}


\begin{Cons}\label{DrinVp0+Def}\rm
When $\wp$ is a maximal ideal and $\phi$ has height~$h$ over~$K$, we set
$$V_{\wp,\minusoneoverh}(\phi)\ :=\ T_{\oneoverh}(M_{F,\wp})^\vee \otimes_A\Omega_A.$$
Recall from Proposition \ref{DrinFp0IC} that $\CM_{F,\wp}$ is an $F_\wp$-isocrystal for which $\oneoverh$ is a slope with multiplicity~$h$. Thus \eqref{TWDef} implies that $T_{\oneoverh}(M_{F,\wp})$ is a one-dimensional right vector space over a central division algebra $D_\wp := \End(\SAa{\kalg}{\oneoverh})$ with Hasse invariant $[\minusoneoverh]$ over~$F_\wp$. Its dual is therefore naturally a one-dimensional left vector space over~$D_\wp$. Since $\Omega_A$ is an invertible $A$-module, the same thus also holds for $V_{\wp,\minusoneoverh}(\phi)$.
For any basis element $v\in V_{\wp,\minusoneoverh}(\phi)$ the natural left action of $\gamma\in\WK$ is then given by the formula ${}^\gamma v = \rho(\gamma)v$ for a unique continuous homomorphism $\rho\colon \WK \to (D_\wp^\opp)^\times$. 
\end{Cons}

Also, for every place $\Fp$ of~$F$ we set
\UseTheoremCounterForNextEquation
\begin{equation}\label{DrinUFpDef}
U_\Fp(\phi)\ :=\ U(M_{F,\Fp})^\vee \otimes_A\Omega_A.
\end{equation}

\begin{Prop}\label{DrinUFpIsom}
This is an $\ICF$-isocrystal of rank $r$ over $\kalg$ with a continuous left action of~$\WK$. Moreover there are natural $\WK$-equivariant isomorphisms of $\ICF$-isocrystals over~$\kalg$
$$U_\Fp(\phi)\ \cong\ 
\left\{\begin{array}{ll}
\unit_{\kalg} \otimes_{F_\Fp} V_\Fp(\phi)
& \hbox{if $\Fp\not=\wp,\infty$,} \\[5pt]
\SAa{\kalg}{\minusoneoverr}^{\mkern+2mu\vee} \otimes_{D_\infty} V_{\infty,\oneoverr}(\phi)
& \hbox{if $\Fp=\infty$,} \\[5pt]
\bigl[ \unit_{\kalg} \otimes_{F_\wp} V_\wp(\phi) \bigr]
\oplus
\bigl[ \SAa{\kalg}{\oneoverh}^{\mkern+2mu\vee} \otimes_{D_\wp} V_{\wp,\minusoneoverh}(\phi) \bigr]
& \hbox{if $\Fp=\wp$.}
\end{array}\right.$$
\end{Prop}

\begin{Proof}
Combining Proposition \ref{WMBigWeilRepGrad} with the formula \eqref{WMSmallWeilRepIsom}, the definition of $U_\Fp(\phi)$ yields a natural $\WK$-equivariant isomorphism
$$U_\Fp(\phi)\ \cong\ 
\bigoplus_{\alpha\in\BQ} \SA{\kalg}^{\mkern+2mu\vee} \centrot_{\raisebox{-1pt}{$\scriptstyle\, D_\alpha$}} T_\alpha(M_{F,\Fp})^\vee  \centrot_{\raisebox{-1pt}{$\scriptstyle A\;$}}\,\Omega_A.$$
In the case $\Fp\not=\wp,\infty$ the desired formula follows from Propositions \ref{DrinFpIC} and \ref{DrinFpTate}. In the case $\Fp=\infty$ it follows from Proposition \ref{DrinFinftyIC} and Construction~\ref{DrinVinftyDef}. In the case $\Fp=\wp$ it follows from Propositions \ref{DrinFp0IC} and \ref{DrinFpTate} and Construction~\ref{DrinVp0+Def}.
\end{Proof}

\medskip
Finally, observe that $\phi$ comes from a Drinfeld $A$-module defined over some finitely generated normal $\BFq$-subalgebra $R\subset K$ that contains the image $c^*(A)$, for which we can also assume that $\Quot(R)=K$.
It is well-known (e.g. by Goss \cite[Thm.\,3.2.3 (b)]{Goss1992}) that the Tate modules $T_\Fp(\phi)$ form a compatible system of Galois representations for $\Fp\not=\wp,\infty$. We can now generalize this to a compatible system of Weil group representations for all places of~$F$. For $\Fp=\infty$ this reproduces the result of J.-K. Yu \cite[Thm.\,3.4]{Yu2003}.

\begin{Thm}\label{CompatibleSystemDM}
Consider any place $\Fp$ of $F$ and any closed point $x \in \Spec R$. If $\Fp\not=\wp$ assume that $c(x)\ne \Fp$, and if $\Fp=\wp$ assume that $\epsilon_\wp(x)\not=0$. Then the inertia group $I_x$ acts trivially on $U_\Fp(\phi)$, and the arithmetic Frobenius element $\Frob_x$ has characteristic polynomial $\charact_{\CM_{x,F}}$. In particular this polynomial has coefficients in $A$ and is independent of~$\Fp$.
\end{Thm}

\begin{Proof}
Since the action of $\WK$ on $U(M_{F,\Fp})^\vee$ is contragredient to the action on $U(M_{F,\Fp})$, the construction \eqref{DrinUFpDef} implies that the characteristic polynomial of $\Frob_x$ on $U_\Fp(\phi)$ is equal to the characteristic polynomial of $\Frob_x^{-1}$ on $U(M_{F,\Fp})$. For $\Fp\not=\wp,\infty$ the statements are therefore direct consequences of Theorem \ref{CompatibleSystemUAMot}. 
For $\Fp=\infty$ Proposition \ref{DrinFinftyIC} says that Proposition \ref{GlobICFpGoodRed} holds with $\epsilon_\infty=1$, and for $\Fp=\wp$ Proposition \ref{DrinFp0IC} says that Proposition \ref{GlobICFpGoodRed} holds with the present choice of $\epsilon_\wp$. The remaining cases thus result from Theorem \ref{CompatibleSystemU}.
\end{Proof}

%
%
%
%
%
%
%


\end{document}